\documentclass[twoside,10pt,halfparskip,headsepline,fleqn]{scrartcl} 

\usepackage[utf8]{inputenc}          
\usepackage[T1]{fontenc}             
\usepackage{amsmath}
\usepackage{amsfonts}                
\usepackage{amssymb}                 
\usepackage{mathrsfs}				
\usepackage{textcomp}                
\usepackage{typearea}                
\usepackage{scrlayer-scrpage}        
\usepackage{array}                   
\usepackage{enumitem}               
\usepackage{titlesec}                
\usepackage[amsmath, thmmarks]{ntheorem} 
\usepackage{mathtools}				
\usepackage{mdframed}				
\usepackage{xcolor}					
\usepackage{hyperref}				
\usepackage{cite}

\typearea[10mm]{12} 				
\titleformat{\section}{\normalfont\normalcolor\Large\bfseries}{\thesection}{.66em}{}
\titleformat{\subsection}{\normalfont\normalcolor\large\bfseries}{\thesubsection}{.66em}{}

\clearpairofpagestyles 			
\ohead{\normalfont \headmark} 
\automark[subsection]{section}	

\makeatletter
\newtheoremstyle{StandardMitBezInKlammern}%
 {\item[\rlap{\vbox{\hbox{\hskip\labelsep \theorem@headerfont
  ##1\ ##2\theorem@separator}\hbox{\strut}}}]}%
 {\item[\rlap{\vbox{\hbox{\hskip\labelsep \theorem@headerfont
  ##1\ ##2 (##3)\theorem@separator}\hbox{\strut}}}]}%
\makeatother
\theoremstyle{StandardMitBezInKlammern}
\theorembodyfont{\rmfamily}
\theoremsymbol{\ensuremath{\diamond}}
\definecolor{vlgray}{gray}{0.85}
\newtheorem{defsatzusw}{}[subsection]
\newmdtheoremenv[backgroundcolor=vlgray,linecolor=black,linewidth=1.5pt,innertopmargin=10pt]{definition}[defsatzusw]{Definition}
\newmdtheoremenv[backgroundcolor=vlgray,linecolor=black,linewidth=1.5pt,innertopmargin=10pt]{lemma}[defsatzusw]{Lemma}
\newmdtheoremenv[backgroundcolor=vlgray,linecolor=black,linewidth=1.5pt,innertopmargin=10pt]{corollary}[defsatzusw]{Corollary}
\newmdtheoremenv[backgroundcolor=vlgray,linecolor=black,linewidth=1.5pt,innertopmargin=10pt]{proposition}[defsatzusw]{Proposition}
\newmdtheoremenv[backgroundcolor=vlgray,linecolor=black,linewidth=1.5pt,innertopmargin=10pt]{theorem}[defsatzusw]{Theorem}
\newmdtheoremenv[backgroundcolor=vlgray,linecolor=black,linewidth=1.5pt,innertopmargin=10pt]{definition/lemma}[defsatzusw]{Definition/Lemma}
\newtheorem{remark}[defsatzusw]{Remark}

\theoremstyle{nonumberplain}
\theoremsymbol{\ensuremath{\square}}
\newtheorem{proof}{Proof}

\numberwithin{equation}{subsection} 	

\newenvironment{mathex}[1][LLLLLLLLLLLL]{
  \[%
  \newcolumntype{L}{>{\displaystyle\setlength{\arraycolsep}{4pt}}l}%
  \newcolumntype{C}{>{\displaystyle\setlength{\arraycolsep}{4pt}}c}%
  \newcolumntype{R}{>{\displaystyle\setlength{\arraycolsep}{4pt}}r}%
  \setlength{\arraycolsep}{1.5pt}%
  \begin{array}{>{\vspace*{1ex}}#1}%
}{
  \end{array}%
  \vspace*{-1ex}%
  \]%
  \ignorespacesafterend%
}

\newcommand{\N}{\mathbb{N}}
\newcommand{\Nz}{\mathbb{N}_0}
\newcommand{\Z}{\mathbb{Z}}

\newcommand{\R}{\mathbb{R}}
\newcommand{\Rn}{{\R^n}}
\newcommand{\Rnp}{{\R^{n+1}}}
\newcommand{\Rnpp}{{\R^{n+1}_+}}

\newcommand{\C}{\mathcal{C}}
\newcommand{\CicRn}{{\mathcal{C}_c^\infty(\Rn)}}
\newcommand{\CicRnRn}{{\mathcal{C}_c^\infty(\Rn,\Rn)}}

\newcommand{\const}{C}

\newcommand{\Lor}[3]{{L^{(#1,#2)}(#3)}}		

\newcommand{\Rj}{{\mathcal{R}_j}}			
\newcommand{\Hil}{{\mathcal{H}}}				
\newcommand{\fL}[1]{{(-\Delta)^\frac{#1}{2}}}	
\newcommand{\Fourier}{{\mathcal{F}}}			

\begin{document}

\thispagestyle{empty}
\begin{center}

\Huge\textbf{Estimates for commutators of\\ fractional differential operators\\ via harmonic extension\\}
\begin{large}
\vspace{2.5cm}
		by\\
	Jonas Ingmanns\\
\vspace{2.5cm}
	Master Thesis in Mathematics\\
\vspace{2.5cm}
		presented to the\\
	Faculty of Mathematics, Computer Science and Natural Sciences\\
	of the Rheinisch-Westfälische Technische Hochschule Aachen\\
\vspace{2cm}
	in November 2020\\
\vspace{1cm}
	Prof. Dr. Heiko von der Mosel
\end{large}

\end{center}
\newpage

\thispagestyle{empty}
\quad\newpage

\tableofcontents
\thispagestyle{plain}
\newpage

\thispagestyle{empty}
\quad\newpage

\addcontentsline{toc}{section}{Notation}
\section*{Notation}
\begin{mathex}
 	\lesssim					&\quad	\text{We write $\mathcal{I}_1 \lesssim \mathcal{I}_2$ when there exists a constant 
 									$C>0$ such that $\mathcal{I}_1 \leq C\,\mathcal{I}_2$,}\\
 							&\quad	\text{where $C$ may depend on the parameters such as $n,s,p,q,\ldots$ but not on the}\\
 							&\quad	\text{variable quantities such as functions $f,g,\phi$ or variables $x,y,t$.}\\
 	\const					&\quad	\text{denotes such a constant as described just above}\\
 	\Rnpp					&\quad	\text{the upper half-space, i.e. }\lbrace(x,t)\in\Rnp\,;\,t>0\rbrace\\
 	\quad\\

 	\Fourier	(f),\widehat{f}	&\quad	\text{$n$-dimensional Fourier transform,}
 										\quad \Fourier(f)(\xi)= \widehat{f}(\xi) = \int_{\Rn}f(x)e^{-2\pi i x\cdot\xi}\,dx\\
 	\mathcal{M}(f)			&\quad	\text{maximal function, see (\ref{eq:maxfunc - definition})}\\
 	P^1_t(f)					&\quad	\text{classical harmonic extension, see Definition \ref{def:P1t Poisson operator}}\\
 	P^s_t(f)					&\quad	s\text{-harmonic extension, see Definition \ref{def:Pst generalized Poisson operator}}\\
 	\Rj(f),~\Rj[f]			&\quad	\text{Riesz transform, see Definition \ref{def:Rt}}\\
 	\Hil(f),	~\Hil[f]			&\quad	\text{Hilbert transform, the 1-dimensional Riesz transform}\\
 	\fL{s}(f)				&\quad	\text{fractional Laplacian, see Definition \ref{def:fL fL}}\\
 	I^s(f)					&\quad	\text{Riesz potential, see Definition \ref{def:Rp}}\\
 	\dot{J}^\sigma(f)		&\quad	\text{lifting operator for homogeneous Triebel-Lizorkin and Besov-Lipschitz spaces}\\
 	J^\sigma(f)				&\quad	\text{lifting operator for inhomogeneous Triebel-Lizorkin and Besov-Lipschitz spaces}\\
 	\quad\\
 	
 	\mathcal{P}(\Rn)			&\quad	\text{space of polynomials over }\Rn\\
 	\C^\infty(\Omega)		&\quad	\text{space of infinitely differentiable functions}\\
 	\C^\infty_c(\Omega)		&\quad	\text{space of compactly supported }\C^\infty\text{-functions}\\
 	\mathcal{S}(\Rn)			&\quad	\text{Schwartz space, rapidly decreasing $C^\infty$-functions}\\
 	\mathcal{S}'(\Rn)		&\quad	\text{space of tempered distributions, i.e. the topological dual of }\mathcal{S}(\Rn)\\
 	L^p(\Omega)				&\quad	\text{Lebesgue space}\\
 	L^{(p,q)}(\Omega)		&\quad	\text{Lorentz space, see Subsection \ref{subsec:tt - lorentz spaces}}\\
 	\dot{B}^\alpha_{p,q}(\Rn)
 			&\quad	\text{homogeneous Besov-Lipschitz space, see Subsection \ref{subsec:tt - triebel lizorkin besov}}\\
 	\dot{F}^\alpha_{p,q}(\Rn)
 			&\quad	\text{homogeneous Triebel-Lizorkin space, see Subsection \ref{subsec:tt - triebel lizorkin besov}}\\
 	B^\alpha_{p,q}(\Rn)		
 			&\quad	\text{inhomogeneous Besov-Lipschitz space, see Subsection \ref{subsec:tt - triebel lizorkin besov}}\\
 	F^\alpha_{p,q}(\Rn)		
 			&\quad	\text{inhomogeneous Triebel-Lizorkin space, see Subsection \ref{subsec:tt - triebel lizorkin besov}}\\
 	W^m_p(\Rn)				&\quad	m\in\N,~\text{Sobolev space}\\
 	W^s_p(\Rn)				&\quad	s\notin\N,~\text{Slobodeckij space, see (\ref{eq:tt def Slobodeckij})}\\
 	BMO(\Rn)					&\quad	\text{space of functions of bounded mean oscillation, see (\ref{eq:tt def BMO})}\\
 	\mathcal{H}_p			&\quad	\text{Hardy space, for $p=1$ its topological dual is }BMO\\    
\end{mathex}

\newpage

\thispagestyle{empty}
\quad\newpage

%
%

\section{Introduction}\label{sec:intro}
The goal of this work is to present an interesting method for proving a variety of commutator estimates involving Riesz transforms, fractional Laplacians and Riesz potentials. 
This thesis is based in large parts on the recent work \cite{LeSchi20} of Lenzmann and Schikorra, in which they proposed this method.
What makes the method interesting, is its accessibility. 
Except for some elementary transformations, the only tools needed are harmonic extensions from $\Rn$ to $\Rnpp$ and integration by parts.
That is, as long as we accept some general estimates from $\Rnpp$ back to $\Rn$ as blackboxes.
Next to the theoretical background for these blackbox estimates, we also provide a number of applications for this method.
Many of the commutator estimates we show are already well known, but originally they were either proven via an individual approach or required a lot more effort.

Before describing the outline of this thesis, we give a short sketch on how the method works. 
Usually, we estimate the commutators either in $L^p$ or in the Hardy space $\mathcal{H}^1$.
These space norm estimates are typically obtained via the respective dual space characterization. 
Therefore, we start with some integral term over $\Rn$.
We run through three steps.
First, we interpret $\Rn$ as boundary of the higher dimensional half-space $\Rnpp$. 
We can then extend every occurring function to $\Rnpp$ by solving the Poisson equation with the respective function as boundary value, obtaining the harmonic extension.
Second, via a suitable integration by parts we obtain an integral over $\Rnpp$.
Third, via some elementary transformations and estimates we reach a suitable form, to which we can apply one of the blackbox estimates. 
During these transformations, we will usually observe ``cancellation effects'' due to the commutator's structure, which are the reason why we obtain better estimates than if we would just estimate each component of the commutators separately.

We will further illustrate this approach with a first application in this section, showing a Hardy-space estimate for the Jacobian $\det(\nabla u)$ as well as an intermediate estimate.
The remaining sections can be divided into two parts. 
In Sections \ref{sec:classical harmonic extension} and \ref{sec:s harmonic extension}, we show the application of the method, proving various commutator estimates.
Sections \ref{sec:operators} and \ref{sec:tt} provide the theoretical background. 

In Section \ref{sec:classical harmonic extension}, we first recall the classical harmonic extension and collect a variety of facts for later use. 
We then prove estimates for three different commutators involving classical differential operators or the half-Laplacian $\fL{1}$.
\begin{itemize}[noitemsep]
\item Subsection \ref{sec:divcurl by CLMS}: The div-curl estimate by Coifman, Lions, Meyer, and Semmes. 
\item Subsection \ref{sec:Rt comms}: The Coifman-Rochberg-Weiss commutator estimate for Riesz transforms.
\item Subsection \ref{sec:1d double comm}: An $L^1$-estimate for a double-commutator complementing Subsection \ref{sec:Rt comms}.
\end{itemize}
In Section \ref{sec:s harmonic extension}, we introduce the $s$-harmonic extension, also known as $s$-Poisson extension. 
This is a generalization of the classical harmonic extension, which allows us to treat commutators involving fractional Laplacians $\fL{s}$. 
After collecting some results for this generalized harmonic extension in Subsection \ref{sec:Pst}, we show four more applications for the method.
\begin{itemize}[noitemsep]
\item Subsection \ref{sec:fL comms}: Estimates for commutators of multiplication and fractional Laplacians.
\item Subsection \ref{sec:Rp comms}: The Chanilo commutator estimate for Riesz potentials of order $<1$.
\item Subsection \ref{sec:fLr 3-t-comm}: Fractional Leibniz rules.
\item Subsection \ref{sec:fLr 3-t-comm hardy}: The Da Lio-Rivi\`ere three-term commutator estimate.
\end{itemize}
For these two sections, we closely follow \cite{LeSchi20}. 
However, we choose a different structure and go into more detail regarding the basic properties of the extensions in  Subsections \ref{sec:P1t} and \ref{sec:Pst}.
Additionally, we were able to shorten some of the proofs.

The goal of the second part, Sections \ref{sec:operators} and \ref{sec:tt}, is to obtain the blackbox estimates, which are essential for the proofs in Sections \ref{sec:classical harmonic extension} and \ref{sec:s harmonic extension}.
In Section \ref{sec:operators}, we collect some well known facts about Riesz transforms, fractional Laplacians and Riesz potentials and prove some specific, elementary Lemmas in preparation for Section \ref{sec:tt}.
If the reader is not familiar with these operators, we recommend to check out this section first.

In Section \ref{sec:tt}, we prove the blackbox estimates by collecting a variety of building blocks, which can be combined to obtain these blackbox estimates. 
There are two different sources for these building blocks. 
First, in \cite{BuiCan17}, Bui and Candy showed that Triebel-Lizorkin and Besov-Lipschitz spaces can be characterized with Poisson like kernels.
From these characterizations we obtain $BMO$, H\"older and fractional Sobolev space estimates for the harmonic extensions.
Second, from Stein's book (\cite{Ste93}) we acquire square function estimates as well as pointwise estimates in terms of the maximal function.
Additionally, we introduce Lorentz spaces together with an interpolation theorem, with which we obtain finer Lorentz estimates from $L^p$-estimates.
The Section is structured in the following way.
\begin{itemize}[noitemsep]
\item Subsection \ref{subsec:tt - triebel lizorkin besov}: Triebel-Lizorkin and Besov-Lipschitz spaces and the Bui-Candy result.
\item Subsection \ref{subsec:tt - trilibeli space char}: Triebel-Lizorkin and Besov-Lipschitz space characterizations.
\item Subsection \ref{subsec:tt - lorentz spaces}: Lorentz spaces and an interpolation theorem.
\item Subsection \ref{subsec:tt - building blocks}: $BMO$, fractional Sobolev and H\"older space estimates.
\item Subsection \ref{subsec:tt - square function estimates}: Square function estimates.
\item Subsection \ref{subsec:tt - maximal function estimates}: The maximal function and pointwise estimates.
\item Subsection \ref{subsec:tt - blackbox estimates}: The blackbox estimates.
\end{itemize}
While the main ideas are also based on \cite{LeSchi20}, in these last two sections we take a wider approach, going into much more detail regarding the tools and the building blocks necessary for the blackbox estimates. 
By gathering some additional building blocks, we are able to further generalize the final blackbox estimates while at the same time presenting them clearer and more intuitive for applications.

Before diving into the first example, we would like to thank A. Schikorra for his quick and comprehensive answers to our questions throughout the work on this thesis.

\subsection{A first example: Jacobian estimates}\label{sec:JacEs}
For a function $u\colon \Rn\rightarrow\Rn$, we have the Jacobian $\det(\nabla u)$.
It is obvious, that 
\begin{equation*}
	\|\det(\nabla u)\|_{L^1(\Rn)}\lesssim \|\nabla u\|_{L^n(\Rn)}^n.
\end{equation*}
Making use of the determinants inherent structure, we obtain better regularity, that is
\begin{equation*}
	\|\det(\nabla u)\|_{\mathcal{H}^1(\Rn)}\lesssim \|\nabla u\|_{L^n(\Rn)}^n,
\end{equation*}
which is a consequence of (\ref{eq:JacEs1}) due to the duality of the Hardy space and $BMO$, see \cite{FeSte72}.
This estimate was originally due to Coifman, Lions, Meyer, and Semmes, see \cite[Theorem II.1.1), p.250]{CLMS93}.
In the same work, they showed that this Jacobian estimate actually is a special case of the commutator estimate we show in Subsection \ref{sec:Rt comms}, confer \cite[Section III.1, pp.257-258]{CLMS93}.
We stated above that with our method the cancellation effects responsible for the better estimate can be easily observed during the third step of the proof. 
This estimate is an exception though, since we already use the special determinant structure in the second step. 
By coding the Jacobian as differential form, we are able to apply Stokes' Theorem instead of a standard integration by parts.
For a definition of the occurring space norms, see Subsection \ref{subsec:tt - building blocks}, i.e. (\ref{eq:tt def Slobodeckij}) and (\ref{eq:tt def BMO}).

\begin{theorem}[Jacobian estimate]\label{th:JacEs}
Let $\phi\in\CicRn$ and $u=(u^1,\ldots,u^n)\in\CicRnRn$. 
Then the following estimate holds.
\begin{equation}\label{eq:JacEs1}
	\int_\Rn \phi \det(\nabla u) \lesssim [\phi]_{BMO}\|\nabla u\|^n_{L^n(\Rn)}
\end{equation}
Moreover, let $0<s_i<1$, $1<p_i<\infty$ for $i=0,\ldots,n$ satisfy
\begin{equation}\label{eq:JacEs2Req}
	\sum_{i=0}^n s_i = n, \quad
	\sum_{i=0}^n \frac{1}{p_i} = 1.
\end{equation}
Then
\begin{equation}\label{eq:JacEs2}
	\int_\Rn \phi \det(\nabla u) \lesssim [\phi]_{W^{s_0,p_0}} [u^1]_{W^{s_1,p_1}}\ldots[u^n]_{W^{s_n,p_n}}.
\end{equation}
\end{theorem}

\begin{proof}[Estimate (\ref{eq:JacEs2})]
As explained above, this proof will be divided into 3 steps. 
First, we extend the occurring functions to $\Rnpp$. 
Here, let $\Phi\colon \Rnpp\rightarrow\R$, $U\colon \Rnpp\rightarrow\Rn$ be the usual harmonic extensions to $\Rnpp$ of $\phi$ and $u$, meaning they fulfill
\begin{mathex}
	&\begin{cases}
		\Delta_\Rnp \Phi \equiv (\Delta_x+\partial_{tt})\Phi = 0	&\text{in } \Rnpp,\\
		\Phi(x,0) = \phi(x)										&\text{in } \Rn,
	\end{cases}\\
	&\begin{cases}
		\Delta_\Rnp U \equiv (\Delta_\Rnp U_i)_i = 0				&\text{in } \Rnpp,\\
		U(x,0) = u(x)											&\text{in } \Rn.
	\end{cases}
\end{mathex}
In order for the integration by parts formulas to be applicable, the extensions need to sufficiently decrease for $|(x,t)|\rightarrow \infty$. 
Choosing zero-boundary data at infinity, $\Phi$ and $U$ are obtained via the Poisson operator given in Definition \ref{def:P1t Poisson operator}, i.e. $\Phi(x,t)\coloneqq P^1_t\phi(x)$ and $U(x,t)\coloneqq P^1_t u(x)$.

For the second step, we integrate by parts to obtain an integral over the upper half-space in terms of $\Phi$ and $U$. 
Here, via Stokes' theorem, see Theorem \ref{th:Stoke Rnpp}, we have
\begin{equation}\label{eq:partInt JacEs2}
	\mathcal{I}\coloneqq 	\left|\int_\Rn \phi~\det_{n\times n}(\nabla_\Rn u)  \right|
						=	\left|\int_\Rnpp \det_{(n+1)\times (n+1)}(\nabla_\Rnp \Phi,\nabla_\Rnp U)  \right|
\end{equation}
since, according to Lemma \ref{lem:JacEs diff form encoding} below, we can rewrite both sides as 
\begin{equation*}
	\int_\Rn \phi~\det_{n\times n}(\nabla_\Rn u) 
	= \int_{\partial\Rnpp} \phi~du^1\wedge du^2 \wedge\ldots\wedge du^n,
\end{equation*}
\begin{equation*}
	\int_\Rnpp \det_{(n+1)\times (n+1)}(\nabla_\Rnp \Phi,\nabla_\Rnp U) 
	= \int_\Rnpp d\Phi \wedge dU^1 \wedge dU^2 \wedge\ldots\wedge dU^n.
\end{equation*}
Additionally, in order to apply Stokes' theorem to $w\coloneqq \Phi \, dU^1 \wedge dU^2 \wedge\ldots\wedge dU^n$, we have to confirm that $\phi~du^1\wedge du^2 \wedge\ldots\wedge du^n$ is integrable on $\Rn$ and that $w$ and its first order derivatives are integrable on $\Rnpp$.
The integrability on $\Rn$ is obvious since $\phi, u^1,\ldots, u^n$ are compactly supported. 
Regarding the integrability on $\Rnpp$, $w$ and its derivatives are bounded according to Corollary \ref{cor:P1t boundedness of P1t derivatives}.
Additionally, the product of multiple harmonic extensions and their derivatives decays sufficiently fast thanks to Corollary \ref{cor:P1t decay of P1t and derviatives}.
Note that $d$ represents the $n+1$-dimensional exterior derivative on the left side of (\ref{eq:JacEs2}) and the $n$-dimensional one on the right side.
But since the pullback of $dt$ to the boundary is zero, we still have
\begin{equation*}
	\int_{\partial\Rnpp} \Phi \, dU^1 \wedge dU^2 \wedge\ldots\wedge dU^n = \int_{\Rn} \phi~du^1\wedge du^2 \wedge\ldots\wedge du^n.
\end{equation*}

In the third and last step, we apply one of the trace theorems from Section \ref{sec:tt} to obtain the integral (semi-)norms on $\Rn$.
In this case, the seminorms $[\cdot]_{W^{s_i,p_i}}$ can be obtained via (\ref{eq:tt Slobodeckij via partial_tFs}) and (\ref{eq:tt Slobodeckij via nabla_xFs}) from Proposition \ref{prop:tt slobodeckij space char}. 
Combined, these two estimates yield 
\begin{equation*}
	\left(\int_\Rn \int_0^\infty |t^{1-\frac{1}{p_0}-s_0}\nabla_{\Rnp}\Phi(x,t)|^{p_0}~dt~dx\right)^\frac{1}{p_0} 
	\approx [\phi]_{W^{s_0,p_0}}
\end{equation*}
and, respectively for $i=1,\ldots,n$,
\begin{equation*}
	\left(\int_\Rn \int_0^\infty |t^{1-\frac{1}{p_i}-s_i}\nabla_{\Rnp} U^i(x,t)|^{p_i}~dt~dx\right)^\frac{1}{p_i} 
	\approx [u^i]_{W^{s_i,p_i}}.
\end{equation*}
We only need to transform $\mathcal{I}$ into a fitting form.
Estimating the determinant against its arguments' norms and using Hölder's inequality, we obtain
\begin{mathex}[LCL]
	\mathcal{I} 	&\lesssim &
			\int_\Rnpp |\nabla_{\Rnp} \Phi(x,t)||\nabla_{\Rnp} U^1(x,t)|\ldots |\nabla_{\Rnp} U^n(x,t)|\,dx\,dt \\
	&\overset{(\ref{eq:JacEs2Req})}{=}&
			\int_\Rn \int_0^\infty 	t^{n+1-\sum_{i=0}^n\frac{1}{p_i}-\sum_{i=0}^n s_i}
								|\nabla_{\Rnp} \Phi(x,t)||\nabla_{\Rnp} U^1(x,t)|\ldots |\nabla_{\Rnp} U^n(x,t)|\,dt\,dx \\
	&\leq 	&\left(\int_\Rn \int_0^\infty |t^{1-\frac{1}{p_0}-s_0}\nabla_{\Rnp}\Phi(x,t)|^{p_0}~dt~dx\right)^\frac{1}{p_0}\\
		&	&\cdot\left(\int_\Rn \int_0^\infty |t^{1-\frac{1}{p_i}-s_1}\nabla_{\Rnp} U^1(x,t)|^{p_1}\,dt\,dx\right)^\frac{1}{p_1}\\
		&	&\vdots \\
		&	&\cdot\left(\int_\Rn \int_0^\infty |t^{1-\frac{1}{p_n}-s_n}\nabla_{\Rnp} U^n(x,t)|^{p_n}\,dt\,dx\right)^\frac{1}{p_n}\\
	&\approx	&[\phi]_{W^{s_0,p_0}} [u^1]_{W^{s_1,p_1}}\ldots[u^n]_{W^{s_n,p_n}}.
\end{mathex}
Thus, we have shown that (\ref{eq:JacEs2}) holds.
\end{proof}

\begin{proof}[Estimate (\ref{eq:JacEs1})]
First, let $\Phi$ and $U$ be the same harmonic extensions of $\phi$ and $u$ as in the proof of (\ref{eq:JacEs2}) above.
In the second step, we start with the same integration by parts and again have
\begin{mathex}
	\mathcal{I}\coloneqq 	
					\left|\int_\Rn \phi~\det_{n\times n}(\nabla_\Rn u)  \right|
				=	\left|\int_\Rnpp \det_{(n+1)\times (n+1)}(\nabla_\Rnp U,\nabla_\Rnp \Phi)  \right|
\end{mathex}
This time, we need an additional derivative to apply the trace theorem we are aiming for. Therefore, we integrate by parts in $t$-direction to obtain an additional $t$-derivative,
\begin{mathex}
	\mathcal{I}	
			&= 	&\left| \int_\Rn \left(
					\lim_{T\rightarrow \infty} \int_{\frac{1}{T}}^T \det(\nabla_{\Rnp} U(x,t), \nabla_{\Rnp} \Phi(x,t))~dt
				\right) dx \right| \\
			&=	&\left| \int_\Rn \lim_{T\rightarrow \infty} \left(
					\left[t\det(\nabla U(x,t), \nabla \Phi(x,t))\right]_\frac{1}{T}^T 
					- \int_{\frac{1}{T}}^T t\partial_t\det(\nabla U(x,t), \nabla \Phi(x,t))~dt
				\right) dx \right|
\end{mathex}
Since the harmonic extensions $U$ and $\Phi$ are bounded, see Corollary \ref{cor:P1t decay of P1t and derviatives}, we have 
\begin{equation*}
	\lim_{t\rightarrow \infty} t|\nabla U^1(x,t)| \ldots |\nabla U^n(x,t)| |\nabla\Phi(x,t)| = 0.
\end{equation*}
Thanks to Corollary \ref{cor:P1t boundedness of P1t derivatives}, the decay-estimate for the harmonic extensions, we obtain 
\begin{equation*}
	\lim_{t\rightarrow 0} t|\nabla U^1(x,t)| \ldots |\nabla U^n(x,t)| |\nabla\Phi(x,t)| = 0
\end{equation*}
for the other limit. 
Therefore, we conclude the second step having shown that
\begin{equation*}
	\mathcal{I} = \left| \int_\Rnpp t \partial_t \det(\nabla_\Rnp U^1,\ldots,\nabla_\Rnp U^n,\nabla_\Rnp \Phi) \,dx\,dt \right|.
\end{equation*}

For the third step, it suffices to show that
\begin{equation}\label{eq:JacEs1 req trace theorem}
	\mathcal{I} \lesssim \int _\Rnpp t|\nabla_\Rnp \Phi|\,|\nabla_\Rnp \nabla_x U|\,|\nabla_\Rnp U|^{n-1} \,dx\,dt.
\end{equation}
Basically, we need to show that the added $t$-derivative can be converted to a derivative in $x$ and that, when it hits $\Phi$, the derivative can be redistributed to instead hit a $U$-term.
Should we have (\ref{eq:JacEs1 req trace theorem}), we can apply Proposition \ref{prop:tt - BMO est} for $s=1$ to obtain
\begin{equation*}
	\mathcal{I}	\lesssim	[\phi]_{BMO(\Rn)}\|(-\Delta)^\frac{1}{2}u\|^n_{L^n(\Rn)}
			\approx	[\phi]_{BMO(\Rn)}\|\nabla u\|^n_{L^n(\Rn)}
\end{equation*}
where the second estimate is due to Lemma \ref{prop:fLl vs grad}.

Let us prove (\ref{eq:JacEs1 req trace theorem}).
Recalling the Leibniz formula for the determinant, we observe that
\begin{mathex}
	\mathcal{I}	&\leq		&\left| \sum_{i=1}^n \int_\Rnpp t \det(\nabla_\Rnp U^1,\ldots,\partial_t\nabla_\Rnp U^i,\ldots,\nabla_\Rnp U^n,\nabla_\Rnp \Phi) \right|\\
				&			&+ \left| \int_\Rnpp t \det(\nabla_\Rnp U^1,\ldots,\nabla_\Rnp U^n,\partial_t\nabla_\Rnp \Phi) \right| \\
				&\eqqcolon	&\mathcal{I}_1 + \mathcal{I}_2.
\end{mathex}
Therefore, we have to consider two cases.
Regarding $\mathcal{I}_1$, with $\partial_{tt}U = -\Delta_x U$ we have
\begin{mathex}
	|\partial_t\nabla_\Rnp U^i| 
	\leq 		|\nabla_\Rnp^2 U| 
	&\approx 	&|\partial_{tt}U| + \sum_{i=1}^n ( |\partial_t \partial_{x_i}U| + |\partial_{x_i} \partial_t U|) + |\nabla_x^2 U|\\
	&\leq 		&|\nabla_x^2 U| + 2|\nabla_\Rnp \nabla_x U| + |\nabla_x^2 U| \\
	&\lesssim	&|\nabla_\Rnp \nabla_x U|.
\end{mathex}
Thus, we obtain
\begin{mathex}
	\mathcal{I}_1	&\leq \sum_{i=1}^n \int_\Rnpp t\,|\nabla_\Rnp U|^{n-1}\,|\partial_t \nabla_\Rnp U|\,|\nabla_\Rnp \Phi| \\
					&\lesssim \int _\Rnpp t\,|\nabla_\Rnp U|^{n-1}\,|\nabla_\Rnp \nabla_x U|\,|\nabla_\Rnp \Phi|.
\end{mathex}
Regarding $\mathcal{I}_2$, we have to redistribute the additional derivative from $\nabla_\Rnp \Phi$ to a $U$-term.
We rename the variables $(z_1,\ldots,z_{n+1})=(x_1,\ldots,x_n,t)$.
Considering the Leibniz formula for the determinant, we obtain
\begin{mathex}
	\mathcal{I}_2	
	&=		&\left| \int_\Rnpp t \det(\nabla_\Rnp U^1,\ldots,\nabla_\Rnp U^n,\partial_t\nabla_\Rnp \Phi) \,dx\,dt\right|\\
	&\leq	&\sum_{i_1,\ldots,i_n=1}^{n+1}\sum_{k=1}^n \left|
				 \int_\Rnpp z_{n+1}\,\partial_{z_{i_1}}U^1 \ldots \partial_{z_{i_n}}U^n\,\partial_{z_{n+1}}\partial_{z_k}\Phi\,dx\,dt \right|\\
	&		&+ \sum_{i_1,\ldots,i_n=1}^{n+1} \left|
				 \int_\Rnpp z_{n+1}\,\partial_{z_{i_1}}U^1 \ldots \partial_{z_{i_n}}U^n\,\partial_{z_{n+1}}\partial_{z_{n+1}}\Phi\,dx\,dt \right|.
\end{mathex}
By harmonicity, $\partial_{z_{n+1}}\partial{z_{n+1}}\Phi=-\sum_{l=1}^n\partial_{z_l z_l}\Phi$ and we have
\begin{mathex}
	\mathcal{I}_2	
	&\leq	&\sum_{i_1,\ldots,i_n=1}^{n+1}\sum_{k=1}^n \left|
				 \int_\Rnpp z_{n+1}\,\partial_{z_{i_1}}U^1 \ldots \partial_{z_{i_n}}U^n\,\partial_{z_{n+1}}\partial_{z_k}\Phi\,dx\,dt \right|\\
	&		&+ \sum_{i_1,\ldots,i_n=1}^{n+1} \sum_{l=1}^n \left|
				 \int_\Rnpp z_{n+1}\,\partial_{z_{i_1}}U^1 \ldots \partial_{z_{i_n}}U^n\,\partial_{z_l}\partial_{z_l}\Phi\,dx\,dt \right|.
\end{mathex}
We redistribute the additional derivative with an integration by parts in $z_k$ for the first term, and in $z_l$ for the second term.
We denote $(z_1,\ldots,z_{i-1},z_{i+1},\ldots z_{n+1})$ by $\tilde{z_i}\in\R^n_+$ for $1\leq i \leq n$ and obtain
\begin{mathex}
	\mathcal{I}_2	
	&\leq	&\sum_{i_1,\ldots,i_n=1}^{n+1}\sum_{k=1}^n \left|
				\int_{\R^n_+} \left( \lim_{Z\rightarrow\infty} 
				 	\int_{-Z}^Z \partial_{z_k}(z_{n+1}\,\partial_{z_{i_1}}U^1\ldots\partial_{z_{i_n}}U^n)\,\partial_{z_{n+1}}\Phi ~dz_k					\right) d\tilde{z_k} \right|\\
	&		&+ \sum_{i_1,\ldots,i_n=1}^{n+1} \sum_{l=1}^n \left|
				 \int_{\R^n_+} \left( \lim_{Z\rightarrow\infty} 
				 	\int_{-Z}^Z \partial_{z_l}(z_{n+1}\,\partial_{z_{i_1}}U^1\ldots\partial_{z_{i_n}}U^n)\,\partial_{z_{n+1}}\Phi ~dz_l					\right) d\tilde{z_l} \right|\\
	&=		&\sum_{i_1,\ldots,i_n=1}^{n+1}\sum_{k=1}^n \left|
				 \int_\Rnpp z_{n+1}\,\partial_{z_k}(\partial_{z_{i_1}}U^1\ldots\partial_{z_{i_n}}U^n)\,\partial_{z_{n+1}}\Phi \,\right|\\
	&		&+ \sum_{i_1,\ldots,i_n=1}^{n+1} \sum_{l=1}^n \left|
				 \int_\Rnpp z_{n+1}\,\partial_{z_l}(\partial_{z_{i_1}}U^1 \ldots \partial_{z_{i_n}}U^n)\,\partial_{z_l}\Phi \,\right|.
\end{mathex}
The boundary terms for the above integration-by-parts in $z_k$ and $z_l$ direction disappear, since $k,l\leq n$ and the harmonic extension decays to zero sufficiently fast as $|z_j|\rightarrow \infty$, $j\leq n$, see Corollary \ref{cor:P1t decay of P1t and derviatives}.
We conclude that 
\begin{equation*}
	\mathcal{I}_2 \lesssim	\int _\Rnpp z_{n+1}\,|\nabla_\Rnp U|^{n-1}\,|\nabla_\Rnp \nabla_x U|\,|\nabla_\Rnp \Phi|.
\end{equation*}
Thus, we have established (\ref{eq:JacEs1 req trace theorem}) and therefore proven (\ref{eq:JacEs1}).
\end{proof}

Rounding off this first application of the method, we justify (\ref{eq:partInt JacEs2}), the main  integration-by-parts formula used in the above proofs.

\begin{lemma}\label{lem:JacEs diff form encoding}
Let $U\subset\Rn$ be open, $v^1,\ldots,v^n\in\C^\infty(U)$. Then
\begin{mathex}
	dv^1\wedge\ldots\wedge dv^n = \det(\nabla v^1,\ldots \nabla v^n)~dx^1\wedge\ldots\wedge dx^n.
\end{mathex}
\end{lemma}
	
\begin{proof}
Using the Einstein summation convention, with the rules for the exterior product and the exterior derivative we obtain
\begin{mathex}
	dv^1\wedge\ldots\wedge dv^n 
	&=&\left(\frac{\partial v^1}{\partial x^i}~dx^i\right) \wedge\ldots\wedge \left(\frac{\partial v^1}{\partial x^i}~dx^i\right)\\
	&=&\sum_{\sigma\in\mathcal{S}_n} 
			\frac{\partial v^1}{\partial x^{\sigma(1)}}\ldots\frac{\partial v^n}{\partial x^{\sigma(n)}} 
			~dx^{\sigma(1)}\wedge\ldots\wedge dx^{\sigma(n)}\\
	&=&\sum_{\sigma\in\mathcal{S}_n} 
			\operatorname{sgn}(\sigma)\frac{\partial v^1}{\partial x^{\sigma(1)}}\ldots\frac{\partial v^n}{\partial x^{\sigma(n)}} 
			~dx^1\wedge\ldots\wedge dx^n\\
	&=&\det(\nabla v^1,\ldots \nabla v^n)~dx^1\wedge\ldots\wedge dx^n.
\end{mathex}
\end{proof}

\newpage
\section{Commutator estimates via classical harmonic extension}\label{sec:classical harmonic extension}
In this section, we will prove the div-curl estimate by Coifman-Lions-Meyer-Semmes as well as an estimate for a commutator of Riesz transforms and pointwise multiplication.
We conclude the section with an $L^1$-estimate for a double commutator of the same type in one dimension.
First though, we recall some facts about the harmonic extension.

\subsection{Harmonic extension to $\Rnpp$ via the Poisson operator $P^1_t$}\label{sec:P1t}
A harmonic extension of a function $f\colon\Rn\rightarrow\R$ to $\Rnpp$  is a solution of the Poisson equation,
\begin{equation}\label{eq:P1t Poisson equation}
	\left\lbrace	\begin{array}{ll}
					\Delta_{\Rnpp} F(x,t) \equiv (\partial_{tt} + \Delta_x)F(x,t) = 0 	\quad & \text{in } \Rnpp, \\
					\lim_{t\rightarrow 0} F(x,t) = f(x)						\quad & \text{almost everywhere on } \Rn.
				\end{array}\right.
\end{equation}
Note that the Poisson equation has no unique solution since linear functions with zero-boundary-value could be added to a solution. Therefore, we demand that additionally $\lim_{|(x,t)|\rightarrow\infty} F(x,t) = 0$. 
Under this assumption, the solution is unique due to the maximum principle for harmonic functions.
We obtain this solution via the Poisson operator.

\begin{definition}[The Poisson operator]\label{def:P1t Poisson operator}
The Poisson extension operator $P^1_t$ is given via the convolution
\begin{equation*}
	P^1_t f(x)	\coloneqq C_n \int_{\Rn} \frac{t}{\left(|x-y|^2+t^2\right)^\frac{n+1}{2}} f(y) \,dy
				= C_n (p^1_t * f) (x)
\end{equation*}
for $f\in L^1(\Rn,\R^m)+L^\infty(\Rn,\R^m)$, where the kernel $p^1_t$ is given by 
\begin{equation*}
	p^1_t (x) \coloneqq \frac{t}{\left(|x|^2+t^2\right)^\frac{n+1}{2}}.
\end{equation*}
\end{definition}

the kernels $p^1_t(x)= t^{-n}p^1_1(t^{-1}x)=(p^1_1)_t(x)$ are dilations of $p^s_1$ and therefore for all $t>0$ we have
\begin{equation}\label{eq:P1t L1 norm of p1t}
	\|p^1_t\|_{L^1(\Rn)} = \int_\Rn p^1_t = \int_\Rn p^1_1 \eqqcolon \frac{1}{C_n}.
\end{equation} 
For the computation of $C_n$, see \cite[Example 2.1.13, pp.92-93]{Gra14}. 
The  $L_\infty$-norm of $p^1_t$ is given by
\begin{equation}
	\|p^1_t\|_{L^\infty(\Rn)} = \frac{t}{(0^2+t^2)^\frac{n+1}{2}} = t^{-n}. 
\end{equation}

For $f\in L^\infty(\Rn)\cap \C(\Rn)$, the Poisson operator solves the Poisson equation in the sense that $F^1(x,t)\coloneqq P^1_t f(x)$ fulfills (\ref{eq:P1t Poisson equation}).
We might also interpret the kernels $(p^1_t)_{t>0}$ as function $k^1$ on $\Rnpp$ with $k^1(x,t)\coloneqq p^1_t(x)$. 
For the details on why $F^1$ fulfills the Dirichlet boundary condition, see again\cite[Example 2.1.13]{Gra14}.
The harmonicity of $F^1$ is inherited from $k^1$ (cf. Section \ref{sec:Pst}) since the convolution commutates with the derivations, see the proofs of Lemma \ref{lem:P1t derivatives of P1t 2} and Lemma \ref{lem:Pst derivatives of Pst 1} for details. 

We make further use of this commutative property to describe the derivatives of the harmonic extension in two lemmas. 
In the first lemma we apply the derivatives to the convolution kernel.

\begin{lemma}\label{lem:P1t derivatives of P1t 1}
Let  $k_i\in\Nz$ for $0\leq i \leq n$. Set $\tilde{k}\coloneqq k_1+\ldots+k_n$ and $k\coloneqq \tilde{k}+k_0$.
For any $f\in L^\infty(\Rn)$ with $F^1(x,t)\coloneqq P^1_t f(x)$ we have
\begin{equation*}
	\partial_t^{k_0}\partial_{x_1}^{k_1}\dots\partial_{x_n}^{k_n} F(x,t) 
	= C_{n,s} \int_\Rn
				\left(\sum_{j=0}^{k}
						\frac{ t^{1-k_0} q_{j}(x-y)}{(|x-y|^2+t^2)^{\frac{n+1+\tilde{k}+j}{2}}}
				\right) f(y)\,dy,
\end{equation*} 
where the $q_{j}$ are homogeneous polynomials of degree $j$, which depend on $k_0,\dots,k_n$. 
In particular, we have $F^1\in\C^\infty(\Rnpp)$.
\end{lemma}

\begin{proof}
See the detailed proof for the generalized harmonic extensions, Lemma \ref{lem:Pst derivatives of Pst 1}.
\end{proof}

As a direct consequence of this first lemma, we get a result about the decay of the harmonic extension and its derivatives.

\begin{corollary}[Decay estimate for the harmonic extension]\label{cor:P1t decay of P1t and derviatives}
Let $f\in L^\infty(\Rn)$ and $k\in\Nz$.
Denote with $F^1(x,t)\coloneqq P^1_t f(x)$ the harmonic extension.
Then 
\begin{equation}\label{eq:P1t t-decay of P1t for Linfty}
	\sup_{(x,t)\in\Rnpp} t^k|\nabla_\Rnpp^k F^1(x,t)| \lesssim \|f\|_{L^\infty(\Rn)}.
\end{equation}
If additionally $f\in L^1(\Rn)$, then
\begin{equation}\label{eq:P1t t-decay of P1t for L1}
	\sup_{(x,t)\in\Rnpp} t^{n+k}|\nabla_\Rnpp^k F^1(x,t)| \lesssim \|f\|_{L^1(\Rn)}.
\end{equation}
Regarding the decay in $x$-direction, if for some $C>0$, $k>0$ we have $|f(x)|\leq C|x|^{-l}$, then
\begin{equation}\label{eq:P1t x-decay of P1t}
	t^{k}|\nabla_\Rnpp^k F^1(x,t)| \lesssim |x|^{-l} + t\|f\|_{L^1(\Rn)}|x|^{-n-1}.
\end{equation}
\end{corollary}

\begin{proof}
Regarding (\ref{eq:P1t t-decay of P1t for Linfty}), according to Lemma \ref{lem:P1t derivatives of P1t 1} or (\ref{eq:Pst pst derivatives}), for any $k_0,\dots k_n\in\Nz$ and $\tilde{k}\coloneqq k_1+\ldots+k_n$ we have 
\begin{equation*}
	|\partial_t^{k_0}\partial_{x_1}^{k_1}\dots\partial_{x_n}^{k_n} k^1(x,t) |
	=
	\left| \sum_{j=0}^{k} \frac{ t^{1-k_0} q_{j}(x)}{(|x|^2+t^2)^{\frac{n+1+\tilde{k}+j}{2}}}\right|
	\leq
	C t^{1-k_0-\tilde{k}}  (|x|^2+t^2)^{-\frac{n+1}{2}}.
\end{equation*}
Therefore, 
\begin{equation*}
	t^k|\nabla_\Rnpp^k F^1(x,t)| \leq 
	C \int_\Rn		t^1 (|x-y|^2+t^2)^{-\frac{n+1}{2}} |f(y)|\,dy 
	\leq C\|f\|_{L^\infty(\Rn)} \|p^1_1\|_{L^1(\Rn)}.
\end{equation*}
Regarding (\ref{eq:P1t t-decay of P1t for L1}), we analogously obtain
\begin{equation*}
	|\partial_t^{k_0}\partial_{x_1}^{k_1}\dots\partial_{x_n}^{k_n} k^1(x,t)|
	\leq C t^{-k_0-\tilde{k}-n}
\end{equation*}
and therefore 
\begin{equation*}
	t^{n+k}|\nabla_\Rnpp^k F^1(x,t)| \leq 
	C \int_\Rn	 |f(y)|\,dy 
	= C\|f\|_{L^1(\Rn)}.
\end{equation*}
Regarding the decay in $x$-direction, we have 
\begin{mathex}
	|\nabla_\Rnpp^k F^1(x,t)| &\leq &
	\int_{\Rn\setminus B_\frac{|x|}{2}(0)} |\nabla_\Rnpp^k k^1(y,t)| |f(x-y)| \,dy 
	+ \int_{B_\frac{|x|}{2}(0)} |\nabla_\Rnpp^k k^1(y,t)| |f(x-y)| \,dy\\
	&\leq&
	 C t^{1-k} \|f\|_{L^1(\Rn)}  \left(\frac{|x|}{2}\right)^{-n-1}
	 + C t^{-k}\|p^1_1\|_{L^1(\Rn)} \left(\frac{|x|}{2}\right)^{-l}
\end{mathex}
and so obtain (\ref{eq:P1t x-decay of P1t}).
\end{proof}

Now, we want to apply the $x$-derivatives to $f$. For the derivatives in $t$-direction, we need the Fourier transform of the kernel $p^1_t$. 

\begin{proposition}[Fourier transform of the Poisson Kernel]\label{prop:P1t kernel fourier symbol}
The Fourier transform of the Poisson kernel is 
\begin{equation*}
	\Fourier(p^1_t)(\xi) = e^{-2\pi t|\xi|}.
\end{equation*}
\end{proposition}

\begin{proof}
This result can be achieved via a direct computation, see for example Theorem 1.11. of \cite{Hao16}. Alternatively, (\ref{eq:P1t Poisson equation}) can be transformed into an ordinary differential equation by applying the $n$-dimensional Fourier transform, see for example Section A.1 of \cite{LeSchi20}.
\end{proof}

With this result, we can now turn to the second lemma about the derivatives of the harmonic extension. 
For the definition of the half-Laplacian $\fL{1}$, see Definition \ref{def:fL fL}.

\begin{lemma}\label{lem:P1t derivatives of P1t 2}
Let $k_i\in\Nz$ for $0\leq i \leq n$ and $f\in W^{\infty,\infty}(\Rn)$.
Set $F^1(x,t)\coloneqq P^1_t f(x)$.
Then 
\begin{equation*}
	\partial_t^{k_0}\partial_{x_1}^{k_1}\dots\partial_{x_n}^{k_n} F^1(x,t) 
	= (-1)^{k_0} P_t^1 (\fL{k_0}\partial_{x_1}^{k_1}\dots\partial_{x_n}^{k_n} f)(x).
\end{equation*}
In particular, we obtain 
\begin{equation}\label{eq:P1t partial t on boundary}
	\lim_{t\rightarrow 0} \partial_t F^1(x,t) = -\fL{1}f(x) \quad \text{on } \Rn.
\end{equation}
\end{lemma} 

\begin{proof}
Assuming $k_0=0$, we show Lemma \ref{lem:P1t derivatives of P1t 2} via induction over $\tilde{k}\coloneqq k_1+\ldots+k_n$.
In a second step, we then deal with the derivatives in $t$ direction.
For $k=0$ the statement is obvious.
Since $\partial_{x_1}^{k_1}\dots\partial_{x_n}^{k_n} f \in W^{\infty,\infty}(\Rn)$, it is enough to show that $\partial_{x_l}P_t^1 f = P_t^1(\partial_{x_l}f)$.
Thanks to Lebesgue's dominated convergence theorem, we have
\begin{mathex}
	\lim_{h\rightarrow 0} D_{h,x_l} P_t^1 f(x)	&=& \lim_{h\rightarrow 0}C_n \int_{^\Rn} p^1_t(y)\frac{f(x+he_l-y)-f(x-y)}{h}  dy\\
												&=& C_n\int_{Rn} p^1_t(y) \partial_{x_l}f(x-y) = P_t^1 (\partial_{x_l}f)(x)
\end{mathex}
where $p^1_t \|\partial_{x_l}f\|_{L^\infty(\Rn)}$ is an integrable majorant thanks to the mean value theorem.

For the derivatives in $t$-direction, we now only need to handle the case $k_0=1$. The reason is that $\partial_{tt}F^1(x,t) = -\Delta_x F^1(x,t) = P_t^1((-\Delta)f)(x)$ due to harmonicity, the linearity of $P_t^1$ and the result for $x$-derivatives we just proved. 
We show that $-\partial_{t} F^1(x,t) = P_t^1(\fL{1}f)(x)$ by applying the $n$-dimensional Fourier transform to both sides.
For the right side, we obtain 
\begin{mathex}
	\Fourier(- p_t^1\ast\fL{1}f)(\xi) 	= - C_n \Fourier(p^1_t)(\xi) \Fourier(\fL{1}f)(\xi)
									= - C_n e^{-2\pi t|\xi|}\, 2\pi|\xi| \Fourier(f)(\xi)
\end{mathex}
with Proposition \ref{prop:P1t kernel fourier symbol} and Definition \ref{def:fL fL}.
For the left side, we obtain
\begin{mathex}
	\Fourier(\partial_t F^1(x,t))(\xi)	= \frac{d}{dt} \Fourier(F^1(x,t)) = \frac{d}{dt} C_n e^{-t|\omega||\xi|}\Fourier(f)(\xi)
										= - 2\pi C_n |\xi| e^{-2\pi t|\xi|} \Fourier(f)(\xi)
\end{mathex}
where we can switch the Fourier transform and $\partial_t$ according to Lemma \ref{lem:P1t switching Fourier and partial_t} below. 
The local boundedness condition is fulfilled due to Lemma \ref{lem:P1t derivatives of P1t 1}.
\end{proof}

While as a consequence of Lemma \ref{lem:P1t derivatives of P1t 1} we got a corollary about the decay of the harmonic extension's derivatives, out of this second lemma we get the boundedness of these derivatives, at least combined with an $L^\infty$ estimate for the half-Laplacian, see for example Lemma \ref{lem:fL Linfty estimate}.

\begin{corollary}[$L^\infty$-estimate for the harmonic extension]\label{cor:P1t boundedness of P1t derivatives}
Let $k_i\in\Nz$ for $0\leq i \leq n$ and $f\in\CicRn$. Then 
\begin{equation*}
	\|\partial_t^{k_0}\partial_{x_1}^{k_1}\dots\partial_{x_n}^{k_n} P_t^1 f\|_{L^\infty(\Rn)} \leq 
	\|\fL{k_0}\partial_{x_1}^{k_1}\dots\partial_{x_n}^{k_n} f \|_{L^\infty(\Rn)}.
\end{equation*}
\end{corollary}

\begin{proof}
The result follows immediately with Lemma \ref{lem:P1t derivatives of P1t 2} and $C_n\|p^1_t\|_{L^1(\Rn)}=1$. 
\end{proof}

Wrapping this section up, we still need to make sure that we were allowed to switch the $n$-dimensional Fourier transform and the $t$-derivative in the proof of Lemma \ref{lem:P1t derivatives of P1t 2}.
The following lemma provides a condition under which this switch is possible.

\begin{lemma}\label{lem:P1t switching Fourier and partial_t}
Let $F\in \C^1(\Rnpp)$. 
Assume $t_0>0$. Suppose that there exists $0<\epsilon<t_0$ such that $\partial_t F$ and $|x|^{n+1} \partial_t F(x,t)$ are bounded on $\Rn \times B_\epsilon(t_0)$. 
Then
\begin{equation*}
	\left.\frac{d}{dt}\right|_{t=t_0} \Fourier(F(\cdot,t)) = \Fourier(\partial_t F(\cdot,t_0)).
\end{equation*}
\end{lemma}

\begin{proof}
Since $\partial_t F$ and $(x,t)\rightarrow |x|^{n+1} \partial_t F(x,t)$ are bounded, there exists $\const > 0$ such that for all $x\in\Rn$, $t\in B_\epsilon(t_0)$
\begin{equation*}
	\partial_t F(x,t) \leq \min\left\lbrace C, \frac{C}{|x|^{n+1}}\right\rbrace.
\end{equation*}
Applying Lebesgue's dominated convergence theorem with $\min\left\lbrace C, \frac{C}{|x|^{n+1}}\right\rbrace$ as majorant, we obtain
\begin{mathex}
	\left.\frac{d}{dt}\right|_{t=t_0} \Fourier(F(\cdot,t))(\xi) 
	&=& \left.\frac{d}{dt}\right|_{t=t_0} \left(\int_\Rn F(x,t)\,e^{-2\pi i x\cdot\xi}\,dx \right)\\
	&=& \lim_{h\rightarrow 0} \int_\Rn \frac{F(x,t_0+h)-F(x,t_0)}{h}\,e^{-2\pi i x\cdot\xi}\,dx\\
	&=& \int_\Rn \partial_tF(x,t_0)\,e^{-2\pi i x\cdot\xi}\,dx
	&= \Fourier(\partial_t F)(\xi).
\end{mathex}
\end{proof}

\subsection{The general div-curl estimate by Coifman-Lions-Meyer-Semmes}\label{sec:divcurl by CLMS}
The Jacobian estimate (\ref{eq:JacEs1}) is actually a special case of the general div-curl estimate below, which was proven by Coifman-Lions-Meyer-Semmes, see \cite[Theorem II.1, pp. 250-251]{CLMS93}. 
While curl is not explicitly mentioned in the formulation we chose, on $\Rn$ every vectorfield $v$ with $\operatorname{curl}(v)=0$ can be written as $v=\nabla f$ for a suitably chosen scalar field $f$. 
In addition to the original result, we again obtain an intermediate estimate.

\begin{theorem}[Coifman-Lions-Meyer-Semmes]\label{th:div-curl original result}
Assume $\phi,f\in\CicRn$, $g\in\CicRnRn$ with 
\begin{equation*}
	\operatorname{div}(g) = \sum_{i=1}^n\partial_i g_i = 0.
\end{equation*}
Suppose $1<p_1,p_2<\infty$, $1\leq q_1,q_2\leq\infty$ with $\frac{1}{p_1}+\frac{1}{p_2} = \frac{1}{q_1}+\frac{1}{q_2} =1$. 
Then
\begin{equation}\label{eq:divcurl original BMO estimate}
	\int_{\Rn}\sum_{i=1}^n \partial_i f \cdot g_i\,\phi 
	\lesssim  
	[\phi]_{BMO} \|\nabla f\|_{L^{(p_1,q_1)}} \|g\|_{L^{(p_2,q_2)}}.
\end{equation}
Suppose $0<s_1,s_2,s_3<1$ with $s_1+s_2+s_3=2$ and $1<p_1,p_2,p_3<\infty$, $1\leq q_1,q_2,q_3\leq\infty$ with $\frac{1}{p_1}+\frac{1}{p_2}+\frac{1}{p_3} = \frac{1}{q_1}+\frac{1}{q_2}+\frac{1}{q_3} = 1$.
Then
\begin{equation*}
	\int_{\Rn}\sum_{i=1}^n \partial_i f \cdot g_i\,\phi 
	\lesssim  
	\|(-\Delta)^\frac{s_1}{2}\phi\|_{L^{(p_1,q_1)}} \|(-\Delta)^\frac{s_2}{2}f\|_{L^{(p_2,q_2)}} \|I^{1-s_3}g\|_{L^{(p_3,q_3)}}.
\end{equation*}
\end{theorem}

For our method, it is advantageous to code the functions as differential forms, since this opens up Stokes' Theorem as a means to obtain the integral over $\Rnpp$.
Using Stokes' Theorem instead of a simpler integration by parts already makes use of the ``div-curl'' expression's inherent structure and therefore replaces the cancellation effects we would otherwise observe in more detail running through the elementary transformations.

Clearly, $f$ can be interpreted as $0$-form, coding $\nabla f$ as $df$. 
Likewise, $g$ can be interpreted as $(n-1)$-form, $g\in\C_c^\infty\left(\bigwedge^{n-1}\Rn\right)$.
By suitably choosing the signs for the components of $g$ in the encoding, we can then rewrite the left side of the estimates as
\begin{equation*}
	\int_{\Rn}\sum_{i=1}^n \partial_i f \cdot g_i\,\phi  = \int_\Rn \phi\,df\wedge g.
\end{equation*}
Due to the way we chose the signs for $g$, we can identify the divergence with the exterior derivative, and therefore have $dg=0$. 
Thus, by the Poincar\'e Lemma on differential forms, there is an $(n-2)$-form $h\in\C_c^\infty\left(\bigwedge^{n-2}\Rn\right)$ such that $g=dh$.

In the language of differential forms, the above theorem can be stated as follows.

\begin{theorem}\label{th:div-curl diffform result}
Let $l\in\lbrace 0,\ldots,n-2 \rbrace$.
Assume $\phi\in\CicRn$, $f\in\C_c^\infty(\bigwedge^l\Rn)$ and $h\in\C_c^\infty(\bigwedge^{n-l-2}\Rn)$.\\
Suppose $1<p_1,p_2<\infty$, $1\leq q_1,q_2\leq\infty$ with $\frac{1}{p_1}+\frac{1}{p_2} = \frac{1}{q_1}+\frac{1}{q_2} =1$. 
Then
\begin{equation}\label{eq:div-curl BMO}
	\int_{\Rn} \phi \, df \wedge dh
	\lesssim
	[\phi]_{BMO}\|\nabla f\|_{L^{(p_1,q_1)}} \|\nabla h\|_{L^{(p_2,q_2)}}.
\end{equation}
The second, intermediate estimate has the following form.
Suppose $0<s_1,s_2,s_3<1$ with $s_1+s_2+s_3=2$ and $1<p_1,p_2,p_3<\infty$, $1\leq q_1,q_2,q_3\leq\infty$ with $\frac{1}{p_1}+\frac{1}{p_2}+\frac{1}{p_3} = \frac{1}{q_1}+\frac{1}{q_2}+\frac{1}{q_3} = 1$.
Then
\begin{equation}\label{eq:div-curl fracSob}
	\int_{\Rn} \phi \, df \wedge dh
	\lesssim
	\|\fL{s_1}\phi\|_{L^{(p_1,q_1)}} \|\fL{s_2} f\|_{L^{(p_2,q_2)}} \|\fL{s_3} h\|_{L^{(p_3,q_3)}}.
\end{equation}
\end{theorem}

Before proving the theorem, let us briefly explain the norms in (\ref{eq:div-curl BMO}) and (\ref{eq:div-curl fracSob}).
We can write any $l$-form $f$ as 
\begin{equation*}
	f = \sum_{1\leq i_1<\ldots<i_l\leq n} f_{i_1,\ldots,i_l}\,dx^{i_1}\wedge\ldots\wedge dx^{i_l}.
\end{equation*}
We can extend all of our norms to $l$-forms by applying them to $f_{i_1,\ldots,i_l}$.
In particular, we have
\begin{equation*}
	\|\nabla f\|_{L^{(p,q)}} \coloneqq \sum_{1\leq i_1<\ldots<i_l\leq n} \|\nabla f_{i_1,\ldots,i_l}\|_{L^{(p,q)}}
\end{equation*}
and
\begin{equation*}
	\|\fL{s} f\|_{L^{(p,q)}} \coloneqq \sum_{1\leq i_1<\ldots<i_l\leq n} \|\fL{s} f_{i_1,\ldots,i_l}\|_{L^{(p,q)}}.
\end{equation*}
In the same manner, we may apply the Poisson operator $P^1_t$ to differential forms,
\begin{equation*}
	P^1_t f \coloneqq \sum_{1\leq i_1<\ldots<i_l\leq n} P^1_t f_{i_1,\ldots,i_l}\,dx^{i_1}\wedge\ldots\wedge dx^{i_l}.
\end{equation*}
Then $F(x,t)=P^1_t f(x)$ is an $l$-form on $\Rnpp$.

Notice that in general we only have
\begin{equation*}
	\|g\|_{L^{(p,q)}}=\|dh\|_{L^{(p,q)}}\lesssim \|\nabla h\|_{L^{(p,q)}}
\end{equation*}
and
\begin{equation*}
	\|I^{1-s}g\|_{L^{(p,q)}}=\|d (I^{1-s}h)\|_{L^{(p,q)}}
	\lesssim \|\nabla I^{1-s}h\|_{L^{(p,q)}}\approx \|\fL{1} I^{1-s}h\|_{L^{(p,q)}}= \|\fL{s}h\|_{L^{(p,q)}}.
\end{equation*}
Therefore, Theorem \ref{th:div-curl original result} clearly implies Theorem \ref{th:div-curl diffform result} for $l=0$.
The other direction is not as obvious. 
However, one can argue with the Gaffney inequality that Theorem \ref{th:div-curl diffform result} indeed also implies Theorem \ref{th:div-curl original result}.
We will not go into further detail though, since these considerations are out of our original scope. 

\begin{proof}[Theorem \ref{th:div-curl diffform result}]
Let $\Phi(x,t)=P^1_t\phi(x)$, $F(x,t) = P^1_tf(x)$ and $H(x,t)=P^1_th(x)$ be the harmonic extensions of $\phi$, $f$, $h$. 
By Stokes' Theorem on differential forms, see Theorem \ref{th:Stoke Rnpp}, we have
\begin{equation*}
	\mathcal{I}\coloneqq \left|\int_\Rn \phi\,df\wedge dh\, \right|
				= \left|\int_{\partial\Rnpp} \Phi\,dF\wedge dH\, \right|
				= \left|\int_{\Rnpp} d\Phi\wedge dF\wedge dH\, \right|.
\end{equation*}
The integrability conditions for Theorem \ref{th:Stoke Rnpp} are met due to $\phi,f,h\in\CicRn$ and the decay and boundedness of the harmonic extensions, see Corollary \ref{cor:P1t decay of P1t and derviatives} and Corollary \ref{cor:P1t boundedness of P1t derivatives}.

Regarding the intermediate estimate, we have
\begin{equation*}
	 \mathcal{I}\lesssim \int_\Rnpp |\nabla_\Rnp \Phi|\,|\nabla_\Rnp F|\,|\nabla_\Rnp H|
	 			= \int_\Rnpp t^{2-s_1-s_2-s_3}|\nabla_\Rnp \Phi|\,|\nabla_\Rnp F|\,|\nabla_\Rnp H|.
\end{equation*}
Thus, (\ref{eq:div-curl fracSob}) follows with Proposition \ref{prop:tt - Lp est}(a,b).

In order to obtain (\ref{eq:div-curl BMO}) from this term with Proposition \ref{prop:tt - BMO est}, we need another derivative on either $|\nabla_\Rnp F|$ or $|\nabla_\Rnp H|$.
An integration by parts in $t$ grants this additional derivative.
For an integrable function $E\in\C^1(\Rnpp)\cap L^\infty(\Rn)$ with integrable derivatives and $\lim_{t\rightarrow\infty}tE(x,t)=0$ we have
\begin{equation*}
	\int_\Rnpp E(x,t) \,dx\,dt = \int_\Rn \left[tE(x,t)\right]_0^\infty dx - \int_\Rnpp t\,\partial_t E(x,t)\,\,dx\,dt
			=- \int_\Rnpp t\,\partial_t E(x,t)\,\,dx\,dt
\end{equation*}
Applying this formula to $\mathcal{I}$, we have
\begin{equation*}
	\mathcal{I}= \left|\int_\Rnpp t\,\partial_t\left(d\Phi\wedge dF\wedge dH\right)\,\right|
\end{equation*}	
due to the decay and the boundedness of the harmonic extensions.
We claim that 
\begin{equation}\label{eq:div-curl BMO suitable form}
	\mathcal{I}\lesssim 
	\int_\Rnpp t\,|\nabla_\Rnp\Phi|\,\left(|\nabla_x\nabla_\Rnp F|\,|\nabla_\Rnp H| + |\nabla_x\nabla_\Rnp H|\,|\nabla_\Rnp F|\right).
\end{equation}
Essentially, we will show that we might redistribute the additional derivative if it hits $\Phi$, and replace the second $t$-derivative with $x$-derivatives. Then (\ref{eq:div-curl BMO}) follows with Proposition \ref{prop:tt - BMO est}. 

We rename the variables to $(z_1,\ldots,z_{n+1})$ from $(x_1,\ldots x_n,t)$ and have
\begin{equation*}
	d\Phi\wedge dF\wedge dH = \sum_{k,i,j=1}^{n+1}\sum_{I,J} 
									\partial_{z_k}\Phi\,\partial_{z_i}F_I\,\partial_{z_j}H_J
							\,dz^k\wedge dz^i \wedge dz^I \wedge dz^j\wedge dz^{J}
\end{equation*}
where the second sum is over all $I=(i_1,\ldots,i_l)$ with $i_1<i_2<\ldots<i_l$ and $J=(j_1,\ldots,j_{n-2-l})$ with $j_1<j_2<\ldots<j_{n-2-l}$. 
Distributing the $t$-derivative, we obtain three terms,
\begin{mathex}
	\partial_t(d\Phi\wedge dF\wedge dH) &=&
	\sum_{k,i,j=1}^{n+1}\sum_{I,J} \partial_{z_k}\Phi\,\partial_{z_{n+1}}\left(\partial_{z_i}F_I\,\partial_{z_j}H_J\right)
																\,dz^k\wedge dz^i \wedge dz^I \wedge dz^j\wedge dz^{J}\\
	&& + \sum_{i,j=1}^{n+1}\sum_{k=1}^n\sum_{I,J} \partial_{z_{n+1}}\partial_{z_k}\Phi\,\partial_{z_i}F_I\,\partial_{z_j}H_J
																\,dz^k\wedge dz^i \wedge dz^I \wedge dz^j\wedge dz^{J}\\
	&& + \sum_{i,j=1}^{n+1}\sum_{I,J} \partial_{z_{n+1}}\partial_{z_{n+1}}\Phi\,\partial_{z_i}F_I\,\partial_{z_j}H_J																		\,dz^{n+1}\wedge dz^i \wedge dz^I \wedge dz^j\wedge dz^{J}\\
	&\eqqcolon& \mathcal{J}_1+\mathcal{J}_2+\mathcal{J}_3.
\end{mathex}
With the harmonicity of the extensions in mind, i.e. $\partial_{z_{n+1}}\partial_{z_{n+1}}F_I = -\sum_{r=1}^n\partial_{z_r}\partial_{z_r}F_I$, we see that $\int_\Rnpp t\mathcal{J}_1$ can already be estimated as in (\ref{eq:div-curl BMO suitable form}).
Regarding the second term, we have $z_k=x_k$ and can therefore redistribute the additional $z_k$-derivative via an integration by parts in $z_k$ since both boundary terms disappear, see Corollary \ref{cor:P1t decay of P1t and derviatives}.
\begin{equation*}
	\left|\int_\Rnpp t\mathcal{J}_2\right|  = 
		\left|\int_\Rnpp z_{n+1}\sum_{i,j=1}^{n+1}\sum_{k=1}^n\sum_{I,J} 
					\partial_{z_{n+1}}\Phi\,\partial_{z_k}\left(\partial_{z_i}F_I\,\partial_{z_j}H_J\right)
																\,dz^k\wedge dz^i \wedge dz^I \wedge dz^j\wedge dz^{J}\,
		\right|
\end{equation*}
Obviously, this term can be estimated as in (\ref{eq:div-curl BMO suitable form}).
For the third term, we use the harmonicity of $\Phi$, that is $\partial_{z_{n+1}}\partial_{z_{n+1}}\Phi = -\sum_{r=1}^n\partial_{z_r}\partial_{z_r}\Phi$.
With an integration by parts in $z_r=x_r$ we obtain
\begin{equation*}
	\left|\int_\Rnpp t\mathcal{J}_3\right|  = 
		\left|\int_\Rnpp z_{n+1}\sum_{i,j=1}^{n+1}\sum_{I,J}\sum_{r=1}^n
						 \partial_{z_r		}\Phi\,\partial_{z_r}\left(\partial_{z_i}F_I\,\partial_{z_j}H_J	\right)								\,dz^{n+1}\wedge dz^i \wedge dz^I \wedge dz^j\wedge dz^{J}\,
		\right|
\end{equation*}
and have therefore shown that $\mathcal{I}$ can be estimated just as we claimed.
\end{proof}

\subsection{The Coifman-Rochberg-Weiss commutator estimate}\label{sec:Rt comms}
In this section, we examine the commutators between pointwise multiplication and the Riesz transforms $\Rj$, $j\in\lbrace 1,\ldots,n\rbrace$.
For $\phi\in\CicRn$, the commutator $[\Rj,\phi]$ is defined by
\begin{equation*}
	[\Rj,\phi](f)\coloneqq \Rj[\phi f] - \phi\Rj[f], \qquad f\in\CicRn.
\end{equation*}
Coifman-Rochberg-Weiss showed that these commutators are bounded operators, mapping $L^p(\Rn)$, $1<p<\infty$, onto itself with $[\phi]_{BMO}$ providing an upper bound for the operator norm.
They did not only prove this for the Riesz transforms, but for all Calder\'{o}n-Zygmund operators, see \cite[Section 2, Theorem 1, p.613]{CRW76}.
While we only prove this fact for the Riesz transforms, due to our method we easily obtain two additional intermediate estimates.

Actually, the div-curl estimate (\ref{eq:divcurl original BMO estimate}) and therefore also the Jacobian estimate (\ref{eq:JacEs1}) are special cases of the following $BMO$-estimate, see \cite[Section III.1, pp.257-258]{CLMS93}.

\begin{theorem}[Coifman-Rochberg-Weiss]\label{th:Rt comms CRW theorem}
Let $\phi,f\in\CicRn$ and $j\in\lbrace 1,\ldots,n\rbrace$. 
Suppose $p\in(1,\infty)$. Then
\begin{equation}\label{eq:Rt comms BMO}
	\| [\Rj,\phi](f) \|_{L^p(\Rn)} \lesssim [\phi]_{BMO} \|f\|_{L^p(\Rn)}.
\end{equation}
Moreover, for $p\in(1,\infty)$, $\sigma\in(0,1)$
\begin{equation}\label{eq:Rt comms Hölder}
	\| [\Rj,\phi](f) \|_{L^p(\Rn)} \lesssim [\fL{\sigma}\phi]_{BMO} \|I^\sigma f\|_{L^p(\Rn)}.
\end{equation}
Suppose $p,p_1,p_2\in(1,\infty)$ and $q_1,q_2\in [1,\infty]$ with $\frac{1}{p_1}+\frac{1}{p_2}=\frac{1}{q_1}+\frac{1}{q_2}=\frac{1}{p}$ and $\sigma\in[0,1)$. Then
\begin{equation}\label{eq:Rt comms fL&Rp}
	\| [\Rj,\phi](f) \|_{L^p(\Rn)} \lesssim \|\fL{\sigma}\phi\|_{L^{(p_1,q_1)}(\Rn)} \|I^\sigma f\|_{L^{(p_2,q_2)}(\Rn)}.
\end{equation}
\end{theorem}

The definition of the Riesz transforms and a collection of results regarding these operators can be found in Subsection \ref{subsec:Rt}.
Next to these results, we will use some specific interactions between Riesz transforms and the harmonic extension for the proof of the above theorem.

\begin{lemma}\label{lem:Rt comms P1t/Rt interaction}
Let $f\in\CicRn$ and $j\in\lbrace 1,\ldots,n\rbrace$. 
Denote with $F(x,t)= P^1_tf(x)$ the harmonic extension of $f$.
We write $\tilde{\mathcal{R}}_j[F](x,t) \coloneqq P^1_t\Rj f$.
Then 
\begin{equation*}
	\partial_t \tilde{\mathcal{R}}_j[F] = \partial_{x_j} F, \qquad
	\partial_{tt} \tilde{\mathcal{R}}_j[F] = \partial_t \partial_{x_j} F, \qquad
	\Delta_x \tilde{\mathcal{R}}_j[F] = -\partial_t \partial_{x_j} F.
\end{equation*}
\end{lemma}

\begin{remark}
Note that when considering the Poisson operator $P^1_t$ for fixed $t>0$, since they are Fourier multipliers, $P^1_t$ and the $n$-dimensional Riesz transform $\Rj$ commutate.
We connect this with a brief warning: While $\mathcal{\tilde{R}}_j$ can be interpreted as operator acting on $\Rnpp$, this is \textit{not} an actual $(n+1)$-dimensional Riesz transform. 
\end{remark}

\begin{proof}[Lemma \ref{lem:Rt comms P1t/Rt interaction}]
Due to Lemma \ref{lem:Pst derivatives of Pst 2} we have 
\begin{equation*}
	\partial_t \tilde{\mathcal{R}}_j[F](x,t) = - P^1_t(\fL{1}\Rj f) =  P^1_t(\partial_{x_j} f) = \partial_{x_j}F,
\end{equation*}
where we used that $\fL{1}\Rj f= -\partial_{x_j} f$. 
This can easily be seen for the Fourier transforms.
The Fourier symbol of $\fL{1}$ is $2\pi|\xi|$ while the symbol of $\Rj$ is $-i\frac{\xi_j}{|\xi|}$, see Definition \ref{def:fL fL} and Proposition \ref{prop:Rt symbol}. Thus,
\begin{equation*}
	\Fourier(\fL{1}\Rj f)(\xi) = - 2\pi i|\xi| \frac{\xi_j}{|\xi|} \Fourier(f)(\xi)
	= -2\pi i\xi_j \Fourier(f) = \Fourier(-\partial_{x_j}f)(\xi).
\end{equation*}
The second equation of the lemma immediately follows from the first equation, the third by harmonicity from the second.
\end{proof}

Equipped with these additional tools, we can now prove Theorem \ref{th:Rt comms CRW theorem}.

\begin{proof}[Theorem \ref{th:Rt comms CRW theorem}] 
By $L^p$-duality, it suffices to show that for $p'\in(1,\infty)$ with $\frac{1}{p}+\frac{1}{p'}=1$ and for any $g\in\CicRn$
we have
\begin{equation}\label{eq:Rt comms estimate via duality}
	\mathcal{I}\coloneqq \left|\int_\Rn [\Rj,\phi](f)\,g\,\right| \lesssim 
	\left\lbrace
	\begin{array}{ll}
		{[}\phi{]}_{BMO}\,\|f\|_{L^p(\Rn)}\, \|g\|_{L^{p'}(\Rn)},\\
		{[}\fL{\sigma}\phi{]}_{BMO}\,\|I^\sigma f\|_{L^p(\Rn)}\, \|g\|_{L^{p'}(\Rn)},\\
		\|\fL{\sigma}\phi\|_{L^{(p_1,q_1)}(\Rn)}\, \|I^\sigma f\|_{L^{(p_2,q_2)}(\Rn)}\,\|g\|_{L^{p'}(\Rn)}.
	\end{array}
	\right.
\end{equation}
Let $\Phi(x,t)=P^1_t\phi(x)$, $F(x,t)=P^1_t f(x)$, $G(x,t)=P^1_t g(x)$ be the harmonic extensions of $\phi,f,g$.
Due to Lemma \ref{lem:Rt integration by parts}, the integration-by-parts property of the Riesz transforms, we have
\begin{equation*}
	\mathcal{I}	= \left| \int_\Rn \Rj[\phi f]g - \phi\Rj[f]g \right|
				= \left| \int_\Rn \phi f\Rj[g] + \phi\Rj[f]g \right|
\end{equation*}
Applying integration by parts twice in $t$-direction, we obtain
\begin{equation*}
	\mathcal{I}	
	= \left| \int_{\Rnpp} \partial_t\left(\Phi F\tilde{\mathcal{R}}_j[G] + \Phi\tilde{\mathcal{R}}_j[F]G\right) \right|
	= \left| \int_{\Rnpp} t \partial_{tt}\left(\Phi F\tilde{\mathcal{R}}_j[G] + \Phi\tilde{\mathcal{R}}_j[F]G\right) \right|
\end{equation*}
since for all $x\in\Rn$ the harmonic extensions satisfy the boundary condition of the Dirichlet problem,
\begin{equation*}
	\lim_{t\rightarrow 0} \left(\Phi F\tilde{\mathcal{R}}_j[G] + \Phi\tilde{\mathcal{R}}_j[F]G\right)(x,t)
	= \phi(x) f(x)\Rj[g](x) + \phi(x)\Rj[f](x)g(x).
\end{equation*}
All other boundary terms disappear since due to Corollary \ref{cor:P1t decay of P1t and derviatives} and Corollary \ref{cor:P1t boundedness of P1t derivatives} we have
\begin{equation*}
	\lim_{t\rightarrow \infty}  \left(\Phi F\tilde{\mathcal{R}}_j[G] + \Phi\tilde{\mathcal{R}}_j[F]G\right)(x,t) = 0,
\end{equation*}
\begin{equation*}
	\lim_{t\rightarrow \infty}  t\partial_t\left(\Phi F\tilde{\mathcal{R}}_j[G] + \Phi\tilde{\mathcal{R}}_j[F]G\right)(x,t)
	=\lim_{t\rightarrow 0} t\partial_t\left(\Phi F\tilde{\mathcal{R}}_j[G] + \Phi\tilde{\mathcal{R}}_j[F]G\right)(x,t) 
	=0.
\end{equation*}
For the rest of the proof, we will establish that
\begin{equation}\label{eq:Rt comms sutiable form for trace theorems}
	\mathcal{I} \lesssim 
	\max_{\tilde{F}\in\lbrace F,\tilde{\mathcal{R}}_jF\rbrace} \max_{\tilde{G}\in\lbrace G,\tilde{\mathcal{R}}_jG\rbrace}
	t |\nabla_\Rnpp\Phi|\,\left( |\nabla_\Rnpp \tilde{F}|\,|\tilde{G}|+ |\nabla_\Rnpp \tilde{G}|\,|\tilde{F}| \right).
\end{equation}
Then due to Proposition \ref{prop:Rt Lp bounded}, the $L^p$-boundedness of the Riesz transforms, we immediately obtain (\ref{eq:Rt comms estimate via duality}) with Proposition \ref{prop:tt - BMO est}(c,d), Proposition \ref{prop:tt - Hölder est}(e,f,g) and Proposition \ref{prop:tt - Lp est}(a,b) respectively.
For (\ref{eq:Rt comms fL&Rp}) we need the $L^{(p,q)}$-boundedness of the Riesz transforms, which we easily obtain from the $L^p$-boundedness via interpolation, see Theorem \ref{th:tt lorentz interpolation}. 

Let us prove (\ref{eq:Rt comms sutiable form for trace theorems}).
Distributing the two $t$-derivatives, we obtain three terms,
\begin{mathex}
	\mathcal{I} &\leq& 
		\left|\int_{\Rnpp} t \left(\partial_t\Phi\, \partial_t(F\,\tilde{\mathcal{R}}_j[G]) 
								+ \partial_t\Phi\, \partial_t(\tilde{\mathcal{R}}_j[F]\,G)\right) \right|\\
	&& +	\left|\int_{\Rnpp} t\left(\partial_{tt}\Phi\,(F\,\tilde{\mathcal{R}}_j[G] +\tilde{\mathcal{R}}_j[F]\,G)\right)\right|
	 +\left|\int_{\Rnpp} t\left(\Phi\,\partial_{tt}(F\,\tilde{\mathcal{R}}_j[G] +\tilde{\mathcal{R}}_j[F]\,G)\right)\right|\\
	&&\eqqcolon \mathcal{I}_1+\mathcal{I}_2+\mathcal{I}_3.
\end{mathex}
The term $\mathcal{I}_1$ can already be estimated as in (\ref{eq:Rt comms sutiable form for trace theorems}) by further distributing the derivative to $F$ and $G$.
For the second term, we apply the harmonicity of the extensions.
With $\partial_{tt}\Phi = -\Delta_x\Phi = -\nabla_x\cdot\nabla_x\Phi$ we have
\begin{equation*}
	\mathcal{I}_2 
	= 	\left|\int_{\Rnpp} t \left(\Delta_x\Phi (F\tilde{\mathcal{R}}_j[G] + \tilde{\mathcal{R}}_j[F]G)\right)\right|
	= 	\left|\int_{\Rnpp} t \left(\nabla_x\Phi\cdot\nabla_x(F\tilde{\mathcal{R}}_j[G] + \tilde{\mathcal{R}}_j[F]G)\right)\right|
\end{equation*}
Again, by further distributing the derivatives we can estimate this term as in (\ref{eq:Rt comms sutiable form for trace theorems}).

Regarding the third term, we will observe multiple cancellation effects. The goal is to isolate additional $x$-derivatives, which we then may redistribute to $\Phi$ via partial integration.
Further distributing the derivatives, we have 
\begin{mathex}[LCL]
	\mathcal{I}_3
	&=& \left| \int_\Rn t\Phi\left(\partial_t\left( F\partial_t\tilde{\mathcal{R}}_j[G] +\partial_t\tilde{\mathcal{R}}_j[F]G \right)
	+ \partial_t \left(\partial_t F\tilde{\mathcal{R}}_j[G] + \tilde{\mathcal{R}}_j[F]\partial_t G\right)\right)  \right|\\ 
	&\eqqcolon& \left| \int_\Rn t\Phi(\mathcal{J}_1+\mathcal{J}_2)\right|
\end{mathex}
After applying Lemma \ref{lem:Rt comms P1t/Rt interaction}, we observe the first cancellation effect;
\begin{equation*}
	\mathcal{J}_1 = \partial_t\left( F\,\partial_{x_j}G +\partial_{x_j}F\, G \right) = \partial_t\partial_{x_j}(FG) = \partial_{x_j}\partial_{t}(FG).
\end{equation*}
For the second integrand, with $\partial_{tt} F = -\Delta_x F$ and again Lemma \ref{lem:Rt comms P1t/Rt interaction} we obtain
\begin{mathex}
	\mathcal{J}_2 &=&	-\Delta_x F\,\tilde{\mathcal{R}}_j[G] - \tilde{\mathcal{R}}_j[F]\,\Delta_x G 
						+\partial_t F\, \partial_{x_j}G + \partial_{x_j}F\, \partial_t G\\
	&=&	-\Delta_x F\,\tilde{\mathcal{R}}_j[G] - \tilde{\mathcal{R}}_j[F]\,\Delta_x G 
		+\partial_{x_j}(\partial_tF\,G+ F\,\partial_tG) + \Delta_x\tilde{\mathcal{R}}_j[F]\,G + F\,\tilde{\mathcal{R}}_j[G].
\end{mathex}
Thanks to the cancellation of the term $\nabla_x F\cdot \nabla_x \tilde{\mathcal{R}}_j[G] + \nabla_x\tilde{\mathcal{R}}_j[F]\cdot \nabla_x G$ we have
\begin{equation*}
	\mathcal{J}_2 =	-\nabla_x\cdot\left(\nabla_x F\,\tilde{\mathcal{R}}_j[G]+ \tilde{\mathcal{R}}_j[F]\,\nabla_x G\right)
					+\nabla_x\cdot\left(\nabla_x \tilde{\mathcal{R}}_j[F]\,G+ F\,\nabla_x \tilde{\mathcal{R}}_j[G]\right)
					+\partial_{x_j}(\partial_tF\, G+ F\,\partial_tG).
\end{equation*}
Combining the calculations of $\mathcal{J}_1$ and $\mathcal{J}_2$,
\begin{mathex}
	\mathcal{J}_1 + \mathcal{J}_2 &=& \partial_{x_j}\partial_{t}(F\,G + \partial_tF\, G+ F\,\partial_tG)\\
					&& +\nabla_x\cdot\left(-\nabla_x F\,\tilde{\mathcal{R}}_j[G]- \tilde{\mathcal{R}}_j[F]\,\nabla_x G
					 					 +\nabla_x \tilde{\mathcal{R}}_j[F]\,G+ F\,\nabla_x \tilde{\mathcal{R}}_j[G]\right).
\end{mathex}
Plugging this into $\mathcal{I}_3$, via integration by parts in $x$-direction we can redistribute $\partial_{x_j}$ and $\nabla_x$ to $\Phi$ since the boundary terms disappear due to decay of the harmonic extensions, confer Corollary \ref{cor:P1t decay of P1t and derviatives}. 
Thus, we estimate $\mathcal{I}_3$ as in (\ref{eq:Rt comms sutiable form for trace theorems}).
\end{proof}

\subsection{1-dimensional $L^1$-estimate for a double-commutator}\label{sec:1d double comm}
In the previous subsection, we were not able to obtain an $L^1$-estimate for the commutator $[\Rj,\phi](f)$. 
This provides us with an opportunity to showcase an application for the interaction between harmonic extensions and fractional Laplacians. In this subsection, we will prove a replacement estimate, where instead of $[\Rj,\phi](f)$ itself we estimate the commutators of these commutators. 
We only obtain this replacement estimate in the one-dimensional case $n=1$ though, since we will us some properties that are exclusive to the $1$-dimensional Riesz transform $\mathcal{R}_1$.
This $1$-dimensional Riesz transform is called the Hilbert transform, for which we write $\mathcal{H}$. 

Regarding the estimates in Theorem \ref{th:Rt comms CRW theorem}, the proof for an $L^1$-version of (\ref{eq:Rt comms BMO}) and (\ref{eq:Rt comms Hölder}) already fails since the trace theorems from Section \ref{sec:tt} only allow to estimate against $\|f\|_{L^p(\Rn)}$ for $p>1$. 
This is not a concern for (\ref{eq:Rt comms fL&Rp}) though, since we could still choose $p_1,p_2\in(1,\infty)$.
However, we also do not have the $L^\infty$-boundedness of the Riesz transforms, which we would need to obtain $\|g\|_{L^\infty}$ from (\ref{eq:Rt comms sutiable form for trace theorems}).
Therefore, while the theorem below at first glance might look like an easy consequence from (\ref{eq:Rt comms fL&Rp}), with the proof in Subsection \ref{sec:Rt comms} we can not obtain a result which covers the estimates below.

\begin{theorem}\label{th:1d double comm estimate}
Assume $\phi,f\in\C_c^\infty(\R)$.
Suppose $s_1,s_2\in(0,1)$ with $s_1+s_2=1$ and $p,p'\in(1,\infty)$, $q,q'\in[1,\infty]$ with $\frac{1}{p}+\frac{1}{p'}=\frac{1}{q}+\frac{1}{q'}=1$.
Then
\begin{equation}\label{eq:1d double comm est}
	\|[\Hil,\phi](\fL{1}f)-[\Hil,f](\fL{1}\phi) \|_{L^1(\R)}
	\lesssim
	\|\fL{s_1}\phi\|_{L^{(p,q)}(\R)}\, \|\fL{s_2}f\|_{L^{(p',q')}(\R)}
\end{equation}
and
\begin{equation}\label{eq:1d double comm est of H}
	\left\| \Hil\left([\Hil,\phi](\fL{1}f) + [\Hil,f](\fL{1}\phi)\right) \right\|_{L^1(\R)}
	\lesssim
	\|\fL{s_1}\phi\|_{L^{(p,q)}(\R)}\, \|\fL{s_2}f\|_{L^{(p',q')}(\R)}.
\end{equation}
\end{theorem}

Just as for the Riesz transforms, given the harmonic extension $F(x,t)= P^1_t f(x)$ of $f\in\C_c^\infty(\R)$ we denote $\tilde{\mathcal{H}}F(x,t) = P^1_t \mathcal{H}f (x)$.
Before the main proof, we collect some additional facts about the Hilbert transform.

\begin{lemma}\label{lem:1d comm hilbert transform properties}
Let $f\in\CicRn$.
Then
\begin{equation*}
	\mathcal{H}\mathcal{H}f = -f.
\end{equation*}
Denote with $F(x,t)=P^1_tf$ the harmonic extension of f.
Then 
\begin{equation*}
	F_t \equiv \partial_tF = -\tilde{\mathcal{H}}F_x.
\end{equation*}
These two properties specifically hold for $\mathcal{H}$. 
Of course, as in Lemma \ref{lem:Rt comms P1t/Rt interaction} we also have
\begin{equation*}
	\tilde{\mathcal{H}}F_t = F_x.
\end{equation*}
\end{lemma}

\begin{proof}[Lemma \ref{lem:1d comm hilbert transform properties}]
The first equation is due to the Fourier symbol of $\mathcal{H}$.
With Proposition \ref{prop:Rt symbol} we have
\begin{equation*}
	\Fourier(\mathcal{H}\mathcal{H}f)(\xi) 
	= \left(-i\frac{\xi}{|\xi|}\right)\left(-i\frac{\xi}{|\xi|}\right)\Fourier(f)(\xi)
	= -\frac{\xi^2}{|\xi|^2}\Fourier(f)(\xi)
	= \Fourier(-f)(\xi).
\end{equation*}
Regarding the second equation, due to Lemma \ref{lem:P1t derivatives of P1t 2} we have 
\begin{equation*}
	F_t(x,t) = -P^1_t(\fL{1}f)(x)
\end{equation*}
as well as 
\begin{equation*}
	\tilde{\mathcal{H}}F_x(x,t) = P^1_t(\mathcal{H}\partial_x f)(x) = P^1_t(\fL{1}f)(x).
\end{equation*}
The last equation can easily be verified on the side of the Fourier transforms;
\begin{equation*}
	 \Fourier(\mathcal{H}\partial_x f)
	 = \left(-i\frac{\xi}{|\xi|}\right)(2\pi i \xi)\Fourier(f)(\xi)
	 = 2\pi \frac{\xi^2}{|\xi|}\Fourier(f)(\xi)
	 = \Fourier(\fL{1}f)(\xi).
\end{equation*}
As already mentioned above, see Lemma \ref{lem:Rt comms P1t/Rt interaction} for the last equation.
\end{proof}

\begin{proof}[Estimate (\ref{eq:1d double comm est})]
Let $g\in\CicRn$. 
Redistributing the Hilbert transform to $g$ with Lemma \ref{lem:Rt integration by parts}, we obtain 
\begin{mathex}[LCL]
	\mathcal{I}
	&\coloneqq& \left|\int_\Rn \left([\Hil,\phi](\fL{1}f)-[\Hil,f](\fL{1}\phi)\right)g\,\right| \\
	&=&		\left|\int_\Rn \phi\,\Hil\fL{1}f\,g + \phi\,\fL{1}f\,\Hil g - \Hil\fL{1}\phi\,f\,g - \fL{1}\phi\,f\,\Hil g\,\right|.
\end{mathex}
In order to prove (\ref{eq:1d double comm est}), by $L^1$-duality it suffices to show that 
\begin{equation}\label{eq:1d comm duality estimate}
	\mathcal{I}	\lesssim
	\|\fL{s_1}\phi\|_{L^{(p,q)}(\R)}\,\|\fL{s_2}f\|_{L^{(p',q')}(\R)}\,\|g\|_{L^\infty(\R)}.
\end{equation}
As usual, let $\Phi(x,t)=P^1_t \phi(x)$, $F(x,t)=P^1_t f(x)$, $G(x,t)=P^1_t g(x)$ be the harmonic extensions of $\phi,f,g$.
Now, regarding the integration by parts, the situation differs a little bit from what we have grown used to. 
There are fractional Laplacians in the boundary term. 
Recalling Lemma \ref{lem:P1t derivatives of P1t 2} though, these can be viewed as the limits of the harmonic extension's $t$-derivatives.
Therefore, via integration by parts in $t$ direction,
\begin{equation*}\label{eq:1d comm integration by parts}
	\mathcal{I} = \left|\int_{\R_+^2} \partial_t\left(
								\Phi\,\tilde{\Hil}F_t\, G + \Phi \,F_t\,\tilde{\Hil}G 
								- \tilde{\Hil}\Phi_t\, F \,G - \Phi_t \,F \,\tilde{\Hil} G
				\right)\right|
				\eqqcolon \left|\int_{\R_+^2} \mathcal{J}\,\right|.
\end{equation*}
Of course, the other boundary term vanishes due to Corollary \ref{cor:P1t decay of P1t and derviatives}.
We claim that 
\begin{equation}\label{eq:1d comm suitable form for trace theorem}
	\mathcal{I} \lesssim\int_{\R^2_+}|\nabla_{\R^2}\Phi| |\nabla_{\R^2}F| | G|
				=		\int_{\R^2_+}t^{1-s_1-s_2}|\nabla_{\R^2}\Phi| |\nabla_{\R^2}F| | G|.
\end{equation}
Then (\ref{eq:1d double comm est}) follows with Proposition \ref{prop:tt - Lp est}.
Therefore, during the distribution of the additional derivative, the goal is to replace the occurring Hilbert transforms via Lemma \ref{lem:1d comm hilbert transform properties} whenever possible, hoping that the remaining terms with Hilbert transforms cancel each other out.
Applying Lemma \ref{lem:1d comm hilbert transform properties}, we obtain
\begin{mathex}
	\mathcal{J} 
	&=& \partial_t\left(\Phi\, F_x \,G + \Phi\, F_t\,\tilde{\Hil}G - \Phi_x\, F\, G - \Phi_t \,F \,\tilde{\Hil} G \right)\\
	&=& \partial_t\left(\Phi\, F_x \,G - \Phi_x\, F\, G\right) 
		+ \partial_t\left(\Phi\, F_t\,\tilde{\Hil}G - \Phi_t \,F \,\tilde{\Hil} G \right).
\end{mathex}
Distributing the derivative, we observe the cancellation of the term $\Phi_t\, F_t\,\tilde{\Hil}G$ in the second group.
\begin{mathex}
	\mathcal{J} 
	&=& \left(\Phi_t\, F_x \,G - \Phi_x\, F_t\, G\right) + \left(\Phi\, F_{xt} \,G - \Phi_{xt}\, F\, G\right) 
		+ \left(\Phi\, F_x \,G_t - \Phi_x\, F\, G_t\right) \\
	&&	+ \left(\Phi\, F_{tt}\,\tilde{\Hil}G - \Phi_{tt} \,F \,\tilde{\Hil} G \right) 
		+ \left(\Phi\, F_t\,\tilde{\Hil}G_t - \Phi_t \,F \,\tilde{\Hil}G_t \right)
\end{mathex}
Applying Lemma \ref{lem:1d comm hilbert transform properties} as well as the harmonicity of the extensions, i.e. $\partial_{tt}F= -\partial_{xx}F$, we have
\begin{mathex}
	\mathcal{J} 
	&=& \left(\Phi_t\, F_x \,G - \Phi_x\, F_t\, G\right) + \left(\Phi\, F_{xt} \,G - \Phi_{xt}\, F\, G\right) 
		+ \left(\Phi\, F_x \,G_t - \Phi_x\, F\, G_t\right) \\
	&&	+ \left(-\Phi\, F_{xx}\,\tilde{\Hil}G + \Phi_{xx} \,F \,\tilde{\Hil} G \right) 
		+ \left(\Phi\, F_t\,G_x - \Phi_t \,F \,G_x \right).
\end{mathex}
We can rewrite the terms of the second to last group as 
\begin{mathex}[RCL]
	 -\Phi\,F_{xx}\,\tilde{\Hil}G 	&=& -(\Phi\,F_x\,\tilde{\Hil}G)_x + \Phi_x\,F_x\,\tilde{\Hil}G - \Phi\,F_x\,G_t \\
	  \Phi_{xx}\,F\,\tilde{\Hil}G	&=& (\Phi_x \,F \,\tilde{\Hil}G)_x - \Phi_x \,F_x \,\tilde{\Hil}G +\Phi_{x}\,F\,G_t
\end{mathex}
where for the last term we applied Lemma \ref{lem:1d comm hilbert transform properties} respectively.
Plugging these terms back into $\mathcal{J}$, we observe the cancellations of $\Phi_x\,F_x\,\tilde{\Hil}G$ as well as $\Phi\,F_x\,G_t$ and $\Phi_{x}\,F\,G_t$;
\begin{mathex}
	\mathcal{J} 
	&=& \left(\Phi_t\, F_x \,G - \Phi_x\, F_t\, G\right) + \left(\Phi\, F_{xt} \,G - \Phi_{xt}\, F\, G\right) \\
	&&	+ \left(-(\Phi\,F_x\,\tilde{\Hil}G)_x + (\Phi_x \,F \,\tilde{\Hil}G)_x \right) 
		+ \left(\Phi\, F_t\,G_x - \Phi_t \,F \,G_x \right).
\end{mathex}
Repeating this reorganization for the other group with double derivatives on one function,
\begin{mathex}[RCL]
	  \Phi\,F_{xt}\,G	&=& (\Phi\,F_t\,G)_x - \Phi_x\,F_t\,G - \Phi\, F_t \,G_x\\
	- \Phi_{xt}\,F\,G	&=& -(\Phi_t\,F\,G)_x + \Phi_t\,F_x\,G + \Phi_t\,F\,G_x.
\end{mathex}
Thus, with the cancellation of $\Phi\, F_t \,G_x$ and $\Phi_t\,F\,G_x$,
\begin{equation*}
	\mathcal{J} = 2\left(\Phi_t\, F_x \,G - \Phi_x\, F_t\, G\right) 
	+ \left(\Phi\,F_t\,G - \Phi_x \,F \,\tilde{\Hil}G - \Phi\,F_x\,\tilde{\Hil}G + \Phi_x \,F \,\tilde{\Hil}G \right)_x.
\end{equation*}
The second term vanishes when integrating in $x$ due to the decay estimates for harmonic extensions and Riesz transforms, see Corollary \ref{cor:P1t decay of P1t and derviatives} and Lemma \ref{lem:Rt decay}.
Therefore, when we plug $\mathcal{J}$ back into (\ref{eq:1d comm integration by parts}), we get
\begin{equation*}
	\mathcal{I}= 2\left|\int_{\R^2_+}\Phi_t\, F_x \,G - \Phi_x\, F_t\, G \,\right| 
				= 2\left|\int_{\R^2_+}\det\left(\nabla_{\R^2}\Phi,\nabla_{\R^2}F\right)G \,\right|.
\end{equation*}
Concluding this proof, we can estimate $\mathcal{I}$ as in (\ref{eq:1d comm suitable form for trace theorem}) and then obtain (\ref{eq:1d comm duality estimate}) and therefore (\ref{eq:1d double comm est}) with Proposition \ref{prop:tt - Lp est}(a,b).
\end{proof}

\begin{proof}[Estimate (\ref{eq:1d double comm est of H})]
Let $g\in\C^\infty_c(\R)$.
Applying the double commutator to $g$, with Lemma \ref{lem:Rt integration by parts} and $\Hil \Hil= -Id$ we obtain
\begin{mathex}[LCL]
		\mathcal{I}	
		&\coloneqq& \left|\int_\R \Hil\left([\Hil,\phi](\fL{1}f) + [\Hil,f](\fL{1}\phi)\right)g \,\right|\\
		&=& \left|\int_\R \left([\Hil,\phi](\fL{1}f) + [\Hil,f](\fL{1}\phi)\right)\Hil g \,\right|\\
		&=& \left|\int_\R  \phi\,\Hil\fL{1}f\,\Hil g - \phi\,\fL{1}f\, g + \Hil\fL{1}\phi\,f\,\Hil g - \fL{1}\phi\,f\,g  \,\right|.
\end{mathex}
In order to prove (\ref{eq:1d double comm est of H}), again by duality, it suffices to show that
\begin{equation}\label{eq:1d comms with Hil duality est}
	\mathcal{I}	\lesssim	\|\fL{s_1}\phi\|_{L^{(p,q)}(\R)}\,\|\fL{s_2}f\|_{L^{(p',q')}(\R)}\,\|g\|_{L^\infty(\R)}.
\end{equation}
Denote with $\Phi(x,t)=P^1_t\phi(x)$, $F(x,t)=P^1_t f(x)$, $G(x,t)=P^1_t g(x)$ the harmonic extensions.
As in the proof of (\ref{eq:1d double comm est}), with integration by parts in $t$ and due to the interaction of fractional Laplacian and $t$-derivatives,
\begin{equation*}
	\mathcal{I} = 	\left|\int_{\R^2_+} \partial_t\left( 
							\Phi\,\tilde{\Hil}F_t\,\tilde{\Hil}G - \Phi\,F_t\, G 
							+ \tilde{\Hil}\Phi_t\,F\,\tilde{\Hil}G - \Phi_t\,F\,G 
					\right) \right|
				\eqqcolon \left|\int_{\R^2_+} \mathcal{J}\right|.
\end{equation*}
Before distributing the additional $t$-derivative, with Lemma \ref{lem:1d comm hilbert transform properties} we obtain
\begin{equation*}
	\mathcal{J} = \partial_t\left(\Phi\,F_x\,\tilde{\Hil}G - \Phi\,F_t\, G + \Phi_x\,F\,\tilde{\Hil}G - \Phi_t\,F\,G\right).
\end{equation*}
Again, we try to write multiple terms as single $x$-derivative, which will vanish when plugging $\mathcal{J}$ into $\mathcal{I}$. With Lemma \ref{lem:1d comm hilbert transform properties}, more specifically $\tilde{\Hil}G_x= -G_t$, 
\begin{mathex}
	\mathcal{J} &=& \partial_t \left((\Phi\,F)_x\,\tilde{\Hil}G\right) - \partial_t\left((\Phi\,F)_t\, G \right) \\
				&=& \left(\Phi\,F\,\tilde{\Hil}G\right)_{xt} 
					- \partial_t\left(-(\Phi\,F)G_t\right)- \partial_t\left((\Phi\,F)_t\, G \right).
\end{mathex}
Due to the cancellation of the term $(F\,G)_t\,\Phi_t$ and then with harmonicity,
\begin{mathex}
	\mathcal{J} &= \left(\Phi\,F\,\tilde{\Hil}G\right)_{xt}  + (\Phi\,F)\,G_{tt}- (\Phi\,F)_{tt}\,G\\
				&= \left(\Phi\,F\,\tilde{\Hil}G\right)_{xt}  - (\Phi\,F)\,G_{xx}- (\Phi\,F)_{tt}\,G.
\end{mathex}
As mentioned above, the first term vanishes when integrating in $x$ due to the decay estimates for the harmonic extension and the Riesz transforms, see Corollary \ref{cor:P1t decay of P1t and derviatives} and Lemma \ref{lem:Rt decay}.
Via a double integration by parts in $x$-direction we obtain
\begin{equation*}
	\mathcal{I} = \left| \int_{\R^2_+} (\Phi\,F)\,G_{xx} + (\Phi\,F)_{tt}\,G \,\right| 
				= \left| \int_{\R^2_+} (\Phi\,F)_{xx}\,G + (\Phi\,F)_{tt}\,G \,\right| 
				= 2\left| \int_{\R^2_+} \nabla_{\R^2}\Phi\cdot \nabla_{\R^2}F \,G \,\right| 
\end{equation*}
since $(\Phi\,F)_{xx} + (\Phi\,F)_{tt} = 2\nabla_{\R^2}\Phi\cdot \nabla_{\R^2}F$ due to harmonicity. 
Therefore, we have
\begin{equation*}
	\mathcal{I}\lesssim \int_{\R^2_+} |\nabla_{\R^2}\Phi|\,|\nabla_{\R^2}F| \,|G| 
						= \int_{\R^2_+} t^{1-s_1-s_2}|\nabla_{\R^2}\Phi|\,|\nabla_{\R^2}F| \,|G|
\end{equation*}
and thus obtain (\ref{eq:1d comms with Hil duality est}) with Proposition \ref{prop:tt - Lp est}(a,b).
\end{proof}

\newpage
\section{Commutator estimates via $s$-harmonic extension}\label{sec:s harmonic extension}
In Subsection \ref{sec:1d double comm}, we just used the harmonic extension's behavior near the boundary, that is $\lim_{t\rightarrow 0} F^1(x,t) = -\fL{1}f(x)$, to show estimates for a commutator involving the half-Laplacian $\fL{1}$. 
Modifying the extension such that we instead obtain $\fL{s}f$ as boundary term, might allow us to deal with commutators involving arbitrary fractional Laplacians.
In the first subsection, we introduce the so-called $s$-harmonic extensions, which do just that.
Most of the results from Subsection \ref{sec:P1t} can be transferred to the $s$-harmonic extensions. 
Equipped with this new tool, we then prove estimates for four different commutators involving fractional Laplacians and Riesz potentials.

\subsection{$s$-harmonic extension to $\Rnpp$ via the generalized Poisson operator}\label{sec:Pst}
We obtain the $s$-harmonic extension to $\Rnpp$ via a generalized Poisson operator, which is defined as follows.

\begin{definition}[The generalized Poisson operator]\label{def:Pst generalized Poisson operator}
Let $s\in(0,2)$.
The generalized Poisson extension operator $P^s_t$ is given via the convolution
\begin{equation*}
	P^s_t f(x)	\coloneqq C_{n,s} \int_{\Rn} \frac{t^s}{\left(|x-y|^2+t^2\right)^\frac{n+s}{2}} f(y) \,dy
				= C_{n,s} (p^s_t * f) (x)
\end{equation*}
where $f\in L^1(\Rn)+L^\infty(\Rn)$ and the kernel $p^s_t$ is given by 
\begin{equation*}
	p^s_t (x) \coloneqq \frac{t^s}{\left(|x|^2+t^2\right)^\frac{n+s}{2}}.
\end{equation*}
\end{definition} 

The function $F^s(x,t)\coloneqq P^s_t f(x)$ is called the $s$-harmonic extension of $f$ to $\Rnpp$ and, for example for $f\in L^\infty(\Rn)\cap\C(\Rn)$, satisfies the Dirichlet problem
\begin{equation}\label{eq:Pst pde}
	\left\lbrace \begin{array}{llll}
		 \operatorname{div}_{\Rnp}(t^{1-s}\nabla_{\Rnp} F^s(x,t))	&=& 0	\quad&	\text{in }\Rnpp, \\
		 \lim_{t\rightarrow 0} F^s(x,t) 						&=& f(x)	\quad&	\text{on }\Rn,\\
		 \lim_{|(x,t)|\rightarrow\infty} F^s(x,t)				&=&	0.	\quad&
	\end{array}\right.
\end{equation}
As suggested above, the behavior of $\partial_t F^s$ towards the boundary is of central interest to us. 
Caffarelli and Silvestre showed in \cite{CaSi07} that 
\begin{equation}\label{eq:Pst partial t on boundary}
	\lim_{t\rightarrow 0} -t^{1-s}\partial_t F^s(x,t) 		= c\fL{s}f(x)	\quad	\text{on }\Rn
\end{equation}
for a fixed constant $c$ depending on $n$ and $s$.

We further observe that, as for the classical harmonic extension, the $p^s_t$ are dilations of $p^s_1$, i.e. $p^s_t(x)= t^{-n}p^s_1(t^{-1}x) = (p^s_1)_t(x)$. 
Therefore, the boundary condition in (\ref{eq:Pst pde}) is indeed satisfied according to Example 1.2.17 and Theorem 1.2.19 in \cite{Gra14} if $C_{n,s}$ is chosen correctly.
In fact, we have
\begin{equation*}
	\|p^s_t\|_{L^1(\Rn)} = \int_\Rn p^s_t = \int_\Rn p^s_1 \eqqcolon \frac{1}{C_{n,s}}.
\end{equation*} 
Just as we did with the classical Poisson kernel, we might also interpret the kernels $(p^s_t)_{t>0}$ as function $k^s$ on $\Rnpp$ with $k^s(x,t)\coloneqq p^s_t(x)$. 

These kernels $k^s$ are $s$-harmonic, meaning they satisfy the partial differential equation in (\ref{eq:Pst pde}),
\begin{mathex}
	\operatorname{div}\left(t^{1-s}\nabla k^s(x,t)\right) 
		&=&	t^{1-s}\Delta_x p^s_t(x) + \partial_t(t^{1-s}\partial_t k^s)(x,t) \\
		&=& 	t^{1-s}\left(
				(n+s)(n+s+2)\frac{|x|^2 t^s}{(|x|^2+t^2)^{\frac{n+s}{2}+2}} - (n+s)n\frac{t^s}{(|x|^2+t^2)^{\frac{n+s}{2}+1}}
			\right)\\
		& &	-(2+s)(n+s)\frac{|x|^2t}{(|x|^2+t^2)^{\frac{n+s}{2}+1}}+(n+s)n\frac{t^3}{(|x|^2+t^2)^{\frac{n+s}{2}+1}}\\
		&=&	0.
\end{mathex}
The extension via the generalized Poisson operator inherits this property. 
We can therefore confirm that $F^s$ satisfies the partial differential equation in (\ref{eq:Pst pde}).
The arguments for this inheritance are analogous to the ones in the proof of the following lemma, which we use to investigate the behavior of $F^s$ and its derivatives towards infinity in the same fashion as for the classical harmonic extension.

\begin{lemma}\label{lem:Pst derivatives of Pst 1}
Let $s\in(0,2)$, $k_i\in\Nz$ for $0\leq i \leq n$. Set $\tilde{k}\coloneqq k_1+\ldots+k_n$ and $k\coloneqq \tilde{k}+k_0$.
For any $f\in L^\infty(\Rn)$ with $F^s(x,t)\coloneqq P^s_t f(x)$ we have
\begin{equation*}
	\partial_t^{k_0}\partial_{x_1}^{k_1}\dots\partial_{x_n}^{k_n} F^s(x,t) 
	= C_{n,s} \int_\Rn
				\left(\sum_{j=0}^{k}
						\frac{ t^{s-k_0} q_{j}(x-y)}{(|x-y|^2+t^2)^{\frac{n+s+\tilde{k}+j}{2}}}
				\right) f(y)\,dy,
\end{equation*} 
where the $q_{j}$ are homogeneous polynomials of degree $j$, which depend on $k_0,\dots,k_n$. 
In particular, we have $F^s\in\C^\infty(\Rnpp)$.
\end{lemma}

\begin{proof}
We first show that for the derivatives of the kernel we have 
\begin{equation}\label{eq:Pst pst derivatives}
	\partial_t^{k_0}\partial_{x_1}^{k_1}\dots\partial_{x_n}^{k_n} k^s(x,t) 
	= \sum_{j=0}^{k}		\frac{ t^{s-k_0} q_{j}(x)}{(|x|^2+t^2)^{\frac{n+s+\tilde{k}+j}{2}}}.
\end{equation}
Then the lemma will follow with Lebesgue's dominated convergence theorem and via induction.

To this end, we first investigate the interaction between a dilation $f_t(x)=t^{-n}f(t^{-1}x)$ and the derivatives for a sufficiently differentiable function $f$.
Regarding derivatives in $x$-direction, we have 
\begin{equation*}
	\nabla_x^k f_t(x) = \nabla_x^k \left((x,t)\mapsto t^{-n} f(t^{-1}x)\right)(x) 
	= t^{-k}\left(\nabla_x^kf\right)_t(x)
\end{equation*}
Regarding the derivatives in direction of the the dilation parameter,
\begin{equation*}
	\frac{d}{dt} t^{-k}f_t(x) = -t^{-n-k-1} \nabla f (t^{-1}x)\cdot(t^{-1}x) -(n+k)t^{-n-k-1} f(t^{-1}x).
\end{equation*}
Therefore, via induction we obtain 
\begin{equation*}
	\partial_t^{k_0} \left((x,t)\mapsto t^{-k}f_t(x)\right) 
	= t^{-k-n-l}\sum_{|\alpha|\leq k_0} c_\alpha \partial_{x_1}^{\alpha_1}\ldots\partial_{x_n}^{\alpha_n}f(t^{-1}x)(t^{-1}x)^\alpha
\end{equation*}
for some constants $c_\alpha\in\R$ depending on $k_0$.
Specifically for the kernel $p^s_1$, we see that 
\begin{equation*}
	\partial_{x_1}^{k_1}\dots\partial_{x_n}^{k_n} p^s_1(x) = \sum_{i=0}^{\tilde{k}} \frac{\tilde{q}_i(x)}{(1+|x|^2)^\frac{n+s+\tilde{k}+i}{2}}
\end{equation*}
where the $\tilde{q}_i$ are homogeneous polynomials of degree $i$ depending on $k_1,\ldots,k_n$.
Combining these three observations, we obtain
\begin{mathex}
	\partial_t^{k_0}\partial_{x_1}^{k_1}\dots\partial_{x_n}^{k_n} k^s(x,t)
	&=& t^{-k-n-l}\sum_{|\alpha|\leq k_0} c_\alpha (t^{-1}x)^\alpha \partial_{x_1}^{k_1+\alpha_1}\ldots\partial_{x_n}^{k_n+\alpha_n}p_s^1(t^{-1}x)\\
	&=& t^{-k-n-l}\sum_{|\alpha|\leq k_0} \sum_{i=0}^{\tilde{k}+|\alpha|} \frac{(t^{-1}x)^\alpha\tilde{q}_{\alpha,i}(t^{-1}x)}{(1+|x|^2)^\frac{n+s+\tilde{k}+|\alpha|+i}{2}}\\
	&=& t^{-k-n-l}\sum_{j=0}^{\tilde{k}+k_0} \frac{q_{j}(t^{-1}x)}{(1+|x|^2)^\frac{n+s+\tilde{k}+j}{2}}
\end{mathex}
where the $q_j$ are homogeneous polynomials of degree $j$ depending on $k_0,\ldots,k_n$.
Thus, we have shown (\ref{eq:Pst pst derivatives}).

The lemma now follows via induction and Lebesgue's dominated convergence theorem where the majorants can be chosen as follows. 
For the derivatives in $x_l$-direction with $0<h\leq 1$ we have
\begin{mathex}[RRL]
	&|D_{h,x_l}&(\partial_t^{k_0}\partial_{x_1}^{k_1}\dots\partial_{x_n}^{k_n} (k^s(x-y,t) f(y)))| \\
	&=& |f(y)|\left|\frac{\partial_t^{k_0}\partial_{x_1}^{k_1}\dots\partial_{x_n}^{k_n} k^s(x+h e_l-y,t)-\partial_t^{k_0}\partial_{x_1}^{k_1}\dots\partial_{x_n}^{k_n} k^s(x-y,t)}{h}\right|\\
	&\leq& \|f\|_{L^\infty(\Rn)} \sup_{\tilde{x}\in B_1(x)}
					|\partial_{x_l}(\partial_t^{k_0}\partial_{x_1}^{k_1}\dots\partial_{x_n}^{k_n} k^s)(\tilde{x}-y,t)|\\
	&\leq& 
	\begin{cases}
		C\|f\|_{L^\infty(\Rn)}t^{s-k_0}((|x-y|-1)^2+t^2)^{-\frac{n+s+(\tilde{k}+1)}{2}}
				\quad&\text{for } |x-y|>2\\
		\|\partial_{x_l}(\partial_t^{k_0}\partial_{x_1}^{k_1}\dots\partial_{x_n}^{k_n} k^s)(\cdot,t)\|_{L^\infty(\Rn)} 
				\quad& \text{for } |x-y|\leq 2
	\end{cases}
\end{mathex}
where we first use the mean value theorem and then that $|x|^{-j} q_j(x)\leq C$ for each $0 \leq j \leq k$ due to the homogeneity of $q_j$.
For the derivative in $t$-direction we have
\begin{mathex}[RRL]
	&|D_{h,t}&(\partial_t^{k_0}\partial_{x_1}^{k_1}\dots\partial_{x_n}^{k_n} (k^s(x-y,t) f(y)))| \\
	&=& |f(y)|\left|\frac{\partial_t^{k_0}\partial_{x_1}^{k_1}\dots\partial_{x_n}^{k_n} k^s(x-y,t+h)-\partial_t^{k_0}\partial_{x_1}^{k_1}\dots\partial_{x_n}^{k_n} k^s(x-y,t)}{h}\right|\\
	&\leq & \|f\|_{L^\infty(\Rn)} \sup_{\tilde{t}\in\left[\frac{t}{2},2t\right]} 
					|\partial_t(\partial_t^{k_0}\partial_{x_1}^{k_1}\dots\partial_{x_n}^{k_n} k^s)(x-y,\tilde{t})|\\
	&\leq & C\|f\|_{L^\infty(\Rn)}(2t)^{s}\left(\frac{t}{2}\right)^{-k_0} \left((|x-y|^2+\left(\frac{t}{2}\right)^2\right)^{-\frac{n+s+(\tilde{k}+1)}{2}}
\end{mathex}
for $h\in\left[-\frac{t}{2},t\right]$ thanks to the mean value theorem.
It is easy to see that these majorants are integrable.
\end{proof}

As a consequence of this lemma, we get the following result for the decay of the $s$-harmonic extension and its derivatives with the same arguments as in Corollary \ref{cor:P1t decay of P1t and derviatives} for the classical harmonic extension.

\begin{corollary}[Decay estimate for the $s$-harmonic extension]\label{cor:Pst decay of Pst and derviatives}
Let $s\in(0,2)$ $f\in L^\infty(\Rn)$ and $k\in\Nz$.
Denote with $F^s(x,t)\coloneqq P^s_t f(x)$ the $s$-harmonic extension.
Then 
\begin{equation}\label{eq:Pst t-decay of Pst for Linfty}
	\sup_{(x,t)\in\Rnpp} t^k|\nabla_\Rnpp^k F^s(x,t)| \lesssim \|f\|_{L^\infty(\Rn)}.
\end{equation}
If additionally $f\in L^1(\Rn)$, then
\begin{equation}\label{eq:Pst t-decay of Pst for L1}
	\sup_{(x,t)\in\Rnpp} t^{n+k}|\nabla_\Rnpp^k F^s(x,t)| \lesssim \|f\|_{L^1(\Rn)}.
\end{equation}
Regarding the decay in $x$-direction, if further for some $C>0$, $k>0$ we have $|f(x)|\leq C|x|^{-l}$, then
\begin{equation}\label{eq:Pst x-decay of Pst}
	t^{k}|\nabla_\Rnpp^k F^s(x,t)| \lesssim |x|^{-l} + t^s\|f\|_{L^1(\Rn)}|x|^{-n-s}.
\end{equation}
\end{corollary}

We conclude this section with a short collection of some more results, which are analogous to those for the classical Poisson operator.

\begin{proposition}[Fourier transform of the generalized Poisson Kernels]\label{prop:Pst kernel fourier symbol}
Let $s\in(0,2)$. The Fourier transform of the Poisson kernel is 
\begin{equation*}
	\Fourier(p^s_t)(\xi) = c_{n,s} \int_0^\infty \lambda^\frac{s}{2}e^{-\lambda-\frac{|t\xi|^2}{c\lambda}}\frac{d\lambda}{\lambda}.
\end{equation*}
for some constants $c_{n,s}\neq 0$, $c>0$.
\end{proposition}

\begin{proof}
Proposition 7.6 in \cite{Hao16} yields
\begin{equation}\label{eq:Pst fourier transform of ps1}
	\Fourier(p^s_1)(\xi) = (2\pi)^n c_{n,s} \int_0^\infty \lambda^\frac{s}{2}e^{-\lambda-\frac{|2\pi\xi|^2}{4\lambda}}\frac{d\lambda}{\lambda}.
\end{equation}
The result follows immediately since $\Fourier(p^s_t)= t^{-n} \Fourier(\delta^{t^{-1}}p^s_1) = \delta^t\Fourier(p^s_1)$ thanks to elementary properties of the Fourier transform. Here $\delta^t$ denotes $\delta^t g(x) = g(tx)$ for any $t>0$ and measurable functions $g$ on $\Rn$.
\end{proof}

In comparison to the classical Poisson kernel, the kernel of the generalized Poisson operator has a way more complex Fourier transform, which makes it more difficult to quantify the derivatives of the $s$-harmonic extension in $t$-direction. 
For our purposes, due to (\ref{eq:Pst pde}) and (\ref{eq:Pst partial t on boundary}) we already know enough regarding the $t$-derivatives of the extension.
Therefore, we settle for a weaker result compared to the corresponding result for the harmonic extension, Lemma \ref{lem:P1t derivatives of P1t 2}.

\begin{lemma}\label{lem:Pst derivatives of Pst 2}
Let $s\in(0,2)$, $k_i\in\Nz$ for $1\leq i \leq n$ and $f\in\CicRn$.
Set $F^s(x,t)\coloneqq P^s_t f(x)$.
Then 
\begin{equation*}
	\partial_{x_1}^{k_1}\dots\partial_{x_n}^{k_n} F^s(x,t) 
	= P_t^s (\partial_{x_1}^{k_1}\dots\partial_{x_n}^{k_n} f)(x).
\end{equation*}
\end{lemma} 

\begin{proof}
Analogous to Lemma \ref{lem:P1t derivatives of P1t 2}.
\end{proof}

As for the classical harmonic extension, we easily obtain the boundedness of the derivatives in $x$-direction from this lemma.

\begin{corollary}[$L^\infty$-estimate for the $s$-harmonic extenison]\label{cor:Pst boundedness of Pst derivatives}
Let $s\in(0,2)$, $k_i\in\Nz$ for $1\leq i \leq n$ and $f\in\CicRn$. Then 
\begin{equation*}
	\|\partial_{x_1}^{k_1}\dots\partial_{x_n}^{k_n} P_t^s f\|_{L^\infty(\Rn)} \leq 
	\|\partial_{x_1}^{k_1}\dots\partial_{x_n}^{k_n} f \|_{L^\infty(\Rn)}.
\end{equation*}
\end{corollary}

\begin{proof}
Analogous to Corollary \ref{cor:P1t boundedness of P1t derivatives}.
\end{proof}

\subsection{Estimating the commutator $[\fL{s},\phi]$}\label{sec:fL comms}
The motivation for introducing the generalized harmonic extensions was the ability to deal with arbitrary fractional Laplacians in the integral term we want to estimate.
Therefore, an obvious first application is the commutator of fractional Laplacians and pointwise multiplication.
For $s>0$ and $\phi\in\CicRn$ the commutator $[\fL{s},\phi]$ is defined by
\begin{equation*}
	[\fL{s},\phi](f) \coloneqq \fL{s}(\phi f) - \phi\, \fL{s} f, \quad \text{for } f\in\CicRn.
\end{equation*}

\begin{theorem}\label{th:fL comms}
Let $s\in(0,1)$ and $\phi,f\in\CicRn$.\\
Suppose $p\in(1,\infty)$ and $\sigma\in[s,1)$. Then
\begin{equation}\label{eq:fL comms Hölder}
	\| [\fL{s},\phi](f) \|_{L^p(\Rn)} \lesssim [\fL{\sigma}\phi]_{BMO}\|I^{\sigma-s}f\|_{L^p(\Rn)}.
\end{equation}
Suppose $q_1, q_2, p\in(1,\infty)$ with $\frac{1}{q_1}+\frac{1}{q_2}=\frac{1}{p}$ and $\sigma\in[s,1)$. Then
\begin{equation}\label{eq:fL comms fL&Rp}
	\| [\fL{s},\phi](f) \|_{L^p(\Rn)} \lesssim \|\fL{\sigma}\phi\|_{L^{q_1}(\Rn)} \|I^{\sigma-s}f\|_{L^{q_2}(\Rn)}.
\end{equation}
\end{theorem}

\begin{proof}
By $L^p$-duality, it suffices to show that for every $g\in\CicRn$
\begin{equation}\label{eq:fL comms duality estimates}
	 \mathcal{I}\coloneqq \left|\int_\Rn [\fL{s},\phi](f)\,g\,\right| 
	 \lesssim \left\lbrace \begin{array}{ll}
					{[}\fL{\sigma}\phi{]}_{BMO}\|I^{\sigma-s}f\|_{L^p(\Rn)}\|g\|_{L^{p'}(\Rn)},\\
					\|\fL{\sigma}\phi\|_{L^{q_1}(\Rn)} \|I^{\sigma-s}f\|_{L^{q_2}(\Rn)} \|g\|_{L^{p'}(\Rn)}.\\
				\end{array}	\right.
\end{equation}
With the integration by parts formula for the fractional Laplacian, see Lemma \ref{lem:fL integration by parts}, we have
\begin{equation*}
	\mathcal{I} = \left|\int_\Rn \fL{s}(\phi f)\,g - \phi\, \fL{s} f\,g\,\right|
				= \left|\int_\Rn \phi\,f\,\fL{s}g - \phi\, \fL{s} f\,g\,\right|.
\end{equation*}
Let $\Phi(x,t)=P^s_t \phi(x)$, $F(x,t)=P^s_tf(x)$, $G(x,t)=P^s_t f(x)$ be the $s$-harmonic extensions of $\phi$, $f$, $g$.
Then with (\ref{eq:Pst partial t on boundary}), up to a constant $c$ which we omit, we have
\begin{equation*}
	\phi\,f\,\fL{s}g - \phi\, \fL{s} f\,g 
	= \lim_{t\rightarrow 0} \Phi \left(F\,t^{1-s}\partial_t G - t^{1-s}\partial_t F\,G\right).
\end{equation*} 
With a partial integration in $t$, where the other boundary term disappears due to the decay of the $s$-harmonic ectensions, see Corollary \ref{cor:Pst decay of Pst and derviatives},
\begin{equation*}
	\mathcal{I} 
	= \left|\int_\Rnpp \partial_t\left(\Phi \left(F\,t^{1-s}\partial_t G - t^{1-s}\partial_t F\,G\right)\right) \right|
	\eqqcolon \left|\int_\Rnpp \mathcal{J}\, \right|.
\end{equation*}
For the rest of the proof, we will show that 
\begin{equation}\label{eq:fL comms suitable form for trace theorems}
	\begin{split}
	\mathcal{I}	\lesssim& \int_{\Rnpp} t^{1-s}|\nabla_x \Phi|\,(|\nabla_xF|\,|G| + |\nabla_xG|\,|F|)\\
				&+ \int_{\Rnpp} t^{2-s}|\nabla_x\nabla_{\Rnp} \Phi|\,(|\nabla_{\Rnp}F|\,|G| + |\nabla_{\Rnp}G|\,|F|).
	\end{split}
\end{equation}
From this estimate we obtain (\ref{eq:fL comms duality estimates}) by applying Proposition \ref{prop:tt - Lp est} and Proposition \ref{prop:tt - Hölder est}(c,d,e) respectively, sticking to the sequence of the functions suggested here for each summand.  
Note that for Proposition \ref{prop:tt - Hölder est}, in order to obtain $\|g\|_{L^{p'}(\Rn)}$ from $|G|$, we have to rely on variation (e) and therefore obtain no actual H\"older estimate.

Distributing the $t$-derivative, we observe the cancellation of the term $t^{1-s}\Phi\,F_t\,G_t$, where $G_t\equiv\partial_tF$.
\begin{equation*}
	\mathcal{J} = t^{1-s} \Phi_t\,(F\,G_t - F_t\,G) + \Phi\,\left(F\,\partial_t(t^{1-s}G_t) - \partial_t(t^{1-s}F_t)G\right)
\end{equation*}
With the $s$-harmonicity, that is $\operatorname{div}(t^{1-s}\nabla_{\Rnp}F)=0$ and therefore $\partial_t(t^{1-s}F_t)=-t^{1-s}\Delta_xF$, we have
\begin{equation*}
	\mathcal{J} = t^{1-s} \Phi_t\,(F\,G_t - F_t\,G) - t^{1-s}\Phi\,\left(F\,\Delta_xG - \Delta_xF\,G\right).
\end{equation*}
With a second cancellation, this time of the term $\nabla_x F\cdot \nabla_x G$,
\begin{equation*}
	\mathcal{J} = t^{1-s} \Phi_t\,(F\,G_t - F_t\,G) - t^{1-s}\Phi\,\nabla_x\cdot\left(F\,\nabla_xG - \nabla_xF\,G\right).
\end{equation*}
Plugging this back into $\mathcal{I}$, we obtain
\begin{equation*}
	\mathcal{I}\leq \left|\int_\Rnpp t^{1-s} \Phi_t\,(F\,G_t - F_t\,G)\right| 
					+ \left|\int_\Rnpp t^{1-s} \Phi\,\nabla_x\cdot\left(F\,\nabla_xG - \nabla_xF\,G\right) \right|
	\eqqcolon \mathcal{I}_1 + \mathcal{I}_2
\end{equation*}
Regarding the second term, with an integration by parts in $x$ where the boundary terms vanish due to the decay of the $s$-harmonic extensions, see Corollary \ref{cor:Pst decay of Pst and derviatives}, we can estimate the term as in (\ref{eq:fL comms suitable form for trace theorems}),
\begin{equation*}
	\mathcal{I}_2 =
	\left|\int_\Rnpp t^{1-s}\nabla_x\Phi\,\cdot\left(F\,\nabla_xG - \nabla_xF\,G\right) \right|
	\lesssim \int_{\Rnpp} t^{1-s}|\nabla_x \Phi|\,(|\nabla_xF|\,|G| + |\nabla_xG|\,|F|).
\end{equation*}
Regarding the first term, we need an additional derivative, since with $\partial_t \Phi$ via Proposition \ref{prop:tt - Lp est} we only get an estimate against $\fL{\nu}\phi$ for $\nu<s$.
Therefore, we use another integration by parts in $t$ to then obtain
\begin{equation*}
	\mathcal{I}_1 
	= \left|\int_\Rn\left[\frac{t^{s}}{s}\,t^{1-s}\Phi_t\,(F\,t^{1-s}G_t - t^{1-s}F_t\,G) \right]_0^\infty
				-\int_\Rnpp \frac{t^{s}}{s}\, \partial_t\left(t^{1-s}\Phi_t\,(F\,t^{1-s}G_t - t^{1-s}F_t\,G)\right)\right|,
\end{equation*}
where the boundary terms disappear due to (\ref{eq:Pst pde}), (\ref{eq:Pst partial t on boundary}) and Corollary \ref{cor:Pst decay of Pst and derviatives}.
Distributing the additional derivative, we observe the cancellation of the term $t^{2-s}\Phi_t\,F_t\,G_t$. Applying the $s$-harmonicity for the other terms, i.e. $\partial_t(t^{1-s}F)=-t^{1-s}\Delta_xF$, we have
\begin{equation}\label{eq:fL comms Vorbau for Rp comms}
	\mathcal{I}_1 = \frac{1}{s} \left| \int_\Rnpp t^{2-s}\Delta_x\Phi\,(F\,G_t-F_t\,G) 
									+ \int_\Rnpp t^{2-s}\Phi_t(F\,\Delta_xG - \Delta_x F\,G)\right|.
\end{equation}
With yet another cancellation, this time of the term $\nabla_xF\cdot\nabla_xG$,
\begin{equation*}
	\mathcal{I}_1 = \frac{1}{s} \left| \int_\Rnpp t^{2-s}\Delta_x\Phi\,(F\,G_t-F_t\,G) 
									+ \int_\Rnpp t^{2-s}\Phi_t\,\nabla_x\cdot(F\,\nabla_xG - \nabla_x F\,G)\right|.
\end{equation*}
After another integration by parts in $x$ for the second term, we can estimate $\mathcal{I}_1$ as in (\ref{eq:fL comms suitable form for trace theorems}) and thus finish this proof.
\end{proof}

\begin{remark}
At first glance, it might seem possible to obtain an additional estimate for $\|[\fL{s},\phi](f) \|_{L^p(\Rn)}$ from (\ref{eq:fL comms suitable form for trace theorems}) with Proposition \ref{prop:tt - BMO est} instead of Proposition \ref{prop:tt - Lp est}.
However, $[\phi]_{BMO}$ and $\|g\|_{L^{p'}(\Rn)}$ would be fixed as factors, leaving $\|\fL{s}f\|_{L^p(\Rn)}$ as the remaining factor. 
To obtain such an estimate with Proposition \ref{prop:tt - BMO est} is not possible for terms where $F$ occurs without derivative, such as for example $\int_\Rnpp t^{2-s}|\nabla_x\nabla_{\Rnp}\Phi| \,|\nabla_{\Rnp}G|\,|F|$.
\end{remark}

\subsection{The Chanilo commutator estimate for Riesz potentials of order $<1$}\label{sec:Rp comms}
Next, we take a look at the commutators of Riesz potentials and pointwise multiplication. 
For $\phi\in\CicRn$ and $s \in(0,n)$ we define the commutator $[I^s,\phi]$ by
\begin{equation*}
	[I^s,\phi](u) \coloneqq I^s(\phi\,u) - \phi\, I^s u, \quad \text{for all } u\in\CicRn.
\end{equation*}
In \cite{Cha82}, Chanillo proved the following $L^q$-estimate for these commutators, $q\in\left(\frac{n}{n-s},\infty\right)$.

\begin{theorem}[Chanillo]\label{th:Rp comm}
Let $\phi,u\in\CicRn$ and  $1<p<\frac{n}{s}$.
Define $q\in\left(\frac{n}{n-s},\infty\right)$ by $\frac{1}{q} = \frac{1}{p}-\frac{s}{n}$.
Then
\begin{equation*}
	\|[I^s,\phi](u)\|_{L^q(\Rn)}
	\lesssim
	[\phi]_{BMO}\,\|u\|_{L^p(\Rn)}.
\end{equation*}
\end{theorem}

The specific choice of $p$ and $q$ is due to the Hardy-Littlewood-Sobolev theorem, see Theorem \ref{th:Rp convergence and Lq boundedness}.
For this combination we have
\begin{equation*}
	\|I^s u \|_{L^q(\Rn)} \lesssim \|u\|_{L^p(\Rn)}.
\end{equation*}
Since the fractional Laplacians are inverse to the Riesz potentials, plugging in $\fL{s}f$ for $u$ we also obtain
\begin{equation}\label{eq:Rp comms Sobolev estimate for fL(s)}
	\|f \|_{L^q(\Rn)} \lesssim \|\fL{s}f\|_{L^p(\Rn)}.
\end{equation}
In general, we would prefer to state the theorem above in terms of fractional Laplacians instead of Riesz potentials, since only then we are able to make full use of the $s$-harmonic extension's behavior near the boundary.
By duality, Theorem \ref{th:Rp comm} is equivalent to
\begin{equation*}
	\mathcal{I}\coloneqq \left|\int_\Rn \left(I^s(\phi\,u) - \phi\, I^s u\right) \,v\,\right|
	\lesssim [\phi]_{BMO}\|u\|_{L^p(\Rn)}\|v\|_{L^{q'}(\Rn)}
	\qquad\text{for all } v\in L^{q'}(\Rn)
\end{equation*}
with $\frac{1}{q'}+\frac{1}{q}=1$.
By again plugging in $\fL{s}f$ for $u$ and with the integration-by-parts formula for Riesz potentials, see Lemma \ref{lem:Rp integration by parts}, we obtain
\begin{equation*}
	\mathcal{I} = \left|\int_\Rn \phi\,\fL{s}f\,I^sv - \phi\,f\,v\,\right|.
\end{equation*}
Note that $\fL{s}(\CicRn)$ is dense in $L^{q'}(\Rn)=\dot{F}^0_{q',2}(\Rn)$. 
This is the case since $\fL{s}$ is an isomorphism from the Triebel-Lizorkin space $\dot{F}^s_{q',2}(\Rn)$ to $L^{q'}(\Rn)=\dot{F}^0_{q',2}(\Rn)$, see Theorem \ref{th:tt lifting property of trlibeli spaces}, and $\CicRn$ is dense in $\dot{F}^s_{q',2}(\Rn)$, confer \cite[Theorem 5.1.5(ii), p.240]{Tri83}.
Therefore, we may replace the test function $v$ with $\fL{s}g$. Then,
\begin{equation*}
	\mathcal{I} = \left|\int_\Rn \phi\,\fL{s}f\,g - \phi\,f\,\fL{s}g \,\right| 
				= \left|\int_\Rn [\fL{s},\phi](f)\,g\,\right|.
\end{equation*}
Thus, by setting $f\coloneqq I^s u$, we obtain Theorem \ref{th:Rp comm} from the following theorem.
Note that we only prove this theorem for $s\in(0,1)$. 
Lenzmann and Schikorra suggest though, that by iterating the integration-by-parts procedure the result can be extended to $s\in(0,2)$ and even to the full range with a suitable higher order extension replacing the $s$-harmonic extensions, confer \cite[Section 5]{LeSchi20}.

\begin{theorem}
Let $s\in(0,1)$ and $p,q$ as above in Theorem \ref{th:Rp comm}.
Assume $\phi,f,g\in\C_c^\infty(\R)$.
Then for $\frac{1}{q'}+\frac{1}{q}=1$ we have
\begin{equation}\label{eq:Rp comms }
	\left|\int_\Rn [\fL{s},\phi](f)\,g\,\right|
	\lesssim
	[\phi]_{BMO} \,\|\fL{s}f\|_{L^p(\Rn)}\, \|\fL{s}g\|_{L^{q'}(\Rn)}.
\end{equation}
\end{theorem}

\begin{proof}
Let $\Phi(x,t)=P^s_t \phi(x)$, $F(x,t)=P^s_tf(x)$, $G(x,t)=P^s_t g(x)$ be the $s$-harmonic extensions of $\phi$, $f$, $g$.
Following the proof of Theorem \ref{th:fL comms} up to (\ref{eq:fL comms Vorbau for Rp comms}), we have 
\begin{equation*}
	\mathcal{I}\coloneqq \left|\int_\Rn [\fL{s},\phi](f)\,g\,\right| \leq \mathcal{I}_1+\mathcal{I}_2 
\end{equation*}
with
\begin{equation*}
	\mathcal{I}_2 =
	\left|\int_\Rnpp t^{1-s}\nabla_x\Phi\,\cdot\left(F\,\nabla_xG - \nabla_xF\,G\right) \right|
	\lesssim \int_{\Rnpp} t^{1-s}|\nabla_x \Phi|\,(|\nabla_xF|\,|G| + |\nabla_xG|\,|F|)
\end{equation*}
and 
\begin{equation*}
	\mathcal{I}_1 = \frac{1}{s} \left| \int_\Rnpp t^{2-s}\Delta_x\Phi\,(F\,G_t-F_t\,G) 
									+ \int_\Rnpp t^{2-s}\Phi_t(\Delta_x G\,F - \Delta_x F\,G)\right|.
\end{equation*}
Applying Proposition \ref{prop:tt - BMO est}(d) to $\mathcal{I}_2$ we have
\begin{equation}\label{eq:Rp comms I2 est}
	\mathcal{I}_2 \lesssim [\phi]_{BMO} \left(\|\fL{s}f\|_{L^p(\Rn)} \|g\|_{L^{p'}(\Rn)} +\|\fL{s}g\|_{L^{q'}(\Rn)} \|f\|_{L^{q}(\Rn)} \right)
\end{equation}
where $\frac{1}{p'}+\frac{1}{p}=1$. 

Regarding $\mathcal{I}_1$, with an integration by parts in $x$, since $\Delta_x=\nabla_x \cdot \nabla_x$,
\begin{mathex}[LCL]
	\mathcal{I}_1 &=& \frac{1}{s} \left| -\int_\Rnpp t^{2-s}\nabla_x\Phi\cdot\nabla_x(F\,G_t-F_t\,G) 
									+ \int_\Rnpp t^{2-s}\Phi_t(\Delta_x G\,F - \Delta_x F\,G)\right|\\
				&\lesssim& \int_\Rnpp t^{2-s}|\nabla_\Rnp\Phi|\,(|\nabla_x\nabla_\Rnp F|\, |G| + |\nabla_x\nabla_\Rnp G|\, |F|)\\
				&&			+ \int_\Rnpp t^{2-s}|\nabla_\Rnp\Phi|\,(|\nabla_x F|\, |\partial_tG| + |\nabla_x G|\, |\partial_t F|).							
\end{mathex}
Respecting the sequence of the functions in the respective summands, with Proposition \ref{prop:tt - BMO est}(c,d,e) we then obtain 
\begin{equation}\label{eq:Rp comms I1 est}
	\mathcal{I}_1 \lesssim [\phi]_{BMO} \left(\|\fL{s}f\|_{L^p(\Rn)} \|g\|_{L^{p'}(\Rn)} +\|\fL{s}g\|_{L^{q'}(\Rn)} \|f\|_{L^{q}(\Rn)} \right).
\end{equation}
Since
\begin{equation*}
	\frac{1}{q} = \frac{1}{p}-\frac{s}{n} 
	\quad\Leftrightarrow\quad
	1-\frac{1}{q'} = 1-\frac{1}{p'}-\frac{s}{n}
	\quad\Leftrightarrow\quad
	\frac{1}{p'} = \frac{1}{q'}-\frac{s}{n},
\end{equation*}
with (\ref{eq:Rp comms Sobolev estimate for fL(s)}) applied to (\ref{eq:Rp comms I2 est}) and (\ref{eq:Rp comms I1 est}) we obtain
\begin{equation*}
	\mathcal{I} \lesssim [\phi]_{BMO}\|\fL{s}f\|_{L^p(\Rn)} \|\fL{s}g\|_{L^{q'}(\Rn)}.
\end{equation*}
\end{proof}

\subsection{Fractional Leibniz rules: The three-term-commutator $H_s(f,g)$}\label{sec:fLr 3-t-comm}
The classical derivatives satisfy the Leibniz rule, that is
\begin{equation*}
	\nabla(f\,g)-\nabla f\,g - f\,\nabla g \equiv 0.
\end{equation*}
This is not true for all differential operators though. 
For example, for the Laplacian we have
\begin{equation*}
	\Delta(f\,g)-\Delta f\,g - f\,\Delta g = 2\nabla f \cdot \nabla g.
\end{equation*}
In this subsection, we prove two \textit{fractional} Leibniz-rules, which are estimates of the three-term-commutator
\begin{equation*}
	H_s(f,g)\coloneqq \fL{s}(f\,g)-\fL{s}f\,g - f\,\fL{s}g,
\end{equation*}
where $s>0$, $f,g\in\CicRn$.
Such fractional Leibniz-rules were originally introduced by Kenig-Ponce-Vega, confer \cite[Theorem A.9, p.611]{KPV93}. 

\begin{theorem}
Assume $s\in(0,1]$, $p\in(1,\infty)$ and $\phi,f\in\CicRn$. Then
\begin{equation}\label{eq:fLr 3-t-comm BMO}
	\|H_s(\phi,f)\|_{L^p(\Rn)} \lesssim [\phi]_{BMO} \|\fL{s}f\|_{L^p(\Rn)}.
\end{equation}
Moreover, suppose $\sigma\in(0,s)$, $p_1,p_2\in(1,\infty)$, $q,q_1,q_2\in[1,\infty]$ 
with $\frac{1}{p}=\frac{1}{p_1}+\frac{1}{p_2}$, $\frac{1}{q}=\frac{1}{q_1}+\frac{1}{q_2}$.
Then
\begin{equation}\label{eq:fLr 3-t-comm Lorentz}
	\|H_s(\phi,f)\|_{L^{(p,q)}(\Rn)} \lesssim \| \fL{\sigma}\phi\|_{L^{(p_1,q_1)}(\Rn)} \|\fL{s-\sigma}f\|_{L^{(p_2,q_2)}(\Rn)}.
\end{equation}
\end{theorem}

\begin{proof}[Estimate (\ref{eq:fLr 3-t-comm Lorentz})]
By duality, we only need to show that for any $g\in\CicRn$
\begin{equation}\label{eq:fLr 3-t-comm Lorentz duality estimate}
	\mathcal{I}\coloneqq \left|\int_{\Rn} H_s(\phi,f)\,g\,\right|
	\lesssim 	\| \fL{\sigma}\phi\|_{L^{(p_1,q_1)}(\Rn)} \|\fL{s-\sigma}f\|_{L^{(p_2,q_2)}(\Rn)}\|g\|_{L^{(p',q')}(\Rn)}
\end{equation}
where $\frac{1}{p'}+\frac{1}{p}=\frac{1}{q'}+\frac{1}{q}=1$, see \cite[Theorem 1.4.16, p.57]{Gra14} for duality regarding Lorentz spaces.
With the integration-by-parts formula for fractional Laplacians, Lemma \ref{lem:fL integration by parts}, we have
\begin{equation*}
	\mathcal{I} = \left| \int_\Rn \phi\,f\,\fL{s}g -\fL{s}\phi\,f\,g - \phi\,\fL{s}f\,g\,\right|.
\end{equation*}
Let $\Phi(x,t)=P^s_t\phi(x)$, $F(x,t)=P^s_t f(x)$, $G(x,t)=P^s_t g(x)$ be the $s$-harmonic extensions of $\phi$, $f$, $g$.
Then, via an integration by parts in $t$ we obtain
\begin{equation*}
	\mathcal{I}=\left|\int_{\Rnpp} 
		\partial_t\left(t^{1-s}\,\Phi\,F\,\partial_tG - t^{1-s}\,\partial_t\Phi\,F\,G - t^{1-s}\,\Phi\,\partial_tF\,G\right)
				\,\right|
	\eqqcolon \left|\int_{\Rnpp}\mathcal{J}\,\right|
\end{equation*}
as always due to the decay towards infinity and the behavior of the $s$-harmonic extensions near the boundary, Corollary \ref{cor:Pst decay of Pst and derviatives} and (\ref{eq:Pst partial t on boundary}), where we omit the constant $c$.
Distributing the additional derivative, we observe the cancellation of the terms $t^{1-s}\partial_t\Phi\,F\,\partial_tG$ and $t^{1-s}\Phi\,\partial_tF\,\partial_tG$;
\begin{equation*}
	\mathcal{J} = \Phi\,F\,\partial_t(t^{1-s}\partial_tG) 
				- \partial_t(t^{1-s}\partial_t\Phi)\,F\,G - \Phi\,\partial_t(t^{1-s}\partial_t F)\,G
				- 2t^{1-s}\,\partial_t\Phi\,\partial_t F\,G
\end{equation*}
Then, by the $s$-harmonicity of the extensions,
\begin{mathex}
	\mathcal{J} &=& -t^{1-s}\,\Phi\,F\,\Delta_xG + t^{1-s}\,\Delta_x\Phi\,F\,G + t^{1-s}\,\Phi\,\Delta_xF\,G 
						- 2t^{1-s}\,\partial_t\Phi\,\partial_t F\,G\\
				&=& t^{1-s}\,(\Delta_x(\Phi\,F)\,G - \Phi\,F\,\Delta_xG) 
					- 2t^{1-s}\,\nabla_x\Phi\cdot\nabla_xF\,G - 2t^{1-s}\,\partial_t\Phi\,\partial_t F\,G.
\end{mathex}
Plugging this back into $\mathcal{I}$, the first term disappears since $\Delta_x=\nabla_x\cdot\nabla_x$ and therefore, via integration by parts in $x$,
\begin{equation*}
	\int_\Rnpp t^{1-s}\,(\Delta_x(\Phi\,F)\,G - \Phi\,F\,\Delta_xG) 
	= \int_\Rnpp t^{1-s}\,(\nabla_x(\Phi\,F)\cdot\nabla_x G - \nabla_x(\Phi\,F)\cdot\nabla_xG) = 0.
\end{equation*}
Thus, we are left with
\begin{equation}\label{eq:fLr 3-t-comm Lorentz Vorbau for BMO}
	\mathcal{I}\lesssim \left|\int_\Rnpp t^{1-s}\nabla_x\Phi\cdot\nabla_xF\,G\, 
										+ t^{1-s}\partial_t\Phi\,\partial_t F\,G\,\right|
				\lesssim \int_\Rnpp t^{1-s}\,|\nabla_{\Rnp}\Phi|\,|\nabla_{\Rnp}F|\,|G|.
\end{equation}
Applying Proposition \ref{prop:tt - Lp est}, we estimate $\mathcal{I}$ as in (\ref{eq:fLr 3-t-comm Lorentz duality estimate}) and therefore have shown (\ref{eq:fLr 3-t-comm Lorentz}).
Note that we have $\sigma,\,s-\sigma < s=\min\lbrace 1,s\rbrace$ since $\sigma\in(0,s)$.
\end{proof}

\begin{proof}[Estimate (\ref{eq:fLr 3-t-comm BMO})]
Again, by duality, we only need to show that for any $g\in\CicRn$
\begin{equation}\label{eq:fLr 3-t-comm BMO duality estimate}
	\mathcal{I}\coloneqq \left|\int_{\Rn} H_s(\phi,f)\,g\,\right|
	\lesssim [\phi]_{BMO}\, \|\fL{s}f\|_{L^p(\Rn)}\,\|g\|_{L^{p'}(\Rn)}.
\end{equation}
Following the proof of (\ref{eq:fLr 3-t-comm Lorentz}) until (\ref{eq:fLr 3-t-comm Lorentz Vorbau for BMO}), we have
\begin{equation*}
	\mathcal{I}\lesssim \left|\int_\Rnpp t^{1-s}\nabla_x\Phi\cdot\nabla_xF\,G\,\right|
						+ \left|\int_\Rnpp t^{1-s}\partial_t\Phi\,\partial_t F\,G\,\right| 
				\eqqcolon \mathcal{I}_1+\mathcal{I}_2.
\end{equation*}
Regarding the first term, for $s<1$, with Proposition \ref{prop:tt - BMO est} we can already estimate 
\begin{equation*}
	\mathcal{I}_1	\lesssim \int_\Rnpp t^{1-s} |\nabla_x\Phi|\,|\nabla_x F|\,|G| 
					\lesssim [\phi]_{BMO} \|\fL{s}f\|_{L^p(\Rn)}\|g\|_{L^{p'}(\Rn)}.
\end{equation*}
For $s=1$, to obtain an estimate against $\fL{1}f$, we either need a second derivative on $f$ or any derivative on $g$.
Therefore, we again integrate by parts in $t$-direction, obtaining
\begin{equation*}
	\mathcal{I}_1	= \left|\int_\Rnpp \nabla_x\Phi\cdot\nabla_xF\,G\, \right| 
					=  \left|\int_\Rnpp t\partial_t\left(\nabla_x\Phi\cdot\nabla_xF\,G\right)\, \right|.
\end{equation*}
Distributing the $t$-derivative, where we integrate by parts in $x$ when that derivative hits $\nabla_x\Phi$,
\begin{mathex}
	\mathcal{I}_1	&=&	\left|\int_\Rnpp -t\,\partial_t\Phi\,\nabla_x\cdot(\nabla_xF\,G)
									    + t\,\nabla_x\Phi\cdot\partial_t\nabla_xF\,G 
									    + t\,\nabla_x\Phi\cdot\nabla_xF\,\partial_tG \,\right|\\
			&\lesssim&	\int_\Rnpp t\,|\nabla_\Rnp\Phi|\,(|\nabla_\Rnp\nabla_x F|\,|G| + |\nabla_\Rnp G|\,|\nabla_x F|).
\end{mathex}
Applying Proposition \ref{prop:tt - BMO est} while respecting the sequence of functions, we can estimate $\mathcal{I}_1$ as in (\ref{eq:fLr 3-t-comm BMO duality estimate}).
Regarding $\mathcal{I}_2$, we also need a second derivative on $f$ or any derivative on $g$ to estimate against $\fL{s}f$.
With an integration by parts in $t$-direction,
\begin{equation*}
	\mathcal{I}_2 = \left| \int_\Rn\left[\frac{t^s}{s}\,t^{1-s}\partial_t\Phi\,t^{1-s}\partial_t F\,G\, \right]_0^\infty
						- \int_\Rnpp \frac{t^s}{s}\,\partial_t\left(t^{1-s}\partial_t\Phi\,t^{1-s}\partial_t F\,G\right)\,
			\right|
\end{equation*}
Distributing the $t$-derivative, by the $s$-harmonicity of the extensions, i.e. $\partial_t(t^{1-s}\partial_t \Phi)=t^{1-s}\Delta_x\Phi$,
\begin{mathex}
	\mathcal{I}_2 &=&	\frac{1}{s}\left|\int_\Rnpp t^{2-s} \left(\Delta_x\Phi\,\partial_t F\,G 
														+ \partial_t\Phi\,\Delta_x F\,G  
														- \partial_t\Phi\,\partial_t F\,\partial_tG \right)\,\right|\\
		&\lesssim&	\int_\Rnpp t^{2-s}\,|\nabla_\Rnp \Phi|\,(\nabla_\Rnp\nabla_x F|\,|G| + |\nabla_\Rnp G|\,|\partial_t F|),
\end{mathex}
where we used that 
\begin{equation*}
	\int_\Rnpp t^{2-s}\,\Delta_x\Phi\,\partial_t F\,G = - \int_\Rnpp t^{2-s}\nabla_x\Phi\cdot\nabla_x(\partial_t F\,G).
\end{equation*}
With Proposition \ref{prop:tt - BMO est} we can estimate $\mathcal{I}_2$ as in (\ref{eq:fLr 3-t-comm BMO duality estimate}) and therefore obtain (\ref{eq:fLr 3-t-comm BMO}).
\end{proof}

\subsection{A Hardy-space estimate for $H_s(f,g)$}\label{sec:fLr 3-t-comm hardy}
We further investigate the three-term-commutator $H_s(f,g)$ from Subsection \ref{sec:fLr 3-t-comm}, defined as 
\begin{equation*}
	H_s(f,g)\coloneqq \fL{s}(f\,g)-\fL{s}f\,g - f\,\fL{s}g.
\end{equation*}
In the previous subsection, we estimated $L^p$ or Lorentz-norm for this commutator. Da Lio and Rivi\`ere showed the following Hardy-space estimate for the commutator with an additional fractional Laplacian in the case $s=1/2$, see 	\cite[Theorem 1.7]{DLR11}
\begin{theorem}[Da Lio-Rivière]
Assume $f,g\in\CicRn$. Then, with $\mathcal{H}^1(\Rn)$ denoting the Hardy-space,
\begin{equation*}
	\|(-\Delta)^\frac{1}{4}H_\frac{1}{2}(f,g)\|_{\mathcal{H}^1(\Rn)} 
	\lesssim 
	\|(-\Delta)^\frac{1}{4}f\|_{L^2(\Rn)} \|(-\Delta)^\frac{1}{4}g\|_{L^2(\Rn)}.
\end{equation*}
\end{theorem}

We proof a generalization of this result to $s\in(0,1]$.
This generalization is already known, see \cite[Theorem 1.4]{Schi15}. The original proof is substantially longer and more complicated than the proof from \cite{LeSchi20}, which we show here.

\begin{theorem}
Let $s\in(0,1]$.
Assume $f,\phi\in\CicRn$ and $p,p'\in(1,\infty)$, $q,q'\in[1,\infty]$ with $\frac{1}{p'}+\frac{1}{p}=\frac{1}{q'}+\frac{1}{q}=1$. Then
\begin{equation*}
	\| \fL{s}H_s(\phi,f)\|_{\mathcal{H}^1}	\lesssim \|\fL{s}\phi\|_{L^{(p,q)}(\Rn)}\,\|\fL{s}f\|_{L^{(p',q')}(\Rn)}.
\end{equation*}
\end{theorem}

\begin{proof}[s<1]
Assume $s<1$.
Since $BMO$ is the topological dual of $\mathcal{H}^1$, see the well known result in \cite{FeSte72}, we only need to show that for any $g\in \CicRn$
\begin{equation}\label{eq:fLr 3-t-comm hardy duality estimate}
	\mathcal{I}\coloneqq \left| \int_{\Rn}  \fL{s}H_s(\phi,f)\,g \right|
	\lesssim
	 \|\fL{s}\phi\|_{L^{(p,q)}(\Rn)}\, \|\fL{s}f\|_{L^{(p',q')}(\Rn)}\,[g]_{BMO}.
\end{equation}
Let $\Phi(s,t)=P^s_t\phi(x)$, $F(s,t)=P^s_tf(x)$, $\tilde{G}(s,t)=P^s_t(\fL{s}g)(x)$ be the $s$-harmonic extensions of $\phi$, $f$ and $\tilde{g}\coloneqq \fL{s}g$.
For the rest of the proof, we are going to show that
\begin{eqnarray}\label{eq:fLr 3-t-comm hardy suitable form}
	\mathcal{I} &\lesssim&
		\int_\Rnpp t^{1-s}\,|\tilde{G}|\,|\nabla_x\Phi|\,|\nabla_xF|\,\nonumber\\
	&&	+\int_\Rnpp t^{2-s}\,|\tilde{G}|\,(|\nabla_x\nabla_\Rnp\Phi|\,|\nabla_\Rnp F| + |\nabla_x\nabla_\Rnp F|\,|\nabla_\Rnp \Phi|)\\
	&&	+\int_\Rnpp t^{3-s}\,|\nabla_\Rnp\tilde{G}|\,(|\nabla_x\nabla_\Rnp\Phi|\,|\nabla_\Rnp F| + |\nabla_x\nabla_\Rnp F|\,|\nabla_\Rnp \Phi|)\nonumber
\end{eqnarray}
Then, since $\tilde{G}(s,t)=P^s_t(\fL{s}g)(x)$ and $s<1$, applying Proposition \ref{prop:tt - BMO est}(a,b) yields (\ref{eq:fLr 3-t-comm hardy duality estimate}).

Before extending the integral to $\Rnpp$, we redistribute some of the fractional Laplacians with Lemma \ref{lem:fL integration by parts},
\begin{equation*}
	\mathcal{I}= \left|\int_\Rnpp \phi\,f\,\fL{s}(\fL{s}g)-(\phi\,\fL{s}f + \fL{s}\phi\,f)\,\fL{s}g\,\right|.
\end{equation*}
With an integration by parts in $t$ we obtain
\begin{equation*}
	\mathcal{I}
	=	\left|\int_\Rnpp \partial_t \left(t^{1-s}\,\Phi\,F\,\partial_t\tilde{G} 
										- t^{1-s} \Phi\,\partial_tF\,\tilde{G}
										- t^{1-s} \partial_t\Phi\,F\,\tilde{G}\right)\, \right|
\end{equation*}
due to the extensions' decay towards infinity and their behavior near the boundary, that is $\fL{s}\phi=c\lim_{t\rightarrow0}\partial_t\Phi$, where we omit the constant $c$.
Distributing the additional derivative, we observe the cancellation of the term $t^{1-s}\partial_t(\Phi\,F)\,\partial_t\tilde{G}$,
\begin{mathex}
	\mathcal{I} &=&	\left|\int_\Rnpp \partial_t \left(t^{1-s}\,\Phi\,F\,\partial_t\tilde{G} 
													-t^{1-s} \partial_t(\Phi\,F)\,\tilde{G}\right)\, \right|\\
				&=& \left|\int_\Rnpp  (\Phi\,F)\,\partial_t\left(t^{1-s}\partial_t\tilde{G}\right)
										-\partial_t\left(t^{1-s}\partial_t(\Phi\,F)\right)\,\tilde{G}\, \right|.
\end{mathex}
With the $s$-harmonicity of $\tilde{G}$, that is $\partial_t\left(t^{1-s}\partial_t\tilde{G}\right)=-t^{1-s}\Delta_x\tilde{G}$,
\begin{equation*}
	\mathcal{I} =	\left|\int_\Rnpp 
							\left(t^{1-s}\Delta_x(\Phi\,F) + \partial_t\left(t^{1-s}\partial_t(\Phi\,F)\right)\right)\,
						\tilde{G}\, \right|
				= 	\left|\int_\Rnpp L_s(\Phi\,F)\,\tilde{G}\, \right|,
\end{equation*}
where we set 
\begin{equation*}
	L_s(H)\coloneqq \operatorname{div}(t^{1-s}\nabla_\Rnp H) = t^{1-s}\Delta_xH + \partial_t\left(t^{1-s}\partial_tH\right).
\end{equation*}
For this operator, we calculate the product rule
\begin{equation*}
	L_s(\Phi\,F)= L_s(\Phi)\,F + \Phi\,L_s(F) + 2t^{1-s}\,\nabla_x\Phi\cdot\nabla_xF + 2t^{1-s}\,\partial_t\Phi\,\partial_tF.
\end{equation*}
Obviously, with the $s$-harmonicity of $\Phi$ and $F$ we have $L(\Phi)=L(F)=0$ and therefore obtain
\begin{equation*}
	\mathcal{I}\lesssim \left|\int_\Rnpp t^{1-s}\,\nabla_x\Phi\cdot\nabla_xF\,\tilde{G}\,\right|
						+  \left|\int_\Rnpp t^{1-s}\,\partial_t\Phi\,\partial_tF\,\tilde{G}\,\right|
				\eqqcolon \mathcal{I}_1+\mathcal{I}_2.
\end{equation*}
The frist term can already be estimated as in (\ref{eq:fLr 3-t-comm hardy suitable form}).
Regarding the second term, if we wanted to estimate this term as in (\ref{eq:fLr 3-t-comm hardy duality estimate}) via Proposition \ref{prop:tt - BMO est}, we would need an additional derivative either on $\partial_t \Phi$ or $\partial_t F$.
Therefore, we again integrate by parts in $t$,
\begin{equation*}
	\mathcal{I}_2 = \left| \int_\Rn\left[\frac{t^s}{s}\,t^{1-s}\partial_t\Phi\,t^{1-s}\partial_tF\,\tilde{G}\right]_0^\infty
					- \frac{1}{s}\int_\Rnpp t^s\,\partial_t\left(t^{1-s}\partial_t\Phi\,t^{1-s}\partial_tF\,\tilde{G}\right)
				 \right|.
\end{equation*}
Distributing the additional derivative, due to the $s$-harmonicity of $\Phi$ and $F$ we obtain
\begin{mathex}[LCL]
 	\mathcal{I}_2&\lesssim& \frac{1}{s}\left(\left| \int_\Rnpp t^{2-s}\,\Delta_x\Phi\,\partial_tF\,\tilde{G} \,\right|
 									+ \left| \int_\Rnpp t^{2-s}\,\partial_t\Phi\,\Delta_xF\,\tilde{G} \,\right|
 									+ \left| \int_\Rnpp t^{2-s}\,\partial_t\Phi\,\partial_tF\,\partial_t\tilde{G}\,\right|\right)\\
 				&\eqqcolon& \mathcal{I}_{2a}+\mathcal{I}_{2b}+\mathcal{I}_{2c}.
\end{mathex}
The first two terms can already be estimated as in (\ref{eq:fLr 3-t-comm hardy suitable form}).
For the third term, we use yet another integration by parts in $t$.
\begin{mathex}
	s\mathcal{I}_{2c}
	&=&\frac{1}{2s} \left|\int_\Rn \left[t^{2s}\,t^{1-s}\partial_t\Phi\,t^{1-s}\partial_tF\,t^{1-s}\partial_t\tilde{G} \right]_0^\infty
	 - \int_\Rnpp t^{2s}\,\partial_t\left(t^{1-s}\partial_t\Phi\,t^{1-s}\partial_tF\,t^{1-s}\partial_t\tilde{G}\right)\,\right|\\
	 &\lesssim& \frac{1}{2s} \left|\int_\Rnpp t^{3-s} \Delta_x \Phi\,\partial_tF\,\partial_t\tilde{G}\right|
	 			+ \frac{1}{2s}\left|\int_\Rnpp t^{3-s} \partial_t \Phi\,\Delta_xF\,\partial_t\tilde{G}\right|\\
	&& 	+ \frac{1}{2s}\left|\int_\Rnpp t^{3-s} \nabla_x\left(\partial_t\Phi\,\partial_tF\right)\cdot\nabla_x\tilde{G}\right|
\end{mathex}
Thus, $\mathcal{I}_{2c}$ can also be estimated as in (\ref{eq:fLr 3-t-comm hardy duality estimate}).
\end{proof}

\begin{proof}[s=1]
The setup is exactly the same as in the proof for $s<1$. 
It suffices to show (\ref{eq:fLr 3-t-comm hardy duality estimate}).
For $s<1$, we obtained this estimate by applying Proposition \ref{prop:tt - BMO est} to (\ref{eq:fLr 3-t-comm hardy suitable form}). 
The only term in (\ref{eq:fLr 3-t-comm hardy suitable form}), which we are not able to immediately estimate analogously, is the first one, which corresponds to
\begin{equation*}
	\mathcal{I}_1 = \left|\int_\Rnpp t^{1-s}\,\nabla_x\Phi\cdot\nabla_xF\,\tilde{G}\,\right| 
					= \left|\int_\Rnpp \nabla_x\Phi\cdot\nabla_xF\,\tilde{G}\,\right|
\end{equation*}
For $s=1$, we need an additional derivative on $\nabla_x\Phi$ or $\nabla_xF$. 
We deal with $\mathcal{I}_1$ just as we dealt with $\mathcal{I}_2$ for $s<1$.
With an integration by parts in $t$ we obtain
\begin{equation*}
	\mathcal{I}_1 = \left|\int_\Rnpp t\,\partial_t\left(\nabla_x\Phi\cdot\nabla_xF\,\tilde{G}\right)\,\right|.
\end{equation*}
Distributing the additional derivative,
\begin{mathex}[LCL]
	\mathcal{I}_1 &\lesssim& \left|\int_\Rnpp t\,\nabla_x\partial_t\Phi\cdot\nabla_xF\,\tilde{G}\,\right|
							+\left|\int_\Rnpp t\,\nabla_x\Phi\cdot\nabla_x\partial_tF\,\tilde{G}\,\right|
							+\left|\int_\Rnpp t\,\nabla_x\Phi\cdot\nabla_xF\,\partial_t\tilde{G}\,\right|\\
	&\eqqcolon& \mathcal{I}_{1a}+\mathcal{I}_{1b}+\mathcal{I}_{1c}.
\end{mathex}
The first two terms can already be estimated as in (\ref{eq:fLr 3-t-comm hardy suitable form}).
The third term requires another integration by parts, where we use the harmonicity of $\tilde{G}$ when the additional derivative hits $\partial_t\tilde{G}$,
\begin{mathex}
	\mathcal{I}_{1c} 
	&=& \left|\int_\Rnpp t\,\partial_t\left(t\,\nabla_x\Phi\cdot\nabla_xF\,\partial_t\tilde{G}\right)\,\right|\\
	&\lesssim& \left|\int_\Rnpp t^2\,\nabla_x\partial_t\Phi\cdot\nabla_xF\,\partial_t\tilde{G}\,\right|
				+ \left|\int_\Rnpp t^2\,\nabla_x\Phi\cdot\nabla_x\partial_tF\,\partial_t\tilde{G}\,\right|\\
	&&		+ \left|\int_\Rnpp t^2\,\nabla_x\left(\nabla_x\Phi\cdot\nabla_xF\right)\cdot\partial_t\nabla_x\tilde{G}\,\right|.
\end{mathex}
Thus, we can estimate $\mathcal{I}_1$ just as $\mathcal{I}_2$ and therefore are left only with terms where we can obtain (\ref{eq:fLr 3-t-comm hardy duality estimate}) for $s=1$ with Proposition \ref{prop:tt - Lp est}.
\end{proof}

\newpage
\section{Riesz transforms, fractional Laplacians, Riesz potentials}\label{sec:operators}
In this section, we introduce Riesz transforms, the fractional Laplacian, and its inverse, the Riesz potential.
For each operator, we gather some essential and well known results from the literature, mainly from \cite{Ste70}.
Additionally, for the Riesz transforms and the fractional Laplacians we provide specific $L^\infty$- and decay-estimates, which are needed for Section \ref{sec:tt}.
\subsection{Riesz transforms}\label{subsec:Rt}
Riesz transforms are the prototypical singular integral operators, see \cite[Chapter III, §1]{Ste70} for a more detailed introduction.
One of their most common applications is the mediation between multiple differential operators of the same order.

\begin{definition}[Riesz transforms]\label{def:Rt}
For each $j\in\lbrace 1,\ldots, n\rbrace$ the corresponding Riesz transform is given by
\begin{mathex}
	\Rj(f)(x) 	=	\lim_{\epsilon\rightarrow 0} \const_n \int_{|y|\geq\epsilon}\frac{y_j}{|y|^{n+1}}f(x-y) \,dy
\end{mathex}
for all  $f\in L^p(\Rn)$, $1\leq p <\infty$ with $\const_n = \frac{\Gamma\left(\frac{n+1}{2}\right)}{\pi^{\frac{n+1}{2}}}$.
\end{definition}

These singular integral operators are well defined, see Theorem 4 in \cite[Chapter II]{Ste70}. 
That theorem further implies the $L^p$-boundedness of the Riesz transforms.

\begin{proposition}[$L^p$-boundedness of $\Rj$]\label{prop:Rt Lp bounded}
Assume $j\in\lbrace 1,\ldots, n\rbrace$ and $1<p<\infty$.
Then there exists $A_p>0$ such that for each $f\in L^p(\Rn)$
\begin{mathex}
	\|\Rj(f)\|_{L^p(\Rn)} \leq A_p \|f\|_{L^p(\Rn)}.
\end{mathex}
\end{proposition}

An important observation, which in fact justifies the definition of the Riesz transform, regards the interplay with the Fourier transform.
Since \cite{Ste70} uses a different variation of the Fourier transform than e.g. \cite{Gra14}, we will prefer the other sources for  references related to the Fourier transform.
The following result is Proposition 5.1.14. in \cite{Gra14}. 
An analogous statement can be found in Chapter III of \cite{Ste70}.

\begin{proposition}[Fourier symbol of Riesz transforms]\label{prop:Rt symbol}
The $j$th Riesz transform $\Rj$ is given on the Fourier transform side by multiplication with the function $-i \frac{\xi_j}{|\xi|}$.
That is, for any $f\in L^2(\Rn)$ it holds
\begin{mathex}
	\Fourier(\Rj(f))(\xi) = -i \frac{\xi_j}{|\xi|} \Fourier(f)(\xi).
\end{mathex}
\end{proposition}

An easy to see consequence of this proposition is the following lemma.

\begin{lemma}[An integration by parts formula for Riesz transforms]\label{lem:Rt integration by parts}
Assume $f,g\in L^2(\Rn)$ and $j\in\lbrace 1,\ldots,n\rbrace$. Then
\begin{equation*}
	\int_{\Rn} g \Rj f = -\int_{\Rn} \Rj g \, f.
\end{equation*}
\end{lemma}

\begin{proof}
Since the Fourier transform is an isometry on $L^2(\Rn)$, with $\Rj f(x)\in\R$ for almost every $x\in\Rn$ and Proposition \ref{prop:Rt symbol} we have
\begin{mathex}
	\int_{\Rn} g \Rj f 	&= \int_{\Rn} g \overline{\Rj f} 
						= \int_{\Rn} \Fourier(g)\overline{\Fourier(\Rj f)}
						= - \int_{\Rn} \left(-i\frac{\xi_i}{|\xi|}\right) \Fourier(g)(\xi) \overline{\Fourier(f)(\xi)}\,d\xi\\
						&= - \int_{\Rn} \Fourier(\Rj g) \overline{\Fourier(f)}
						= - \int_{\Rn} \Rj g \, f.
\end{mathex}
\end{proof}

Another consequence of Proposition \ref{prop:Rt symbol} is the following application, the already mentioned mediation between differential operators of the same order. 
For example, Riesz transforms can be used to mediate between the partial derivatives and the Laplacian of a function.

\begin{proposition}\label{prop:Rt mediation}
Assume $\phi\in\mathcal{S}(\Rn)$ and $j,k\in\lbrace 1,\ldots,n\rbrace$. Then
\begin{equation}\label{eq:Rt mediation}
	\partial_j \partial_k \phi(x) = \Rj \mathcal{R}_k (-\Delta) \phi(x).
\end{equation}
Together with the $L^p$-boundedness of Riesz transforms, for every $1<p<\infty$ this implies
\begin{mathex}
	\| \nabla^2 \phi \|_{L^p(\Rn)} \leq A_p^2 \|\Delta \phi\|_{L^p(\Rn)}.
\end{mathex}
\end{proposition}

\begin{proof}
For any $f\in\mathcal{S}(\Rn)$ we have
\begin{equation*}
	\Fourier(\partial_j f)(\xi) = 2\pi i \xi_j \Fourier(f)(\xi).
\end{equation*}
Therefore, applying the Fourier transform to the left side of (\ref{eq:Rt mediation}), we have
\begin{equation*}
	\Fourier(\partial_j \partial_k \phi)(\xi) = -4\pi^2 \xi_j \xi_k \Fourier(\phi)(\xi).
\end{equation*}
Thanks to Proposition \ref{prop:Rt symbol}, for the right side we have
\begin{equation*}
	\Fourier(-\Rj \mathcal{R}_k \Delta \phi)(\xi) 
	= \frac{\xi_j \xi_k}{|\xi|^2}\Fourier(\Delta \phi)(\xi)
	= \frac{\xi_j \xi_k}{|\xi|^2}  \left(\sum_{l=1}^n -4\pi^2 \xi_l^2 \right) \Fourier(\phi)(\xi)
	= -4\pi^2 \frac{\xi_j \xi_k}{|\xi|^2}|\xi|^2 \Fourier(\phi)(\xi).
\end{equation*} 
Now, applying the inverse Fourier transform, we have shown (\ref{eq:Rt mediation}). 
\end{proof}

This same type of argument is used in Section \ref{sec:Rt comms} to show similar relations between other differential operators.
For the last two results of this subsection, we do not use the Fourier characterization of the Riesz transforms, but rather fall back on the original definition as a singular integral operator.

\begin{lemma}[The decay of $\Rj f$]\label{lem:Rt decay}
Let $j\in\lbrace 1,\ldots,n\rbrace$.
Assume $f\in C^1(\Rn) \cap L^p(\Rn)$, $1\leq p <\infty$, with $k\in\R$ and $\const>0$ such that for all $x\in\Rn$
\begin{equation*}
	|\nabla f(x)| \leq \frac{C}{|x|^k}.
\end{equation*}
Then there exists a constant $\const_f>0$ such that for all $x\in\Rn$
\begin{equation*}
	|\Rj f(x)| \leq \const_f \left(\frac{1}{|x|^\frac{n}{p}} + \frac{1}{|x|^{k-1}} \right).
\end{equation*}
\end{lemma}

\begin{proof}
Let $x\in\Rn$.
By the definition of the Riesz transforms we have
\begin{mathex}
	|\Rj f(x)| 
	&\leq&		\left|\int_{\Rn\setminus B_\frac{|x|}{2}(x)} \frac{x_i-y_i}{|x-y|^{n+1}}f(y)\,dy \right|
				+ \left| \lim_{\epsilon\rightarrow0}
								\int_{B_\frac{|x|}{2}(x)\setminus B_\epsilon(x)} \frac{x_i-y_i}{|x-y|^{n+1}}f(y)\,dy\right|\\
	&\eqqcolon&	\const_1 + \const_2.
\end{mathex}
Regarding the second term, first note that, since $y\mapsto (x_i-y_i)|x-y|^{-(n+1)}$ is odd, for any $0<r<R\leq\infty$ we have 
\begin{equation}\label{eq:Rt kernel integral vanishes}
	\int_{B_R(x)\setminus B_r(x)} \frac{x_i-y_i}{|x-y|^{n+1}}dy = 0.
\end{equation} 
Applying this observation, we obtain
\begin{mathex}
	\const_2
	&=&	\left| \lim_{\epsilon\rightarrow0}
							\int_{B_\frac{|x|}{2}(x)\setminus B_\epsilon(x)} \frac{x_i-y_i}{|x-y|^{n+1}}(f(y)-f(x))\,dy\right|\\
	&\leq& 	\lim_{\epsilon\rightarrow0} \int_{B_\frac{|x|}{2}(x)\setminus B_\epsilon(x)} \frac{|f(y)-f(x)|}{|x-y|^{n}}\,dy\\
	&\leq&	\sup_{\tilde{x}\in B_\frac{|x|}{2}(x)} |\nabla f(\tilde{x})| 
					\lim_{\epsilon\rightarrow 0}\int_{B_\frac{|x|}{2}(x)\setminus B_\epsilon(x)} \frac{|x-y|}{|x-y|^{n}}\,dy \\
	&\leq&	\frac{2^k C}{|x|^k} \omega_n \int_0^\frac{|x|}{2}	\frac{r^{n-1}}{r^{n-1}}dr\\
	&=&		2^{k-1}C\omega_n\frac{1}{|x|^{k-1}}
\end{mathex}
where $\omega_n$ denotes the $n-1$-dimensional volume of the sphere $S^{n-1}$.

Regarding the first term, with Hölder's inequality and $\frac{1}{p'}+\frac{1}{p}=1$ we obtain
\begin{mathex}
	\const_1
	&\leq&			\int_{\Rn\setminus B_\frac{|x|}{2}(x)} \frac{1}{|x-y|^{n}}|f(y)|\,dy\\
	&\leq&			\|y\mapsto |y|^{-n}\|_{L^{p'}(\Rn\setminus B_{|x|/2}(0))} \|f\|_{L^p(\Rn)}\\
	&=&			C |x|^{-\frac{n}{p}} \|f\|_{L^p(\Rn)}
\end{mathex}
for a constant $C>0$ depending on $n$ and $p$.
\end{proof}

\begin{lemma}[An $L^\infty$-estimate for $\Rj f$]\label{lem:Rt Linfty estimate}
Let $j\in\lbrace 1,\ldots,n\rbrace$.
Assume $f\in C^1(\Rn) \cap L^p(\Rn)$ with $1\leq p <\infty$ and $\nabla f \in L^\infty(\Rn)$.
Then
\begin{equation*}
	\|\Rj f\|_{L^\infty(\Rn)} \lesssim \|f\|_{L^p(\Rn)} + \|\nabla f\|_{L^\infty(\Rn)}.
\end{equation*}
\end{lemma}

\begin{proof}
Let $x\in\Rn$.
Similar to the previous proof, we have 
\begin{mathex}
	|\Rj f(x)| 
	&\leq&		\left|\int_{\Rn\setminus B_1(x)} \frac{x_i-y_i}{|x-y|^{n+1}}f(y)\,dy \right|
				+ \left| \lim_{\epsilon\rightarrow0}
								\int_{B_1(x)\setminus B_\epsilon(x)} \frac{x_i-y_i}{|x-y|^{n+1}}f(y)\,dy\right|\\
	&\eqqcolon&	\const_1 + \const_2.
\end{mathex}
Regarding the first term, with Hölder's inequality,
\begin{equation*}
	\const_1
	\leq			\int_{\Rn\setminus B_1(x)} \frac{1}{|x-y|^{n}}|f(y)|\,dy
	\leq			\|y\rightarrow |y|^{-n}\|_{L^{p'}(\Rn\setminus B_1(0))} \|f\|_{L^p(\Rn)}.
\end{equation*}
where $\frac{1}{p'}+\frac{1}{p}=1$, in particular $p'\in(1,\infty]$ and $\|y\rightarrow |y|^{-n}\|_{L^{p'}(\Rn\setminus B_1(0))}<\infty$.
For the second term we obtain
\begin{equation*}
	\const_2
	=	\left| \lim_{\epsilon\rightarrow0}
							\int_{B_1(x)\setminus B_\epsilon(x)} \frac{x_i-y_i}{|x-y|^{n+1}}(f(y)-f(x))\,dy\,\right|
	\leq	\|\nabla f\|_{L^\infty(\Rn)} \int_{B_1(x)} \frac{|x-y|}{|x-y|^{n}}\,dy,
\end{equation*}
where we used (\ref{eq:Rt kernel integral vanishes}) as well as the mean value theorem.
\end{proof}

\subsection{Fractional Laplacians}
We already saw in the proof of Proposition \ref{prop:Rt mediation} above how the Laplacian interacts with the Fourier transform,
\begin{equation*}
	\Fourier(-\Delta f)(\xi) = (2\pi |\xi|)^2 \Fourier(f)(\xi) \qquad \text{for any $f\in\mathcal{S}(\Rn)$.}
\end{equation*}
This motivates the following definition of the fractional Laplacian operator.

\begin{definition}[Fractional Laplacian]\label{def:fL fL}
Let $s>0$ and $f\in\mathcal{S}(\Rn)$.
Then the fractional Laplacian operator $\fL{s}$ is defined by
\begin{mathex}
	\Fourier(\fL{s}f)(\xi) = (2\pi |\xi|)^s \Fourier(f)(\xi).
\end{mathex}
\end{definition}

Note that for each $s>0$ and $f\in\mathcal{S}(\Rn)$ the function $\xi\mapsto(2\pi |\xi|)^s \Fourier(f)(\xi)$ is still rapidly decreasing, although it might not be differentiable in $\xi=0$ anymore. 
In particular, we have an $L^1$ function and therefore can apply the inverse Fourier transform.
Thus, the fractional Laplacian operator is well defined. 
Throughout this thesis, fractional Laplacians will often be applied to functions that are not strictly of the Schwartz space, see for example Lemma \ref{lem:fL semi group} below.
In most cases, it is easy to see that the function the fractional Laplacian is applied to is in $L^2(\Rn)$ or $L^1(\Rn)$ with a sufficiently decreasing Fourier transform, resulting in the fractional Laplacian being well defined with the same reasoning as for Schwartz functions. 
We will only investigate the well-definedness further, if this is not the case.

For later applications, we will mostly consider the case $s\in (0,2)$.
Under this restriction, the definition can be extended to other function spaces such as $L^p(\Rn)$ with $1\leq p \leq 2$, see \cite{Kwa17}.
Note that there is a variety of different, equivalent ways to introduce and define the fractional Laplacian, again see \cite{Kwa17}. 
For example, the fractional Laplacian may be defined similar to the Riesz transform as a singular integral operator. 
We pick up this equivalent definition towards the end of this chapter.

Since the fractional Laplacians act as Fourier multipliers by definition, there are some easy to see consequences, similar to those we got for the Riesz transforms after Proposition \ref{prop:Rt symbol}.
First, we get an integration by parts formula. 

\begin{lemma}[An integration by parts formula for fractional Laplacians]\label{lem:fL integration by parts}
Let $f,g\in\mathcal{S}(\Rn)$ and $s>0$.
Then
\begin{mathex}
	\int_{\Rn} g \fL{s} f = \int_{\Rn} \fL{s} g \, f.
\end{mathex}
\end{lemma}

\begin{proof}
Analogous to Lemma \ref{lem:Rt integration by parts}. We have $\fL{s} f(x)\in\R$ for all $x\in\Rn$ because of Proposition \ref{prop:fL singular integral definition}, the alternative definition of the fractional Laplacian.
\end{proof}

Second, we observe an interaction between multiple fractional Laplacians which we can describe as a semigroup property. 

\begin{lemma}[Semigroup property of fractional Laplacians]\label{lem:fL semi group}
Let $f\in\mathcal{S}(\Rn)$ and $s_1,s_2>0$.
Then
\begin{mathex}
	\fL{s_1}\fL{s_2}f = \fL{s_1+s_2}f.
\end{mathex}
\end{lemma}

\begin{proof}
Applying the Fourier transform to the equation, via the definition of the fractional Laplacian we get
\begin{mathex}
	\Fourier(\fL{s_1}\fL{s_2}f)(\xi) = (2\pi |\xi|)^{s_1} (2\pi |\xi|)^{s_2} \Fourier(f)(\xi)
	= (2\pi |\xi|)^{s_1+s_2} \Fourier(f)(\xi) = \Fourier(\fL{s_1+s_2}f).
\end{mathex}
Now applying the inverse Fourier transform, we have proven the equation.
\end{proof}

Comparing the Fourier symbols of $\fL{1}$ and $\partial_j$, $1\leq j \leq n$, one would expect a relation similar to Proposition \ref{prop:Rt mediation}, which dealt with the usual Laplacian operator.
And indeed, we get a similar result for the fractional Laplacian with $s=1$, although the proof is not quite as simple.

\begin{proposition}\label{prop:fLl vs grad}
Let $f\in\mathcal{S}(\Rn)$ and $1<p<\infty$.
Then
\begin{equation}
	\|\fL{1} f\|_{L^p(\Rn)} \approx \|\nabla f\|_{L^p(\Rn)}.
\end{equation}
\end{proposition}

\begin{proof}
Both $\|\fL{1} f\|_{L^p(\Rn)}$ and $\|\nabla f\|_{L^(\Rn)}$ are equivalent norms on the Triebel-Lizorkin space $\dot{F}^1_{p,2}$, see \ref{th:tt lifting property of trlibeli spaces} or \cite[Theorem 5.2.3.1]{Tri83} for the original source.
\end{proof}

We conclude this subsection about the fractional Laplacian estimating the $L^\infty$-norm and investigating the decay of $\fL{s}f$.
In order to show these estimates, we make use of one of the alternative definitions of the fractional Laplacian.

\begin{proposition}[Alternative definition of fractional Laplacians]\label{prop:fL singular integral definition}
Let $s\in(0,2)$. 
Then for $f\in\mathcal{S}$
\begin{mathex}
	\fL{s}f(x) 	&=& \const_{n,s} \lim_{\epsilon\rightarrow 0} \int_{\Rn\setminus B_\epsilon(x)} \frac{f(x)-f(y)}{|x-y|^{n+s}}dy\\
				&=& -\frac{1}{2}\const_{n,s} \int_{\Rn}\frac{f(x+y)+f(x-y)-2f(x)}{|y|^{n+s}}dy
\end{mathex}
with $\const_{n,s} = \frac{2^s\Gamma\left(\frac{n+s}{2}\right)}{\pi^\frac{n}{2}\left|\Gamma\left(-\frac{s}{2}\right)\right|} = \left(\int_\Rn \frac{1-\cos(\xi_1)}{|\xi|^{n+2s}}d\xi\right)^{-1}$.
\end{proposition}

\begin{proof}
For the first equation see \cite{Kwa17} or Proposition 3.3 in \cite{DiN12}. 
The second equation follows from rewriting the principal value integral and then using a second order taylor expansion to remove the singularity, see Lemma 3.2 in \cite{DiN12}. 
\end{proof}

Again, this proposition or alternative definition respectively can be extended to $L^p(\Rn)$, $p\in[1,\infty)$, see \cite{Kwa17}.

\begin{lemma}[An $L^\infty$-estimate for $\fL{s}f$]\label{lem:fL Linfty estimate}
Let $s\in(0,2)$ and $f\in\mathcal{S}(\Rn)$.
Then 
\begin{equation}\label{eq:fL Linfty estimate 1}
	\|\fL{s}f\|_{L^\infty(\Rn)} \lesssim \|f\|_{L^\infty(\Rn)} + \|\nabla^2 f \|_{L^\infty(\Rn)}.
\end{equation}
For $s\in(0,1)$, we also have
\begin{equation}\label{eq:fL Linfty estimate 2}
	\|\fL{s}f\|_{L^\infty(\Rn)} \lesssim \|f\|_{L^\infty(\Rn)} + \|\nabla f \|_{L^\infty(\Rn)}.
\end{equation}
\end{lemma}

\begin{proof}
We show (\ref{eq:fL Linfty estimate 2}) first. 
Let $s\in(0,1)$ and $x\in\Rn$. Then
\begin{mathex}
	|\fL{s}f(x)|
	&=&			\left| \lim_{r\rightarrow 0} \const_{n,s} \int_{\Rn\setminus B_r(0)} \frac{f(x+y)-f(x)}{|y|^{n+s}}dy\right|\\
	&\lesssim&	\int_{\Rn\setminus B_1(0)} \frac{|f(x+y)-f(x)|}{|y|^{n+s}} dy 
				+ \lim_{r\rightarrow 0}  \int_{B_1(0)\setminus B_r(0)} \frac{|f(x+y)-f(x)|}{|y|^{n+s}} dy
\end{mathex}
With the mean value theorem we have $|f(x+y)-f(x)|\leq \|\nabla f\|_{L^\infty(\Rn)}|y|$ for the second term.
Therefore, we obtain
\begin{mathex}				
	|\fL{s}f(x)|		
	&\leq&		2\|f\|_{L^\infty(\Rn)}  \int_{\Rn\setminus B_1(0)} \frac{1}{|y|^{n+s}}dy 
				+ \|\nabla f\|_{L^\infty(\Rn)} \lim_{r\rightarrow 0}  \int_{B_1(0)\setminus B_r(0)} \frac{|y|}{|y|^{n+s}} dy\\
	&\leq&		\left(2 \int_{\Rn\setminus B_1(0)} |y|^{-(n+s)}dy + \int_{B_1(0)} |y|^{-(n+s-1)} dy \right) (\|f\|_{L^\infty(\Rn)}+\|\nabla f\|_{L^\infty(\Rn)}).
\end{mathex}
The first integral exists because $s>0$, the second one because $s<1$.
Thus, we have shown (\ref{eq:fL Linfty estimate 2}).

Now, let $s\in(0,2)$ and $x\in\Rn$. Then with a second order Taylor expansion yielding 
\begin{equation}\label{eq:fL taylor expansion}
	|f(x+y)+f(x-y)-2f(x)| = |y^T \nabla^2 f(x+\theta(y-x)) y| \leq \|\nabla^2 f\|_{L^\infty(\Rn)}|y|^2,
\end{equation}
where $\theta\in[-1,1]$, we obtain
\begin{mathex}
|\fL{s}f(x)|
	&=&			\left| \frac{1}{2} \const_{n,s} \int_{\Rn} \frac{f(x+y)+f(x-y)-2f(x)}{|y|^{n+s}}dy\right|\\
	&\lesssim&	4 \|f\|_{L^\infty(\Rn)} \int_{\Rn\setminus B_1(0)} |y|^{-(n+s)}dy 
				+ \|\nabla^2 f\|_{L^\infty(\Rn)}\int_{B_1(0)} |y|^{n+s-2} dy \\
	&\lesssim&	\left(4 \int_{\Rn\setminus B_1(0)} |y|^{-(n+s)}dy  + \int_{B_1(0)} |y|^{n+s-2} dy\right) 
				(\|f\|_{L^\infty(\Rn)} + \|\nabla^2 f\|_{L^\infty(\Rn)}).
\end{mathex}
Since the integrals exist due to $s\in(0,2)$, we have shown (\ref{eq:fL Linfty estimate 1}).
\end{proof}

\begin{remark}
By not dividing $\Rn$ into $\Rn\setminus B_1(0)$ and $B_1(0)$ but rather varying the radius for the cut out ball, one can easily obtain weighted inequalities such as
\begin{equation*}
	\|\fL{s}f\|_{L^\infty(\Rn)} \leq  \const(\epsilon) \|f\|_{L^\infty(\Rn)} + \epsilon\|\nabla^2 f \|_{L^\infty(\Rn)},
\end{equation*}
where $\epsilon>0$ can be chosen arbitrarily and $\const(\epsilon)>0$ is a constant depending on the choice.
An analogous observation can be made for example for Lemma \ref{lem:Rt Linfty estimate}.
\end{remark}

For our last basic result regarding the fractional Laplacian, we investigate its decay. 
The approach for the proof of the following lemma is similar to the one for Lemma \ref{lem:Rt decay}.

\begin{lemma}[The decay of $\fL{s}f$]\label{lem:fL decay}
Let $s\in(0,2)$.
Assume $f\in C^2(\Rn) \cap W^{2,\infty}(\Rn) \cap L^1(\Rn)$ with $k,l\in\R$ and $\const>0$ such that for all $x\in\Rn$
\begin{equation*}
	|f(x)|\leq \frac{C}{|x|^k}, \quad |\nabla^2 f(x)| \leq \frac{C}{|x|^l}.
\end{equation*}
Then there exists a constant $\const_f>0$ such that for all $x\in \Rn$
\begin{equation*}
	|\fL{s}f(x)| \leq \const_f \left(\frac{1}{|x|^{n+s}} + \frac{1}{|x|^{k+s}} + \frac{1}{|x|^{l+s-2}}\right).	
\end{equation*}
\end{lemma}

\begin{proof}
Let $x\in\Rn$.
With Proposition \ref{prop:fL singular integral definition} we have 
\begin{mathex}
	|\fL{s}f(x)| 
	&\lesssim &	\left|\int_{\Rn\setminus B_{\frac{|x|}{2}}(0)} \frac{f(x+y)+f(x-y)-2f(x)}{|y|^{n+s}}dy\right|\\
			&&	+ \left|\int_{B_{\frac{|x|}{2}}(0)} \frac{f(x+y)+f(x-y)-2f(x)}{|y|^{n+s}}dy\right|\\
	&\eqqcolon& \const_1 + \const_2.
\end{mathex}
Let $\omega_n$ denote the $n-1$-dimensional volume of the sphere $S^{n-1}$. For the first term we have
\begin{mathex}	
	\const_1
	&\leq	&	2 \frac{2^{n+s}}{|x|^{n+s}}\|f\|_{L^1(\Rn)} 
				+ 2 |f(x)| \int_{\Rn\setminus B_{\frac{|x|}{2}}(0)} \frac{1}{|y|^{n+s}}dy\\
	&\leq	&	\frac{2^{n+s+1}}{|x|^{n+s}}\|f\|_{L^1(\Rn)} 
				+ \frac{2C}{|x|^k} \omega_n \int_\frac{|x|}{2}^\infty \frac{r^{n-1}}{r^{n+s}}dr\\
	&=		&	\frac{2^{n+s+1}}{|x|^{n+s}}\|f\|_{L^1(\Rn)} + \frac{2C\omega_n}{|x|^k} \frac{1}{s}\frac{2^s}{|x|^s}.
\end{mathex}
Regarding the second term, with the second order Taylor expansion (\ref{eq:fL taylor expansion}) we obtain
\begin{mathex}		
	\const_2
	&\leq	&	\sup_{\tilde{x}\in B_\frac{|x|}{2}(x)} |\nabla^2 f(\tilde{x})| \int_{B_\frac{|x|}{2}(0)} \frac{1}{|y|^{n+s-2}}dy
				\leq \frac{2^l C}{|x|^l}\omega_n \int_0^\frac{|x|}{2} \frac{r^{n-1}}{r^{n+s-2}}dr
				= \frac{2^l C\omega_n}{|x|^l}\frac{1}{2-s}\frac{|x|^{2-s}}{|2|^{2-s}}
\end{mathex}
\end{proof}

Having now covered the basics for fractional Laplacians, we turn to their inverse, the Riesz potentials.

\subsection{Riesz potentials}
Concluding this section, we introduce the Riesz potentials, which are inverse to the corresponding fractional Laplacians.
These operators are defined via singular integrals, see \cite[Chapter V]{Ste70}. 

\begin{definition}[Riesz potential]\label{def:Rp}
Let $0<s<n$ and $f\in\mathcal{S}(\Rn)$.
Then the fractional Laplacian operator $\fL{s}$ is defined by
\begin{equation}\label{eq:Rp definition}
	I^{s}f(x) = \const_{n,s} \int_{\Rn} \frac{f(x-y)}{|y|^{n-s}}
\end{equation}
with $\const_{n,s} = \frac{\Gamma\left(\frac{n-s}{2}\right)}{2^s \pi^\frac{n}{2}\Gamma\left(\frac{s}{2}\right)}$.
\end{definition}

From \cite[Chapter V]{Ste70} we also obtain a result regarding the pointwise convergence of these singular integrals as well as potential $L^p$-boundedness. 

\begin{theorem}[Hardy-Littlewood-Sobolev theorem of fractional integration]\label{th:Rp convergence and Lq boundedness}
Let $0<s<n$, $1\leq p < q < \infty$, $\frac{1}{q}= \frac{1}{p}- \frac{s}{n}$, in particular $p < \frac{n}{s}$.
Then for any $f \in L^p(\Rn)$ the integral (\ref{eq:Rp definition}), defining $I^{s}f$, converges absolutely for almost every $x\in\Rn$.

If, in addition, $1<p$, then there exists $\const_{p,q}>0$ such that for all $f\in L^p(\Rn)$
\begin{equation*}
	\|I^s f \|_{L^q(\Rn)} \leq \const_{p,q} \|f\|_{L^p(\Rn)}
\end{equation*}
\end{theorem}

\begin{proof}
See \cite{Ste70}, Theorem 1 in Chapter V.
\end{proof}

Thanks to this result, we can extend the definition of the Riesz Potentials from $\mathcal{S}(\Rn)$ to $L^p(\Rn)$, $1\leq p< \frac{n}{s}$.
Regarding the following results, we will not consider the original definition though, but rather focus on the interplay with the Fourier transform. 
Like Riesz transforms and fractional Laplacians, Riesz potentials are multiplier operators.

\begin{proposition}[Fourier symbol of Riesz potentials]\label{prop:Rp symbol}
Let $0<s<n$.
The corresponding Riesz potential is given on the Fourier transform side by multiplication with the function $(2\pi\xi)^{-s}$.
That is, for any $f\in\mathcal{S}(\Rn)$ we have
\begin{mathex}
	\Fourier(I^s f)(\xi) = (2\pi |\xi|)^{-s} \Fourier(f)(\xi).
\end{mathex}
\end{proposition}

\begin{proof}
See \cite{Ste70}, Lemma 1 in Chapter V.
\end{proof}


As in the subsections above, we list some easy to see consequences of this result.
First, we can immediately confirm that indeed fractional Laplacians and Riesz potentials are inverse to each other.
Second and third, we again get an integration by parts formula and a semigroup property.

\begin{lemma}
Let $s\in(0,n)$ and $f\in\mathcal{S}(\Rn)$.
Then
\begin{equation*}
	I^s\fL{s} f = \fL{s} I^{s} f = f.
\end{equation*}
\end{lemma}

\begin{proof}
We have $\fL{s}f\in L^1(\Rn)$ due to Lemma \ref{lem:fL Linfty estimate} and Lemma \ref{lem:fL decay}. 
Moreover, we have $\fL{s}f\in L^2(\Rn)$  because obviously $\Fourier(\fL{s}f)\in L^2(\Rn)$. 
Combining these, for any $1\leq p \leq 2$ we have $\fL{s}f\in L^p(\Rn)$.
Therefore, the application of the Riesz potential $I^s$ to $\fL{s}f$ is well defined with Theorem \ref{th:Rp convergence and Lq boundedness}, choosing $p_0\in\left(1,\min\left\lbrace\frac{n}{s},2\right\rbrace\right)$ and $\frac{1}{q_0} = \frac{1}{p_0} - \frac{s}{n}$.
Furthermore, we can find a sequence $(g_k)_{k\in\N}\subset\mathcal{S}(\Rn)$ such that $\|g_k-\fL{s}f\|_{L^p(\Rn)}\overset{k\rightarrow\infty}{\longrightarrow} 0$ for all $1\leq p \leq 2$.
In particular, this implies that 
\begin{equation}\label{eq:Rp inverse to fL (Lp' convergence)}
	\|\Fourier(g_k)-\Fourier(\fL{s}f)\|_{L^{p'}(\Rn)}\overset{k\rightarrow\infty}{\longrightarrow} 0
	\quad \text{for all } 1\leq p \leq 2 \text{ and } \frac{1}{p'}+\frac{1}{p}=1
\end{equation} 
thanks to the Hausdorff-Young inequality.
Additionally, because of Theorem \ref{th:Rp convergence and Lq boundedness}, $I_s g_k$ converges against $I^s(\fL{s} f)$ in $L^{q_0}(\Rn)$ and therefore also in $\mathcal{S}'(\Rn)$. 
For the Fourier transform of $I^s (\fL{s}f)$ in the sense of tempered distributions we then have
\begin{equation*}
	\Fourier(I^s g_k) \overset{\mathcal{S}'}{\longrightarrow} \Fourier(I^s (\fL{s}f)) \quad \text{for } k\rightarrow\infty.
\end{equation*}
Now let $h\in\mathcal{S}(\Rn)$. We have
\begin{equation*}
	\const	\coloneqq \int_\Rn \Fourier(I^s\fL{s}f) h 
			= \lim_{k\rightarrow\infty} \int_\Rn \Fourier(I^s g_k) h 
			= \lim_{k\rightarrow\infty} \int_\Rn (2\pi|\xi|)^{-s}\Fourier(g_k)(\xi) h(\xi) \,d\xi
\end{equation*}
thanks to Proposition \ref{prop:Rp symbol}.
We now observe that $(2\pi|\cdot|)^{-s}h \in L^{\tilde{p}}$ for some $1<\tilde{p}<\min\left\lbrace\frac{n}{s},2\right\rbrace$ since $s<n$.
With (\ref{eq:Rp inverse to fL (Lp' convergence)}), via the duality of $L^p$-spaces and the definition of the fractional Laplacian,
\begin{equation*}
	\const	= \int_\Rn \Fourier(\fL{s}f)(\xi) (2\pi|\xi|)^{-s}h(\xi) \,d\xi 
			= \int_\Rn (2\pi|\xi|)^{s} \Fourier(f)(\xi) (2\pi|\xi|)^{-s}h(\xi) \,d\xi
			= \int_\Rn \Fourier(f)h,
\end{equation*} 
yielding $\Fourier(I^s\fL{s}f)=\Fourier(f)\in\mathcal{S}$ and therefore $I^s\fL{s}f=f$.

On the other hand, with the definition of the fractional Laplacian we have
\begin{equation*}
	\Fourier(\fL{s}I^s f)(\xi) = (2\pi |\xi|)^{s} \Fourier(I^s f)(\xi) 
	= (2\pi |\xi|)^{s} (2\pi |\xi|)^{-s} \Fourier(f)(\xi) = \Fourier(f)(\xi),
\end{equation*}
which also justifies applying the fractional Laplacian to $I^s f$ since the inverse Fourier transform of $(2\pi |\xi|)^{s} \Fourier(I^s f)$ is obviously well defined.
\end{proof}

\begin{lemma}[An integration by parts formula for Riesz potentials]\label{lem:Rp integration by parts}
Let $f,g\in\mathcal{S}(\Rn)$ and $0<s<n$.
Then
\begin{equation*}
	\int_{\Rn} g I^{s} f = \int_{\Rn} I^{s} f \, g.
\end{equation*}
\end{lemma}

\begin{proof}
Analogous to Lemma \ref{lem:Rt integration by parts}.
\end{proof}

\begin{lemma}[Semigroup property of Riesz potentials]\label{lem:Rp semi group}
Let $f\in\mathcal{S}(\Rn)$ and $s_1,s_2>0$ with $s_1+s_2<n$.
Then
\begin{equation*}
	I^{s_1}I^{s_2}f = I^{s_1+s_2}f.
\end{equation*}
\end{lemma}

\begin{proof}
First note that the composition is well defined. With Theorem \ref{th:Rp convergence and Lq boundedness} we have $I^{s_2}f\in L^q(\Rn)$ for all $\frac{1}{q}= \frac{1}{p}- \frac{s}{n}$ where $1<p<\frac{n}{s_2}$ since $f\in L^p(\Rn)$ for all $1\leq p\leq\infty$. 
Because of 
\begin{equation*}
	q=\frac{pn}{n-s_2 p}\overset{p\rightarrow 1}{\longrightarrow} \frac{n}{n-s_2} < \frac{n}{n-(n-{s_1})}=\frac{n}{s_1},
\end{equation*}  
we find $1<p<\frac{n}{s_2}$ such that $q<\frac{n}{s_1}$. 
Then, again thanks to Theorem \ref{th:Rp convergence and Lq boundedness}, $I^{s_1}(I^{s_2}f)$ is well defined.
The rest of the proof is analogous to Lemma \ref{lem:fL semi group}.
\end{proof}

\newpage
\section{Trace Theorems}\label{sec:tt}
The goal of this section is to establish the blackbox estimates we use throughout the previous sections.
These complex estimates are listed and proven in Subsection \ref{subsec:tt - blackbox estimates}. 
The proofs depend on combining multiple smaller building blocks, which we collect in Subsections \ref{subsec:tt - building blocks}-\ref{subsec:tt - maximal function estimates}. 

The first and largest group of these building blocks originate from Triebel-Lizorkin or Besov-Lipschitz space characterizations via the generalized Poisson kernel, which are due to a recent result from Bui-Candy, \cite{BuiCan17}.
Providing the necessary background, in the first subsection we briefly introduce the framework of Tiebel-Lizorkin and Besov-Lipschitz spaces as well as the Bui-Candy result.
In the second subsection, we show that this result is indeed applicable to our situation and obtain a variety of space characterizations.
Since Triebel-Lizorkin or Besov-Lipschitz spaces include fractional Sobolev spaces, the $BMO$-space and Hölder spaces, we obtain a variety of estimates from these space characterizations. 
These are listed and proven in Subsection \ref{subsec:tt - building blocks}.
The second group of estimates originate from the theory of so-called square functions. 
This theory is closely related to the theory of maximal functions, which provides a third group of estimates. 
These estimates are shown in Subsection \ref{subsec:tt - square function estimates} and Subsection \ref{subsec:tt - maximal function estimates} respectively.
In their original form, many of these building blocks are $L^p$-estimates. 
From these we can easily obtain estimates on the much finer Lorentz scale via interpolation. 
Therefore, we introduce these Lorentz spaces and give the corresponding interpolation theorem in Subsection \ref{subsec:tt - lorentz spaces}.

\subsection{A short introduction to Triebel-Lizorkin and Besov-Lipschitz spaces}\label{subsec:tt - triebel lizorkin besov}
The theory of Triebel-Lizorkin and Lipschitz-Besov spaces provides a unified framework for a wide variety of important and frequently used function spaces. 
In this subsection, all functions, distributions etc. are defined on $\Rn$. 
Therefore, we omit $\Rn$ in the notation of the function spaces, i.e. instead of $A(\Rn)$ we write $A$.
First, we give a short list of spaces included in this framework and state some definitions, from which we will procede to the central result from \cite{BuiCan17}. 
At the end of this subsection, we collect some elementary embeddings and lifting properties.

The Triebel-Lizorkin spaces are denoted with $F^\alpha_{p,q}$, the Besov-Lipschitz spaces with $B^\alpha_{p,q}$, where $\alpha\in\R$, $0<p,q\leq\infty$. 
Additionally, there are their homogeneous counterparts, $\dot{F}^\alpha_{p,q}$ and $\dot{B}^\alpha_{p,q}$.
All the function spaces listed below can be identified with the corresponding Triebel-Lizorkin or Besov-Lipschitz spaces, the respective norms are equivalent.

\begin{tabular}{lcll | l}
$\mathscr{C}^s$	&$=$	& $B^s_{\infty,\infty}$	& if $0<s$					&\cite[Theorem 2.5.7.]{Tri83}\\
$\C^s$			&$=$	& $B^s_{\infty,\infty}$	& if $0<s\notin\Z$ 			&\cite[Proof of Theorem 2.5.7.]{Tri83}\\
$\Lambda^s_{p,q}$&$=$	& $B^s_{p,q}$		& if $0<s$, $1\leq p<\infty$, $1\leq q\leq\infty$	&\cite[Theorem 2.5.7.]{Tri83}\\
$W^s_p$			&$=$	& $B^s_{p,p}=F^s_{p,p}$	& if $0<s\notin\Z$, $1\leq p<\infty$	&\cite[Proof of Theorem 2.5.7.]{Tri83}\\
$h^p$			&$=$	& $F^0_{p,2}$ 			& if $0<p<\infty$			&\cite[Theorem 2.5.8.1.]{Tri83}\\
$\mathcal{H}^p$	&$=$	& $\dot{F}^0_{p,2}$		& if $0<p<\infty$			&\cite[Theorem 5.2.4.]{Tri83}\\
$L^p$			&$=$	& $F^0_{p,2}$			& if $1<p<\infty$			&\cite[Proof of Theorem 2.5.6.]{Tri83}\\
$L^p$			&$=$	& $\dot{F}^0_{p,2}$		& if $1<p<\infty$			&\cite[Theorem 5.2.3.1.(ii)]{Tri83}\\
$H^s_p$			&$=$	& $F^s_{p,2}$			& if $s\in\R$, $1<p<\infty$ 	&\cite[Theorem 2.5.6.(i)]{Tri83}\\
$W^m_p$			&$=$	& $F^m_{p,2}$			& if $m\in\N$, $1<p<\infty$ 	&\cite[Theorem 2.5.6.(ii)]{Tri83}\\	
$bmo$			&$=$	& $F^0_{\infty,2}$		&							&\cite[Theorem 2.5.8.2.]{Tri83}\\
$BMO$			&$=$	& $\dot{F}^0_{\infty,2}$	&							&\cite[Theorem 5.2.4.]{Tri83}	
\end{tabular}

Here, $\mathscr{C}^s$ are the Zygmund spaces, $\C^s$ the Hölder spaces, $\Lambda^s_{p,q}$ the Besov (or Lipschitz) spaces, $W^s_p$ the Slobodeckij spaces, $\mathcal{H}^p$ the Hardy spaces and $h^p$ the local Hardy spaces, $H^s_p$ the Bessel-potential spaces and $W^m_p$ the Sobolev spaces.
$BMO$ is the space of functions of bounded mean oscillation, the dual of $\mathcal{H}_1$ (see \cite{FeSte72}). 
$bmo$ is the inhomogeneous counterpart of $BMO$ and accordingly the dual of $h^1$ (see \cite{Go79}).
For definitions of these spaces see \cite[2.2.2.]{Tri83}.
Additionally, the definitions of $\C^s$, $W^s_p$, $H^s_p$ and $BMO$ are stated in the next subsection. 

Since we mainly use the homogeneous spaces in this work, $\dot{F}^\alpha_{p,q}$ and $\dot{B}^\alpha_{p,q}$, we only give their definition in this subsection and refer to \cite[2.3.1. and 2.3.4.]{Tri83} for the definition of the inhomogeneous Triebel-Lizorkin and Besov-Lipschitz spaces.
Working with the homogeneous spaces poses a difficulty though. While their inhomogeneous counterparts are defined as subsets of $\mathcal{S}'$, the space of tempered distributions, the elements of the homogeneous spaces are tempered distributions modulo polynomials, $\mathcal{S}'/\mathcal{P}$. The corresponding spaces can still be identified with for example $L^p$ though, because in each equivalence class there can only be one element with finite $L^p$-norm since for each polynomial $P\neq 0$ we have $\|P\|_{L^p}=\infty$ (cf. \cite[5.2.3. and 5.2.4.]{Tri83}). 

We now state the definition of the homogeneous Triebel-Lizorkin and Besov-Lipschitz spaces following \cite[Section 1]{BuiCan17} (who in turn refer to \cite[5.1.]{Tri83} and \cite[Chapter 3]{Pee76}).
In the following, for a function $\phi\colon \Rn\rightarrow\R$ and $j\in\Z$ we let $\phi_j(x)= 2^{jn}\phi(2^jx)$ denote the dyadic dilation of $\phi$.
Simultaneously, $\phi_t(x)=t^{-n}\phi(t^{-1}x)$ denotes the standard dilation for $t>0$.
The type of dilation we are using should always be clear, depending on whether we are in a discrete or continuous scenario.
For this definition we further use convolutions between Schwartz functions and tempered distributions.
For $f\in\mathcal{S}'$, the convolution $\phi\ast f\in\mathcal{S}'$ with $\phi\in\mathcal{S}$ is defined by
\begin{equation*}
	\langle \phi*f,\psi\rangle_{\mathcal{S}'\times\mathcal{S}} = \langle f,\tilde{\phi}\ast\psi\rangle_{\mathcal{S}'\times\mathcal{S}}, \qquad \text{for all }\psi\in\mathcal{S},
\end{equation*}
where $\tilde{\phi}(x)=\phi(-x)$.
While these convolutions are first defined only as a tempered distribution, we actually have $\phi*f\in\C^\infty$, see \cite[Theorem 2.3.20]{Gra14}.

\begin{definition}[Triebel-Lizorkin and Besov-Lipschitz spaces]\label{def:tt triebel lizorkin besov lipschitz spaces}
Fix $\varphi\in\mathcal{S}$ such that $\operatorname{supp} \widehat{\varphi}= \lbrace 1/2 \leq |\xi| \leq 2 \rbrace$ and for every $\xi\neq 0$
\begin{equation*}
	\sum_{j\in\Z} \widehat{\varphi}(2^{-j}\xi)\widehat{\varphi}(2^{-j}\xi) = 1.
\end{equation*}
Let $0<p,q\leq\infty$ and $\alpha\in\R$.
The homogeneous Besov-Lipschitz spaces are defined as
\begin{equation}\label{eq:tt def Bspq}
	\dot{B}^\alpha_{p,q} = 	\bigg\lbrace f\in \mathcal{S}'/\mathcal{P} \,;\,
								\|f\|_{\dot{B}^\alpha_{p,q}}= 
									\bigg(\sum_{j\in\Z}\left(2^{j\alpha} \|\varphi_j * f\|_{L^p}\right)^q\bigg)^\frac{1}{q}<\infty 
							\bigg\rbrace.
\end{equation}
For $p\neq\infty$, the homogeneous Triebel-Lizorkin spaces are defined as
\begin{equation}\label{eq:tt def Fspq}
	\dot{F}^\alpha_{p,q} = 	\bigg\lbrace f\in \mathcal{S}'/\mathcal{P} \,;\,
								\|f\|_{\dot{F}^\alpha_{p,q}}= 
									\bigg\|\bigg(\sum_{j\in\Z}\left(2^{j\alpha} |\varphi_j * f|\right)^q\bigg)^\frac{1}{q}\bigg\|_{L^p}<\infty 
							 \bigg\rbrace.
\end{equation}
For $p=\infty$, the corresponding Triebel-Lizorkin spaces are defined as
\begin{equation}\label{eq:tt def Fsinftyq}
	\dot{F}^\alpha_{\infty,q} =	\bigg\lbrace f\in \mathcal{S}'/\mathcal{P} \,;\,
									\|f\|_{\dot{F}^\alpha_{\infty,q}}= 
										\sup_Q \bigg(\frac{1}{|Q|}\int_Q \sum_{j\geq-l(Q)}\left(2^{j\alpha} |\varphi_j * f|\right)^q dx \bigg)^\frac{1}{q}<\infty 
								 \bigg\rbrace.
\end{equation}
where the supremum is over all dyadic cubes $Q$ and $l(Q)= \log_2(\text{ side length of $Q$ })$.
\end{definition}

For $q=\infty$, the interpretation is 
\begin{equation*}
	\|f\|_{\dot{B}^\alpha_{p,\infty}}= \sup_{j\in\Z}\left(2^{j\alpha} \|\varphi_j * f\|_{L^p}\right),	\qquad\\
	\|f\|_{\dot{F}^\alpha_{p,\infty}}= \|\sup_{j\in\Z}\left(2^{j\alpha} |\varphi_j * f|\right)\|_{L^p}
\end{equation*}
and for the special case of the Triebel-Lizorkin space $\dot{F}^\alpha_{\infty,\infty}$ 
\begin{equation*}
	\|f\|_{\dot{F}^\alpha_{\infty,\infty}}= \sup_Q \sup_{j\geq-l(Q)}\frac{1}{|Q|}\int_Q 2^{j\alpha} |\varphi_j * f| dx.
\end{equation*}

Further note that the quasi-norms applied to polynomials indeed yield zero since the convolution of any polynomial $\rho\in\mathcal{P}$ with any of the $\varphi_j$ is zero. This can be seen on the side of the Fourier transforms since the support of a polynomials' Fourier transform is always $\lbrace 0 \rbrace$.
We have
\begin{equation*}
	\langle\varphi_j\ast\rho, \psi\rangle_{\mathcal{S}'\times\mathcal{S}} 
	= \langle\rho, \tilde{\varphi_j}\ast\psi\rangle_{\mathcal{S}'\times\mathcal{S}}
	= \langle \Fourier^{-1}\rho, \Fourier\tilde{\varphi_j} \Fourier\psi\rangle_{\mathcal{S}'\times\mathcal{S}} = 0
\end{equation*}
for all $\psi\in\mathcal{S}$.

It is well established that this definition is independent of the choice of the kernel $\varphi$ in the sense that replacing $\varphi$ by another kernel satisfying the same conditions yields equivalent norms, see \cite[Theorem 5.1.5.]{Tri83}. Sustaining the independence of the chosen kernel is also the reason for the deviating definition of $\dot{F}^\alpha_{\infty,q}$. 
The dependence on the kernel when choosing the natural extension to $p=\infty$ has been remarked in \cite[2.3.1.4.]{Tri83}.
For a discussion of the chosen definition for $\dot{F}^\alpha_{\infty,q}$ and the independence of the kernel, see \cite[Section 5]{FraJa90}.
Originally, these definitions were formulated for a a family of kernels that were not necessarily all dilations of one single kernel - that is just the easiest way to realize these conditions (cf. \cite[Chapter 3, pp. 47-49]{Pee76} and \cite[5.1.3]{Tri83}).

The definition via a single kernel and its dyadic dilations is more interesting to us though, since it naturally leads to a continuous version of the Triebel-Lizorkin and Besov-Lipschitz quasi-norms.
The summation over the convolutions with dyadic dilation
s can be replaced by an integral with the dilation parameter as continuous integration variable.
For $\alpha\in\R$, $0<p<\infty$, and $0<q\leq\infty$ we have
\begin{equation*}
	\|f\|_{\dot{B}^\alpha_{p,q}} 
	\approx \bigg(\int_0^\infty (t^{-\alpha}\|\varphi_t * f\|_{L^p})^q\frac{dt}{t}\bigg)^\frac{1}{q},
\end{equation*}
\begin{equation*}
	\|f\|_{\dot{F}^\alpha_{p,q}} 
	\approx \bigg\|\bigg( \int_0^\infty(t^{-\alpha}|\varphi_t * f|)^q\frac{dt}{t}\bigg)^\frac{1}{q}\bigg\|_{L^p}
\end{equation*}
for all $f\in\mathcal{S}'/\mathcal{P}$, see the introduction of \cite{BuiPaTai96} and Section 6(d) for the arguments.
Additionally, Theorem 1 in \cite{BuiTai00} covers Triebel-Lizorkin spaces for $p=\infty$,
\begin{equation*}
	\|f\|_{\dot{F}^\alpha_{\infty,q}} \approx 
	\sup_{x\in\Rn,r>0} \bigg\lbrace \frac{1}{|B(x,r)|}\int_{B(x,r)}\int_0^r(t^{-\alpha}|\varphi_t*f|)^q\frac{dt}{t}\,dy\bigg\rbrace^\frac{1}{q}
\end{equation*}
for all $f\in\mathcal{S}'/\mathcal{P}$.

Over the years, looser conditions for kernels $\psi\in\mathcal{S}$, under which replacing $\varphi$ with $\psi$ in (\ref{eq:tt def Bspq}), (\ref{eq:tt def Fspq}) and (\ref{eq:tt def Fsinftyq}) yields equivalent norms, were developed (see \cite{Tri83}, \cite{BuiPaTai96}, \cite{BuiTai00}).
Of key interest to us though are the results in \cite{BuiCan17}. 
There, suitable conditions for non-smooth kernels $\psi\notin\mathcal{S}$ are established, specifically with fractional derivatives of the Poisson kernel in mind. 

One of the problems arising when working with a non-smooth kernel $\psi\notin\mathcal{S}$ is, that it is not possible to define the convolution $\psi*f$ for arbitrary $f\in\mathcal{S}'$.
To deal with this problem, distributions of bounded growth are used in \cite{BuiCan17}.

\begin{definition}[Distributions of growth $l$]
\begin{enumerate}[label=(\roman*)]
\item A tempered distribution $f\in\mathcal{S}'$ is of growth $l\geq 0$ if for any $\phi\in\mathcal{S}$ we have $\phi*f = \mathcal{O}(|x|^l)$ as $|x|\rightarrow\infty$.
\item  Assume $f$ is a distribution of growth $l$.
Then if $(1+|\cdot|)^l\psi\in L^1$, the convolution $\psi*f\in\mathcal{S}'$ is given as 
\begin{equation*}
	\langle\psi*f,\phi\rangle_{\mathcal{S}'\times\mathcal{S}}  = \int_\Rn \psi(x)(\widetilde{\phi}*f)(x)dx, \quad \phi\in\text{for all }\mathcal{S},
\end{equation*}
where $\widetilde{\phi}(x)=\phi(-x)$.
\end{enumerate}
\end{definition}

While we are able to define the convolutions themselves by restricting the admissible distributions, as further unspecified tempered distributions the convolutions still make no sense in (\ref{eq:tt def Bspq}) or (\ref{eq:tt def Fspq}).
In the main result of \cite{BuiCan17}, additional conditions ensure that these convolutions are continuous functions.
Before we discuss this result, we give an example for distributions of growth $l\in\N$ - polynomials. 

\begin{lemma}\label{lem:tt polynomials as distirbutions of growth l}
Let $\rho\in\mathcal{P}$, $l=\deg(\rho)$. 
Then $\rho$ is a distribution of growth exactly $l$, meaning that $\rho$ is not a distribution of growth $\tilde{l}$ for any $\tilde{l}<l$.
\end{lemma}

\begin{proof}
We first show that $\rho$ is indeed a distribution of growth $l$. 
There exists $C_1>0$ such that $|\rho(x)|\leq C_1(1+|x|)^l$ for all $x\in\Rn$.
Now let $\phi\in\mathcal{S}$. 
Then there exists $C_2>0$ such that $|\phi(x)|\leq C_2 (1+|x|)^{-n-1-l}$ for all $x\in\Rn$.
Combining these two estimates, since for $y\notin B_{x}(x)$ we have $|x-y|\geq |y|/2$,
\begin{mathex}
	|\phi*\rho(x)|
	&\leq&	\int_\Rn |\phi(x-y)\rho(y)| \,dy \\
	&=&		\int_{B_{\frac{x}{2}}(x)} |\phi(x-y)\rho(y)| \,dy 
			+ \int_{\Rn\setminus B_{\frac{x}{2}}(x)} |\phi(x-y)\rho(y)| \,dy\\
	&\leq&	C_1\|\phi\|_{L^1}(1+2|x|)^l 
			+ \int_{\Rn\setminus B_{\frac{x}{2}}(x)} C_1(1+|y|)^{l}\,C_2\left(1+\frac{|y|}{2}\right)^{-n-1-l} \,dy\\
	&\leq&	C_1\|\phi\|_{L^1}(1+2|x|)^l + 2^{n+1+l} C_1 C_2 \int_{\Rn} (1+|y|)^{-n-1}\,dy\\
	&=&		\mathcal{O}(|x|^l)
	\quad \text{as } |x|\rightarrow \infty.
\end{mathex}
Since the growth of this convolution is bounded, it is easy to see with Fubini's theorem that this convolution in the standard integral sense coincides with the convolution in the distributional sense, see \cite[Definition 2.3.13]{Gra14} for the definition.
Therefore, $\rho$ is a distribution of growth $l$. 

Now assume $\tilde{l}<l$.
In order to show that $\rho$ is no distribution of growth $\tilde{l}$, we construct a ray along which the polynomial grows according to its degree. We then examine how the convolution with a compactly supported, positive kernel behaves along that ray.
Let $\rho= \sum_{|\gamma|\leq l} c_\gamma x^\gamma$ where the sum is over multiindices $\gamma\in\Nz^n$.
Denote with $\tilde{\gamma}$ the largest element of $\Gamma_l\coloneqq \lbrace \gamma\in\Nz^n \,;\, |\gamma|= l \text{ and } c_\gamma\neq 0\rbrace$ in terms of the lexicographic order.
Set
\begin{equation*}
	C\coloneqq \frac{1}{|c_{\tilde{\gamma}}|}\sum_{|\gamma|= l} |c_\gamma|,
\end{equation*}
\begin{equation*}
	v\coloneqq \left( C^{(l+1)^{n-i}} \right)_{i=1,\ldots,n} = \left(C^{(l+1)^{n-1}},C^{(l+1)^{n-2}},\ldots,C^{l+1},C\right)\in\Rn.
\end{equation*}
The vector $v$ will be the direction of the aforementioned ray. 
We first observe that 
\begin{equation}\label{eq:tt polynomials gamma=l}
	c_{\tilde{\gamma}}\left|v^{\tilde{\gamma}}\right|-\sum_{|\gamma|=l, \gamma\neq\tilde{\gamma}}|c_\gamma| \left|v^\gamma\right|
	\geq \frac{|c_{\tilde{\gamma}}|}{C} \left|v^{\tilde{\gamma}}\right|.
\end{equation}
This is clear if $\tilde{\gamma}=le_n$, because then $\Gamma_l=\lbrace le_n\rbrace$.
Therefore, assume $\tilde{\gamma}\neq le_n$.
For $\gamma\in\Gamma_l$, $\gamma\neq \tilde{\gamma}$, let $i\in\lbrace 1,\ldots,n\rbrace$ be minimal with the property that $\tilde{\gamma}_i>\gamma_i$. 
Since $C\geq 1$,
\begin{equation*}
	\frac{\left|v^{\tilde{\gamma}}\right|}{\left|v^{\gamma}\right|}
	= 
	\frac{\left(C^{(l+1)^{n-i}}\right)^{\tilde{\gamma}_i-\gamma_i}\left|(v_{i+1,\ldots,n})^{\tilde{\gamma}_{i+1,\ldots,n}}\right|}{\left|(v_{i+1,\ldots,n})^{\gamma_{i+1,\ldots,n}}\right|} 
	\geq 
	\frac{C^{(l+1)^{n-i}}\cdot 1}{\left(C^{(l+1)^{n-(i+1)}}\right)^l}
	= 
	C^{(l+1)^{n-(i+1)}} \geq C.
\end{equation*}
With this estimate and the definition of $C$ we obtain
\begin{equation*}
	|c_{\tilde{\gamma}}|\left|v^{\tilde{\gamma}}\right|
	= \frac{1}{C}\sum_{\gamma\in\Gamma_l} |c_\gamma|\left|v^{\tilde{\gamma}}\right|
	= \frac{|c_{\tilde{\gamma}}|}{C} \left|v^{\tilde{\gamma}}\right| + \sum_{|\gamma|=l, \gamma\neq\tilde{\gamma}}|c_\gamma| \frac{\left|v^{\tilde{\gamma}}\right|}{C}
	\geq \frac{|c_{\tilde{\gamma}}|}{C} \left|v^{\tilde{\gamma}}\right| + \sum_{|\gamma|=l, \gamma\neq\tilde{\gamma}}|c_\gamma|  \left|v^\gamma\right|
\end{equation*}
and therefore prove (\ref{eq:tt polynomials gamma=l}).
Thus, for $s\geq 1$
\begin{equation*}
	|\rho(sv)| 
	\geq  s^l\left(c_{\tilde{\gamma}}\left|v^{\tilde{\gamma}}\right|
	-\sum_{|\gamma|=l, \gamma\neq\tilde{\gamma}}|c_\gamma| \left|v^\gamma\right|\right)
	-\sum_{|\gamma|<l}|c_\gamma| \left|(sv)^\gamma\right|	
	\geq s^l \frac{|c_{\tilde{\gamma}}|}{C} \left|v^{\tilde{\gamma}}\right|
	-s^{l-1}\sum_{|\gamma|<l}|c_\gamma| \left|(v)^\gamma\right|
\end{equation*}
and apparently there is an $R_1>0$, such that for all $s>R_1$
\begin{equation*}
	|\rho(sv)|	\geq s^l \frac{|c_{\tilde{\gamma}}|}{2C} \left|v^{\tilde{\gamma}}\right|
				\geq s^l \frac{|c_{\tilde{\gamma}}|}{2C}.
\end{equation*}
Thus, we have shown that $\rho$ grows sufficiently fast in the chosen direction.
We also need to transfer this growth to neighborhoods along the ray though. 
Since the derivatives of $\rho$ are polynomials of degree $l-1$, there exists $R_2>0$ such that 
\begin{equation*}
	|\nabla \rho(x)|\leq  \frac{|c_{\tilde{\gamma}}|}{4C |v|^l} (1+|x|)^{l}
\end{equation*}
for all $x\in\Rn$, $|x|\geq R_2$.
Thus, for $x\in B_1(sv)$ and $s>R_3 \coloneqq \max\lbrace R_1, |v|^{-1}(R_2+1) \rbrace$ 
\begin{equation}
	|\rho(x)|\geq |\rho(sv)|-\|\nabla \rho\|_{L^\infty(B_1(sv))} 
	\geq s^l \frac{|c_{\tilde{\gamma}}|}{2C} - \frac{|c_{\tilde{\gamma}}|}{4C |v|^l}(1+(|sv|-1))^l 
	= s^l \frac{|c_{\tilde{\gamma}}|}{4C}.
\end{equation}

Let $\phi\in\mathcal{S}$ with $\operatorname{supp}(\phi)= \overline{B_1(0)}$ and $\phi>0$ on $B_1(0)$.
Then, for $s > R_3$,
\begin{equation}
	|(\phi * \rho)(sv)| = \int_{B_1(sv)} \phi(sv-y) |\rho(y)| \,dy \geq \frac{|c_{\tilde{\gamma}}|\|\phi\|_{L^1}}{4C|v|^l} |sv|^l. 
\end{equation}
Therefore, $\rho$ can not be a distribution of growth $\tilde{l}$.
\end{proof}

We now turn to the main result of \cite{BuiCan17} - the characterization of Triebel-Lizorkin and Besov-Lipschitz spaces with Poisson like kernels.
There are three main conditions the kernels $\psi\in L^1(\Rn)$ have to satisfy in order to obtain the norm equivalencies, reading as follows.
The parameters $\Lambda\geq 0$ and $m,r\in\R$ are chosen depending on the characterized space, $[\cdot]$ denotes the integer part of the inserted number.
\begin{description}
\item[(C1)]	\textit{(Cancellation)} 
	$\Fourier(\psi)\in \C^{n+1+[\Lambda]}(\Rn\setminus\lbrace 0\rbrace)$ such that for every $|\kappa|\leq n+1+[\Lambda]$
	\begin{equation*}
		\partial_{\kappa}\Fourier(\psi) = \mathcal{O}(|\xi|^{r-|\kappa|}) \qquad \text{as } |\xi|\rightarrow 0.
	\end{equation*} 

\item[(C2)]\textit{(Tauberian condition)} 
	For every direction $\xi\in S^{n-1}$ there exist $a,b\in\R$ (depending on $\xi$) with $0<2a\leq b$ such that for every $a < t < b$
	\begin{equation*}
		|\Fourier(\psi)(t\xi)|>0.
	\end{equation*}

\item[(C3)] \textit{(Smoothness)} 
	$\Fourier(\psi)\in \C^{n+1+[\Lambda]}(\Rn\setminus\lbrace 0\rbrace)$ such that for every $|\kappa|\leq n+1+[\Lambda]$ 
	\begin{equation*}
		\partial_{\kappa}\Fourier(\psi) = \mathcal{O}(|\xi|^{-n-m}) \qquad \text{as } |\xi|\rightarrow\infty.
	\end{equation*}
\end{description}
Following \cite{BuiCan17}, we split the characterization into two theorems - one for each estimate required for the norm equivalency, stressing the different sets of assumptions.
Note that, since for our method we only need to estimate upwards against special cases of the Triebel-Lizorkin and Besov-Lipschitz quasi-norms, we technically only need the first theorem.
We still include the second one for completeness though.

\begin{theorem}[Theorem 1.1 of \cite{BuiCan17}]\label{th:tt BuiCandy main result 1}
Let $\alpha\in\R$ and $0<p,q\leq\infty$. 
Let $l\geq 0$ with $l> \alpha - \frac{n}{p}$.
Assume $(1+|\cdot |)^l\psi\in L^1$.
\begin{enumerate}[label=(\roman*)]
\item	Let $f\in \dot{B}^\alpha_{p,q}$.
		Assume that $\psi$ satisfies \textbf{(C1)} and \textbf{(C3)} for $\Lambda=\frac{n}{p}$, $r>\alpha$ and $m>\Lambda-\alpha$.
		Then there exists a polynomial $\rho\in\mathcal{P}$ such that $f-\rho$ is a distribution of growth $l$ and we have
		\begin{equation*}
			\bigg(\sum_{j\in\Z}(2^{j\alpha}\|\psi_j*(f-\rho)\|_{L^p})^q\bigg)^\frac{1}{q} 
			\lesssim \|f\|_{\dot{B}^\alpha_{p,q}},
		\end{equation*}
		and the continuous version
		\begin{equation*}
			\bigg(\int_0^\infty (t^{-\alpha}\|\psi_t * (f-\rho)\|_{L^p})^q\frac{dt}{t}\bigg)^\frac{1}{q}
			\lesssim \|f\|_{\dot{B}^\alpha_{p,q}}.
		\end{equation*}
\item	Let $f\in \dot{F}^\alpha_{p,q}$.
		Assume that $\psi$ satisfies \textbf{(C1)} and \textbf{(C3)} for 
		$\Lambda=\max\left\lbrace\frac{n}{p},\frac{n}{q}\right\rbrace$ (with $\Lambda=n$ when $p=q=\infty$),
		$r>\alpha$ and $m>\Lambda-\alpha$.
		Then there exists a polynomial $\rho\in\mathcal{P}$ such that $f-\rho$ is a distribution of growth $l$ and 
		for $p<\infty$ we have
		\begin{equation*}
			\bigg\|\bigg(\sum_{j\in\Z}(2^{j\alpha}|\psi_j*(f-\rho)|)^q\bigg)^\frac{1}{q} \bigg\|_{L^p}
			\lesssim \|f\|_{\dot{F}^\alpha_{p,q}},
		\end{equation*}
		and the continuous version	 
		\begin{equation*}
			\bigg\|\bigg( \int_0^\infty(t^{-\alpha}|\psi_t * (f-\rho)|)^q\frac{dt}{t}\bigg)^\frac{1}{q}\bigg\|_{L^p}
			\lesssim \|f\|_{\dot{F}^\alpha_{p,q}}.
		\end{equation*}
		For $p=\infty$ we have
		\begin{equation*}
			\sup_Q \bigg(\frac{1}{|Q|}\int_Q \sum_{j\geq-l(Q)}\left(2^{j\alpha} |\psi_j * (f-\rho)|\right)^q dx \bigg)^\frac{1}{q}
			\lesssim \|f\|_{\dot{F}^\alpha_{\infty,q}},
		\end{equation*}
		and the continuous version
		\begin{equation*}
			\sup_{x\in\Rn,r>0} \bigg\lbrace \frac{1}{|B(x,r)|}\int_{T(B(x,r))}(t^{-\alpha}|\psi_t*(f-\rho)|)^q\frac{dt}{t}\,dy\bigg\rbrace^\frac{1}{q}
			\lesssim \|f\|_{\dot{F}^\alpha_{\infty,q}},
		\end{equation*}
		where $T(B(x,r))=\lbrace (y,t)\in\Rnpp\, :\,|y-x|<r-t\rbrace$ is the ``tent'' over the ball $B(x,r)$.
\end{enumerate}
\end{theorem}

\begin{proof}
See \cite{BuiCan17}. 
Theorem 1.1 being formulated for the Peetre maximal function applied to the convolutions poses no problem since the maximal function majorizes the convolution, see Remark 1.1(iv) in \cite{BuiCan17}.
The continuous version of the result with integrals instead of sums is confirmed but left for the reader to prove, see the end of the introduction in \cite{BuiCan17}. 
For details regarding the precise formulation and modifications of the proof they refer to \cite{BuiPaTai96}, \cite{BuiPaTai97} and \cite{BuiTai00}.

Note that for the continuous version of the $\dot{F}^\alpha_{\infty,q}$-characterization we integrate over tents while in \cite{BuiTai00} they integrate over cylinders. 
But since $|B(x,r)|= 2^{-n}|B(x,2r)|$ and because of the inclusions $T(B(x,r))\subset B(x,r)\times[0,r] \subset T(B(x,2r))$ we have 
\begin{mathex}
	\sup_{x\in\Rn,r>0} \bigg\lbrace \frac{1}{|B(x,r)|}\int_{T(B(x,r))}(t^{-\alpha}|\varphi_t*f|)^q\frac{dt}{t}\,dy\bigg\rbrace^\frac{1}{q}\\
	\leq \sup_{x\in\Rn,r>0} \bigg\lbrace \frac{1}{|B(x,r)|}\int_{B(x,r)\times[0,r]}(t^{-\alpha}|\varphi_t*f|)^q\frac{dt}{t}\,dy\bigg\rbrace^\frac{1}{q}\\
	\leq \sup_{x\in\Rn,r>0} \bigg\lbrace \frac{1}{2^{-n}|B(x,2r)|}\int_{T(B(x,2r))}(t^{-\alpha}|\varphi_t*f|)^q\frac{dt}{t}\,dy\bigg\rbrace^\frac{1}{q}
\end{mathex}
and therefore see that both versions are equivalent.
\end{proof}

Of course, the polynomials depending on $f$ appearing in the estimates are undesirable. 
We will deal with those later by showing that for our specific kernels we may always choose $\rho=0$, see Theorem \ref{th:tt triebel lizorkin besov char}.

\begin{theorem}[Theorem 1.3 of \cite{BuiCan17}]\label{th:tt BuiCandy main result 2}
Let $\alpha\in\R$, $0<p,q\leq\infty$.
Assume $f$ is a distribution of growth $l\geq 0$.
Suppose $(1+|\cdot |)^l\psi\in L^1$ and $\psi$ satisfies \textbf{(C2)} as well as \textbf{(C3)} for every $m\in\R$ and some $\Lambda\geq l$.
Then for every $j\in\Z$ the convolution $\psi_j*f$ is a continuous function.
Moreover,
\begin{enumerate}[label=(\roman*)]
\item	if \textbf{(C3)} holds for $\Lambda = \max\left\lbrace l,\frac{n}{p}\right\rbrace$ and every $m\in\R$ respectively , then
		\begin{equation*}
			\|f\|_{\dot{B}^\alpha_{p,q}}
			\lesssim \bigg(\sum_{j\in\Z}(2^{j\alpha}\|\psi_j*f\|_{L^p})^q\bigg)^\frac{1}{q},
		\end{equation*}
		\begin{equation*}
			\|f\|_{\dot{B}^\alpha_{p,q}}
			\lesssim  \bigg(\int_0^\infty (t^{-\alpha}\|\psi_t * f\|_{L^p})^q\frac{dt}{t}\bigg)^\frac{1}{q}.
		\end{equation*}
\item	if \textbf{(C3)} holds for $\Lambda = \max\left\lbrace l,\frac{n}{p},\frac{n}{q}\right\rbrace$ and every $m\in\R$ respectively, then for $p<\infty$
		\begin{equation*}
			\|f\|_{\dot{F}^\alpha_{p,q}}
			\lesssim \bigg\|\bigg(\sum_{j\in\Z}(2^{j\alpha}|\psi_j*f|)^q\bigg)^\frac{1}{q} \bigg\|_{L^p},
		\end{equation*}
		\begin{equation*}
			\|f\|_{\dot{F}^\alpha_{p,q}}
			\lesssim  \bigg\|\bigg( \int_0^\infty(t^{-\alpha}|\psi_t * f|)^q\frac{dt}{t}\bigg)^\frac{1}{q}\bigg\|_{L^p}.	
		\end{equation*}
		In the case $p=\infty$ we have
		\begin{equation*}
			\|f\|_{\dot{F}^\alpha_{\infty,q}} \lesssim 
			\sup_Q \bigg(\frac{1}{|Q|}\int_Q \sum_{j\geq-l(Q)}\left(2^{j\alpha} |\psi_j * f|\right)^q dx \bigg)^\frac{1}{q},
		\end{equation*}
		\begin{equation*}
			\|f\|_{\dot{F}^\alpha_{\infty,q}} \lesssim 
			\sup_{x\in\Rn,r>0} \bigg\lbrace \frac{1}{|B(x,r)|}\int_{T(B(x,r))}(t^{-\alpha}|\psi_t*f|)^q\frac{dt}{t}\,dy\bigg\rbrace^\frac{1}{q},
		\end{equation*}
		where again $T(B(x,r))=\lbrace (y,t)\in\Rnpp\, :\,|y-x|<r-t\rbrace$ is the ``tent'' over $B(x,r)$.
\end{enumerate}
\end{theorem}

\begin{proof}
See \cite{BuiCan17}. 
The continuous version of the result with integrals instead of sums is confirmed but left for the reader to prove, see the end of the introduction in \cite{BuiCan17}. 
Because the mentioned ``standard'' Tauberian condition is weaker than \textbf{(C2)}, we do not need to make adjustments regarding the assumptions.
For details regarding the precise formulation and modifications of the proof they refer to \cite{BuiPaTai96}, \cite{BuiPaTai97} and \cite{BuiTai00}.

While in Theorem 1.3 it says ``Suppose $(1+|\cdot |)^l\psi\in L^1$ satisfies the Tauberian condition \textbf{(C2)} and that for every $m\in\R$, the smoothness condition \textbf{(C3)} holds with $\Lambda\geq l$'', considering the results such as Theorem 5.1 leading up to Theorem 1.3, the more precise, less misleading formulation seems to be ``Suppose $(1+|\cdot |)^l\psi\in L^1$ and $\psi$ satisfies [...]''.

With the same reasoning as in the proof of Theorem \ref{th:tt BuiCandy main result 1} we integrate over tents instead of cylinders for the $\dot{F}^\alpha_{\infty,q}$-characterization.
\end{proof}


Concluding this subsection, we list two further results for Triebel-Lizorkin and Besov-Lipschitz spaces which we use throughout this work.

The first result we cover concerns some embeddings between different Triebel-Lizorkin and Besov-Lipschitz spaces. 
In the following proposition, for two quasi-normed spaces $A_1$ and $A_2$,  $A_1\subset A_2$ means that $A_1$ is continuously embedded in $A_2$. 
More specifically, we have $\|a\|_{A_2}\lesssim \|a\|_{A_1}$ for all $a\in A_1$.
This result will be helpful when combining multiple space characterizations for some of the more complex estimates in Section \ref{subsec:tt - blackbox estimates}.

\begin{proposition}[Elementary embeddings]\label{prop:tt elementary embeddings for trilibeli}
\begin{enumerate}[label=(\roman*)]
\item	Let $0<q_0\leq q_1 \leq \infty$ and $\alpha\in\R$. Then
		\begin{equation*}
			B^\alpha_{p,q_0} \subset B^\alpha_{p,q_1} 						\quad\text{and}\quad 
			\dot{B}^\alpha_{p,q_0} \subset \dot{B}^\alpha_{p,q_1} 			\quad \text{if } 0<p\leq\infty
		\end{equation*}
		as well as 
		\begin{equation*}
			F^\alpha_{p,q_0} \subset F^\alpha_{p,q_1} 						\quad\text{and}\quad 
			\dot{F}^\alpha_{p,q_0} \subset \dot{F}^\alpha_{p,q_1} 			\quad \text{if } 0<p<\infty.
		\end{equation*}
\item	Let $0<q_0 \leq \infty$, $0<q_1 \leq \infty$, $\alpha\in\R$ and $\epsilon>0$. Then
		\begin{equation*}
			B^{\alpha-\epsilon}_{p,q_0} \subset B^\alpha_{p,q_1} 						\quad\text{and}\quad 
			\dot{B}^{\alpha-\epsilon}_{p,q_0} \subset \dot{B}^\alpha_{p,q_1} 			\quad \text{if } 0<p\leq\infty
		\end{equation*}
		and
		\begin{equation*}
			F^{\alpha-\epsilon}_{p,q_0} \subset F^\alpha_{p,q_1} 						\quad\text{and}\quad 
			\dot{F}^{\alpha-\epsilon}_{p,q_0} \subset \dot{F}^\alpha_{p,q_1} 			\quad \text{if } 0<p<\infty.
		\end{equation*}
\item	Let $0<q\leq\infty$, $0<p<\infty$ and $\alpha\in\R$. Then
		\begin{equation*}
			B^\alpha_{p,\min\lbrace p,q\rbrace} \subset F^\alpha_{p,q} \subset B^\alpha_{p,\max\lbrace p,q\rbrace}
			\quad\text{and}\quad
			\dot{B}^\alpha_{p,\min\lbrace p,q\rbrace} \subset \dot{F}^\alpha_{p,q} \subset \dot{B}^\alpha_{p,\max\lbrace p,q\rbrace}
		\end{equation*}
\end{enumerate}
\end{proposition}

\begin{proof}
See \cite[Proposition 2.3.2.2]{Tri83}. The proofs for the homogeneous spaces are analogous.
\end{proof}

The second result we cover are the lifting properties of the Triebel-Lizorkin and Besov-Lipschitz spaces.
These lifting properties open up more possible norm equivalencies.
More specifically, based on the identifications of $\dot{F}^\alpha_{p,q}$, $F^\alpha_{p,q}$, $\dot{B}^\alpha_{p,q}$ or $B^\alpha_{p,q}$ with spaces we already know, we can easily formulate equivalent, more intuitive norms for a wider range of Triebel-Lizorkin and Besov-Lipschitz spaces.
For example, the identification of Triebel-Lizorkin spaces with $H^s_p$ and therefore $W^m_p$ listed in the table at the start of the subsection is due to the following theorem.

Before stating the Theorem, we introduce the lifting operators, which will allow us to mediate between different Triebel-Lizorkin and Lipschitz-Besov spaces.
Let $\sigma\in\R$, then 
\begin{mathex}
	J^\sigma f 		&\coloneqq \Fourier^{-1}\left((1+|\xi|^2)^\frac{\sigma}{2}\Fourier(f)(\xi)\right) 
							&\qquad \text{for }f\in\mathcal{S}',\\
	\dot{J}^\sigma f	&\coloneqq \Fourier^{-1}(|\xi|^\sigma\Fourier(f)(\xi) 	
							&\qquad \text{for }f\in\mathcal{Z}'=\mathcal{S}'/\mathcal{P},
\end{mathex} 
where $\mathcal{Z}'$ is the topological dual of 
	$\mathcal{Z}=\lbrace \phi\in\mathcal{S};\quad \partial_\alpha\Fourier(\phi)(0) \quad\text{for all }\alpha\in\Nz^n\rbrace$.
Defining the Fourier transform on $\mathcal{Z}'$ via $\mathcal{Z}$ analogous to the way the Fourier transform on $\mathcal{S}'$ is defined via $\mathcal{S}$, the above definition of $\dot{J}^\sigma$ makes sense.
For more details we refer to \cite[5.1]{Tri83}.
We do not go into more depth regarding the spaces $\mathcal{Z}$ and $\mathcal{Z}'$ since under the circumstances we work with the expressions will already be well defined in the standard sense.

Obviously, up to a constant $\const_\sigma$ we have 
\begin{equation*}
	\dot{J}^\sigma f = \left\lbrace\begin{array}{lll}
								\const_\sigma I^\sigma f 		\qquad	&\text{if } -n<\sigma<0,	\\
								\const_\sigma \fL{\sigma} f 		\qquad	&\text{if } \sigma>0	,
						\end{array}\right.
\end{equation*}
for example for $f\in\mathcal{S}$. 

\begin{theorem}[Lifting property of Triebel-Lizorkin and Besov-Lipschitz spaces]\label{th:tt lifting property of trlibeli spaces}
Let $\alpha,\sigma \in\R$ and $0<q\leq\infty$.
\begin{enumerate}[label=(\roman*)]
\item	If $0<p\leq\infty$, then $J^\sigma$ maps $B^\alpha_{p,q}$ isomorphically onto $B^{\alpha-\sigma}_{p,q}$ while $\dot{J}^\sigma$ maps $\dot{B}^\alpha_{p,q}$ isomorphically onto $\dot{B}^{\alpha-\sigma}_{p,q}$. 
In particular, $\|J^\sigma f\|_{B^{\alpha-\sigma}_{p,q}}$ and $\|\dot{J}^\sigma f\|_{\dot{B}^{\alpha-\sigma}_{p,q}}$ are equivalent quasi-norms on $B^\alpha_{p,q}$ and $\dot{B}^\alpha_{p,q}$ respectively.
Moreover, for $m\in\Nz$
\begin{equation*}
	\sum_{|\kappa|\leq m} \|\partial_\kappa f \|_{B^{\alpha-m}_{p,q}} 
	\quad\text{and}\quad
	\sum_{|\kappa|= m} \|\partial_\kappa f \|_{\dot{B}^{\alpha-m}_{p,q}}
\end{equation*}
are equivalent norms on $B^\alpha_{p,q}$ and $\dot{B}^\alpha_{p,q}$ respectively.

\item	If $0<p<\infty$, then $J^\sigma$ maps $F^\alpha_{p,q}$ isomorphically onto $F^{\alpha-\sigma}_{p,q}$ while $\dot{J}^\sigma$ maps $\dot{F}^\alpha_{p,q}$ isomorphically onto $\dot{F}^{\alpha-\sigma}_{p,q}$.
In particular, $\|J^\sigma f\|_{F^{\alpha-\sigma}_{p,q}}$ and $\|\dot{J}^\sigma f\|_{\dot{F}^{\alpha-\sigma}_{p,q}}$ are equivalent quasi-norms on $F^\alpha_{p,q}$ and $\dot{F}^\alpha_{p,q}$ respectively.
Moreover, for $m\in\Nz$
\begin{equation*}
	\sum_{|\kappa|\leq m} \|\partial_\kappa f \|_{F^{\alpha-m}_{p,q}} 
	\quad\text{and}\quad
	\sum_{|\kappa|= m} \|\partial_\kappa f \|_{\dot{F}^{\alpha-m}_{p,q}}
\end{equation*}
are equivalent norms on $F^\alpha_{p,q}$ and $\dot{F}^\alpha_{p,q}$ respectively.
\end{enumerate}
\end{theorem}

\begin{proof}
See \cite[Theorem 2.3.8.]{Tri83} and \cite[Theorem 5.2.3.1.(i)]{Tri83}.
Note that the equivalency for the standard derivatives is not stated explicitly for the homogeneous spaces, but heavily implied when for the proof of \cite[Theorem 5.2.3.1.(iii)]{Tri83} there is only a reference to the proof of the inhomogeneous version, which makes use of \cite[Theorem 2.3.8.]{Tri83}, which is the inhomogeneous version of the result in question.
\end{proof}

\subsection{Triebel-Lizorkin and Besov-Lipschitz space characterizations}\label{subsec:tt - trilibeli space char}

Since fractional derivatives of the Poisson kernel were the main motivation for \cite{BuiCan17}, it is no surprise that also fractional derivatives of the generalized Poisson kernel $p^s_1$ satisfy the conditions of Theorem \ref{th:tt BuiCandy main result 1} and \ref{th:tt BuiCandy main result 2}, yielding the following result.

\begin{theorem}[Triebel-Lizorkin and Besov-Lipschitz space characterizations]\label{th:tt triebel lizorkin besov char}
Let $\alpha\in\R$, $\beta>\max\lbrace \alpha, 0\rbrace$, $0<p,q\leq\infty$, $s\in(0,2)$ and $f\in\CicRn$.\\
For $\beta\in(0,1]$, $\alpha<\beta$ we have
\begin{equation}\label{eq:tt beli with fracLapl}
	\|f\|_{\dot{B}^\alpha_{p,q}} \approx 
	\left(\int_0^\infty	\left( \int_{\Rn} |t^{-\frac{1}{q}-\alpha+\beta}P^s_t\fL{\beta}f(x)|^p \,dx\right)^\frac{q}{p}\,dt	\right)^\frac{1}{q},
\end{equation}
and for $p\neq\infty$
\begin{equation}\label{eq:tt trili with fracLapl}
	\|f\|_{\dot{F}^\alpha_{p,q}} \approx 
	\left(\int_{\Rn}	\left( \int_0^\infty |t^{-\frac{1}{q}-\alpha+\beta}P^s_t\fL{\beta}f(x)|^q \,dt\right)^\frac{p}{q}\,dx	\right)^\frac{1}{p}.
\end{equation}
Regarding standard derivatives in $x$, we have the following for $\alpha<1$,
\begin{equation}\label{eq:tt beli with nablax}
	\|f\|_{\dot{B}^\alpha_{p,q}} \gtrsim
	\left(\int_0^\infty \left(\int_{\Rn} |t^{-\frac{1}{q}-\alpha+1}\nabla_x P^s_tf(x)|^p \,dx\right)^\frac{q}{p}\,dt\right)^\frac{1}{q},
\end{equation}
and for $p\neq\infty$
\begin{equation}\label{eq:tt trili with nablax}
	\|f\|_{\dot{F}^\alpha_{p,q}} \gtrsim
	\left(\int_{\Rn}	\left( \int_0^\infty |t^{-\frac{1}{q}-\alpha+1}\nabla_x P^s_tf(x)|^q \,dt\right)^\frac{p}{q}\,dx	\right)^\frac{1}{p}.
\end{equation}
Regarding derivatives in $t$, for $\alpha<s$ we have
\begin{equation}\label{eq:tt beli with partialt}
	\|f\|_{\dot{B}^\alpha_{p,q}} \approx 
	\left(\int_0^\infty \left( \int_{\Rn} |t^{-\frac{1}{q}-\alpha+1}\partial_t P^s_tf(x)|^p \,dx\right)^\frac{q}{p}\,dt	\right)^\frac{1}{q},
\end{equation}
and for $p\neq\infty$
\begin{equation}\label{eq:tt trili with partialt}
	\|f\|_{\dot{F}^\alpha_{p,q}} \approx 
	\left(\int_{\Rn}	\left( \int_0^\infty |t^{-\frac{1}{q}-\alpha+1}\partial_t P^s_tf(x)|^q \,dt\right)^\frac{p}{q}\,dx	\right)^\frac{1}{p}.
\end{equation}
\end{theorem}

\begin{proof}
The proof consists of four steps.
First, we recall the Fourier transform of the kernel $p^s_1$ and collect some related facts that are vital for verifying the conditions \textbf{(C1)}-\textbf{(C3)}.
Second, we show that the kernels $\fL{\beta}p^s_1$, $\beta\in(0,1]$, and $\partial_t p^s_1 = \partial_t k^s\big|_{t=1}$ satisfy the conditions of \ref{th:tt BuiCandy main result 1} and \ref{th:tt BuiCandy main result 2} for suitably chosen parameters while the kernels $\partial_{x_j}p^s_1$ only satisfy the conditions of \ref{th:tt BuiCandy main result 1}, hence we only have an estimate instead of an equivalency in (\ref{eq:tt beli with nablax}) and (\ref{eq:tt trili with nablax}).
Third, we argue why we can always choose $\rho=0$ in the context of Theorem \ref{th:tt BuiCandy main result 1}.
Fourth, we show that the expressions from Theorem \ref{th:tt BuiCandy main result 1} and Theorem \ref{th:tt BuiCandy main result 2} for the listed kernels actually match with the expressions of this theorem.

\underline{Step 1:} We recall (\ref{eq:Pst fourier transform of ps1}),
\begin{equation}\label{eq:tt trilibeli char G_s def}
	\Fourier(p^s_1)(\xi) = G_s(\xi)\coloneqq
	\const \int_0^\infty \lambda^\frac{s}{2}e^{-\lambda-\frac{|2\pi\xi|^2}{4\lambda}}\frac{d\lambda}{\lambda}.
\end{equation}
Obviously the function $G_s$ can be defined for all $s\in\R$ since the exponential term dominates the integrand near zero as well as towards infinity and we have $G_s>0$ on $\Rn\setminus\lbrace 0\rbrace$.
This exponential term will further allow us to use Lebesgue's dominated convergence theorem when dealing with different derivatives of $G_s$. 
From \cite[Proposition 7.6(c)]{Hao16} we obtain that
\begin{equation}\label{eq:tt trilibeli char G_s -> infty}
	G_s(\xi)\leq \const e^{-\pi|\xi|}, \quad \text{when } |\xi|\geq \frac{1}{\pi}
\end{equation}
for some constant $\const>0$.
On the other hand we have
\begin{equation}\label{eq:tt trilibeli char G_s -> 0}
	G_s(\xi) = \left\lbrace \begin{array}{llll}
		 \mathcal{O}(|\xi|^{s})							\quad&	\text{if } s<0, \\
		 \mathcal{O}(\operatorname{ln} \frac{2}{|\xi|})	\quad&	\text{if } s=0, \\
		 \mathcal{O}(1)									\quad&	\text{if } s>0
	\end{array}\right.
	\quad \text{as }\xi\rightarrow 0.
\end{equation}
While the result is only formulated for $s>-n$, the proofs presented in \cite{Hao16} are also valid for any $s\in\R$.
With Lebesgue's dominated convergence theorem we see that 
\begin{equation}\label{eq:tt trilibeli char partial G_s}
	\partial_{\xi_j} G_s(\xi) = 
	 \const \int_0^\infty \lambda^\frac{s}{2}e^{-\lambda-\frac{|2\pi\xi|^2}{4\lambda}} \left(\frac{2\pi^2}{\lambda} \xi_j\right) \frac{d\lambda}{\lambda}
	 = 2\pi^2\const \xi_j G_{s-2}(\xi).
\end{equation}
It is easy to show via induction that for any multiindex $\kappa$ we have
\begin{equation}\label{eq:tt trilibeli char partialkappa G_s}
	\partial_\kappa G_s(\xi) = \sum_{\gamma\leq\kappa} c_\gamma \xi^\gamma G_{s-(|\kappa| + |\gamma|)}(\xi)
\end{equation}
for some constants $c_\gamma\in\R$ depending on $\kappa$ where the sum is over all multiindices $\gamma$ that are smaller than or equal to $\kappa$ in each component.
In particular, $G_s\in\C^{\infty}(\Rn\setminus\lbrace 0\rbrace)$.

\underline{Step 2:} We start the second part by checking the integrability of the kernels.
For our first kernel $\fL{\beta}p^s_1$, due to Lemma \ref{lem:fL Linfty estimate} and Lemma \ref{lem:fL decay} we have
\begin{equation*}
	|\fL{\beta}p^s_1(x)|\leq \min\left\lbrace \|p^s_1\|_{W^{1,\infty}(\Rn)}, \frac{C}{|x|^{n+\beta}} \right\rbrace
\end{equation*} 
for some constant $C>0$. 
Confer (\ref{eq:Pst pst derivatives}) to confirm the boundedness as well as the necessary decay estimate for $p^s_1$.
Therefore, we obviously have $(1+|\cdot|)^l\fL{\beta}p^s_1\in L^1(\Rn)$ for any $0\leq l < \beta$.
The second kernel $\partial_t p^s_1$ is given by 
\begin{equation*}
	\partial_t p^s_1(x) = \frac{s}{(|x|^2+1)^\frac{n+s}{2}}-\frac{n+s}{(|x|^2+1)^{\frac{n+s}{2}+1}}.
\end{equation*}
Here, we have $(1+|\cdot|)^l\partial_tp^s_1\in L^1(\Rn)$ for any $0\leq l < s$.
For the kernels $\partial_{x_i}p^s_1$, again with (\ref{eq:Pst pst derivatives}), we obtain
\begin{equation*}
	|\partial_{x_i}p^s_1(x)|\leq \max\left\lbrace \|\partial_{x_i} p^s_1\|_{L^\infty(\Rn)}, \frac{C}{|x|^{n+1}} \right\rbrace
\end{equation*}
which implies $(1+|\cdot|)^l\partial_{x_i}p^s_1\in L^1(\Rn)$ for any $0\leq l < 1$.
These observations are the main reason for the specific restrictions on $\alpha$ since in the context of \ref{th:tt BuiCandy main result 1} we need the integrability for $l>\alpha-\frac{n}{p}$.
We therefore fix 
	$l_{(-\Delta)^{\beta/2}}\in(\max\lbrace 0,\alpha\rbrace,\beta)$ for $\alpha<\beta$ as well as 
	$l_{\partial_t}\in(\max\lbrace 0,\alpha\rbrace,s)$ for $\alpha<s$ and
	$l_{\nabla_x}\in(\max\lbrace 0,\alpha\rbrace,1)$ for $\alpha<1$.
In the context of \ref{th:tt BuiCandy main result 2} we can always choose $l=0$ because we are only interested in functions $f\in\CicRn$, which can be interpreted as distributions of growth $0$ since the convolution of two Schwartz functions is a Schwartz function itself, see \cite[Proposition 2.2.7]{Gra14}.

In order to check \textbf{(C1)}, \textbf{(C2)} and \textbf{(C3)}, we need the Fourier transform of the kernels.
With the definition of the fractional Laplacian and (\ref{eq:tt trilibeli char G_s def}) we have
\begin{equation*}
	\Fourier(\fL{\beta}p^s_t)(\xi) = (2\pi|\xi|)^\beta\Fourier(p^s_1)(\xi) = C|\xi|^\beta G_s(\xi)
\end{equation*}
for some constant $C>0$.
With Lemma \ref{lem:P1t switching Fourier and partial_t}, Proposition \ref{prop:Pst kernel fourier symbol} and dominated convergence we obtain
\begin{equation}\label{eq:tt partial_t pst Fourier}
	\Fourier(\partial_t p^s_1)(\xi)
	= \partial_t c_{n,s} \int_0^\infty \lambda^\frac{s}{2}e^{-\lambda-\frac{|t\xi|^2}{c\lambda}}\frac{d\lambda}{\lambda}
	= c_{n,s} \int_0^\infty \lambda^\frac{s}{2}e^{-\lambda-\frac{|t\xi|^2}{c\lambda}}\left(-\frac{2|\xi|^2}{c\lambda}t\right)\frac{d\lambda}{\lambda}.
\end{equation}
The conditions of the lemma are fulfilled due to (\ref{eq:Pst pst derivatives}).
Setting $t=1$,
\begin{equation*}
	\Fourier(\partial_t p^s_1)(\xi) = C|\xi|^2 G_{s-2}(\xi)
\end{equation*}
for some constant $C>0$.
For the last kernel, with the usual interaction between Fourier transform and derivatives,
\begin{equation*}
	\Fourier(\partial_{x_i}p^s_1)(\xi) = C \xi_i G_s(\xi)
\end{equation*}
for some $C\in\mathbb{C}\setminus\lbrace 0 \rbrace$.
Since $G_s>0$ on $\Rn\setminus\lbrace 0\rbrace$, it is clear that the kernels $\fL{\beta}p^s_1$ and $\partial_t p^s_1$ fulfill condition \textbf{(C2)}.

To conclude the second step of the proof, it remains to check the behavior of the kernels' Fourier transforms and their derivatives for $\xi\rightarrow 0$ as well as for $|\xi|\rightarrow \infty$.
With (\ref{eq:tt trilibeli char partialkappa G_s}) and since $\partial_\kappa|\xi|^\beta = \sum_{\gamma\leq\kappa} c_\gamma |\xi|^{\beta-2|\gamma|}\xi^\gamma$ for some constants $c_\gamma\in\R$, we obtain 
\begin{equation*}
	\partial_\kappa \Fourier(\fL{\beta}p^s_1)(\xi) 
	= \sum_{\eta\leq\kappa} c_\eta \partial_\eta |\xi|^\beta \partial_{\kappa-\eta} G_s(\xi)
	= \sum_{\eta\leq\kappa} \sum_{\gamma\leq\eta} \sum_{\tilde{\gamma}\leq\kappa-\eta}
				c_{\gamma,\tilde{\gamma}}\xi^{\gamma+\tilde{\gamma}}|\xi|^{\beta-2|\gamma|}G_{s-|\eta|-|\tilde{\gamma}|}(\xi)
\end{equation*}
with some constants $c_{\gamma,\tilde{\gamma}}\in\R$.
For the second kernel, via induction,
\begin{equation*}
	\partial_\kappa\Fourier(\partial_t p^s_1)(\xi) 
	= \sum_{|\gamma|\leq |\kappa|+2} c_\gamma \xi^\gamma G_{s-|\kappa|-|\gamma|}(\xi)
\end{equation*}
with some constants $c_\gamma\in\R$ with $c_\gamma= 0$ for $|\kappa| + |\gamma| <2$.
For the third kernel, also via induction,
\begin{equation*}
	\partial_\kappa\Fourier(\partial_{x_i}p^s_1) 
	= \sum_{\gamma\leq \kappa + e_i} c_\gamma \xi^\gamma G_{s-|\kappa|-|\gamma|+1}(\xi)
\end{equation*}
with some constants $c_\gamma\in\mathbb{C}$ and $e_i$ denoting the corresponding standard base vector.
In particular, the Fourier transforms of all three kernels are in $\C^\infty(\Rn\setminus\lbrace 0 \rbrace)$.
Therefore, the regularity condition in \textbf{(C1)} and \textbf{(C3)} is fulfilled.
Furthermore, since $G_\sigma$ decays exponentially fast as $|\xi|\rightarrow\infty$ for any $\sigma\in\R$, see (\ref{eq:tt trilibeli char G_s -> infty}), all three kernels fulfill \textbf{(C3)} for any $\Lambda\geq 0$, $m\in\R$.
It remains to check the behavior of the derivatives as $\xi\rightarrow 0$. 
For this we rely on (\ref{eq:tt trilibeli char G_s -> 0}) and obtain
\begin{equation*}
	\partial_\kappa \Fourier(\fL{\beta} p^s_1)(\xi)
	= \sum_{\eta\leq\kappa} \sum_{\gamma\leq\eta} \sum_{\tilde{\gamma}\leq\kappa-\eta}
				\mathcal{O}(|\xi|^{\beta-|\gamma|+|\tilde{\gamma}|}(1+|\xi|^{s-|\eta|-|\tilde{\gamma}|}) )
	= \mathcal{O}(|\xi|^{\beta-|\kappa|}) 
	\quad\text{as } \xi\rightarrow 0.
\end{equation*}
\begin{equation*}
	\partial_\kappa\Fourier(\partial_tp^s_1) (\xi)
	= \sum_{|\gamma|\leq |\kappa|+2} \mathcal{O}(|\xi|^\gamma|\xi|^{s-|\kappa|-|\gamma|})
	= \mathcal{O}(|\xi|^{s-|\kappa|})
	\quad \text{as } \xi\rightarrow 0.
\end{equation*}
For the kernels $\partial_{x_i}p^s_1$, $s\neq 1$, we have
\begin{equation*}
	\partial_\kappa\Fourier(\partial_{x_i}p^s_1) (\xi)
	= \left\lbrace \begin{array}{llll}
		 \mathcal{O}(|\xi|)											\quad&	\text{if } \kappa=0, \\
		 \sum_{\gamma\leq \kappa+e_i} \mathcal{O}(|\xi|^{|\gamma|}(1+|\xi|^{s-|\kappa|-|\gamma|+1}))
		 = \mathcal{O}(1+|\xi|^{1+s-|\kappa|})						\quad&	\text{if } |\kappa|\geq 1 
	\end{array}\right.
\end{equation*}
as $\xi\rightarrow 0$. 
For $s=1$ the only combination we have to treat differently is $|\kappa|+|\gamma|=2$ because then we have
\begin{equation*}
	c_\gamma \xi^\gamma G_{s-|\kappa|-|\gamma|+1}(\xi) 
	= \mathcal{O}\left(|c_\gamma| |\xi|^{|\gamma|} \operatorname{ln}\frac{2}{|\xi|}\right)
	\quad\text{as } \xi\rightarrow 0.
\end{equation*}
One can easily check that $c_\gamma=0$ for $|\kappa|= 2$ and $\gamma=0$. 
For $\gamma\neq 0$ we have $|\xi|^{|\gamma|} \operatorname{ln}\frac{2}{|\xi|}=\mathcal{O}(1)$ as $\xi\rightarrow 0$ and therefore end up with the same estimates as for $s\neq 1$.
Therefore, the kernel $\fL{\beta}p^s_1$ fulfills \textbf{(C1)} for all $\Lambda\geq 0$ and $r=\beta$, the second kernel $\partial_t p^s_1$ for all $\Lambda\geq 0$ and $r=s$. 
The kernels $\partial_{x_i}p^s_1$ fulfills \textbf{(C1)} for all $\Lambda\geq 0$ and $r=1$.

We conclude that all these kernels fulfill the conditions of Theorem \ref{th:tt BuiCandy main result 1} for $\alpha<\beta$ with $l=l_{(-\Delta)^{\beta/2}}$ as well as $\alpha<s$ with $l=l_{\partial_t}$ and $\alpha<1$ with $l=l_{\nabla_x}$ respectively.
Moreover, $\fL{\beta} p^s_1$ and $\partial_t p^s_1$ fulfill the conditions of Theorem \ref{th:tt BuiCandy main result 2} for functions $f\in\CicRn$ with $\alpha<\beta$ and $\alpha<s$ respectively.

\underline{Step 3:} 
Let $\rho$ be the polynomial from Theorem \ref{th:tt BuiCandy main result 1} for $f\in\CicRn$ and the respective kernel.
$f\in\CicRn$ is a distribution of growth $0$ since again the convolution with a Schwartz function is another Schwartz function. 
Therefore, if the distribution $f-\rho$ is of growth $l$, then $\rho$ has to be a distribution of growth $l$.
In particular, $\rho$ then is of degree at most $\lfloor l \rfloor$ due to Lemma \ref{lem:tt polynomials as distirbutions of growth l}.

For the first kernel, we have $l_{(-\Delta)^{\beta/2}}<\beta\leq 1$, implying that $\rho$ is of degree $0$ and therefore a constant, $\rho=c$. 
We have
\begin{equation*}
	(\fL{\beta}p^s_1 * \rho)(x) = c \int_\Rn \fL{\beta}p^s_1(x-y)\,dy 
	=  		c \int_\Rn \fL{\beta}p^s_1(y) e^{-2\pi i y\cdot 0}\,dy
	\equiv 	c\Fourier(\fL{\beta}p^s_1)(0).
\end{equation*}
Because of $\fL{\beta}p^s_1$ fulfilling the integrability condition as discussed at the beginning of the second step, we have $\fL{\beta}p^s_1\in L^1(\Rn)$ and therefore the continuity of the Fourier transform, see \cite[Exercise 2.2.6]{Gra14}.
Since $\fL{\beta}p^s_1$ further fulfills \textbf{(C1)} for $r=\beta>0$, 
\begin{equation*}
	\fL{\beta}p^s_1 * \rho = c\Fourier(\fL{\beta}p^s_1)(0) = \lim_{\xi\rightarrow 0} c\Fourier(\fL{\beta}p^s_1)(\xi) = 0.
\end{equation*}
For the dilations of $\fL{\beta}p^s_1$ we obtain $(\fL{\beta}p^s_1)_t*\rho = \fL{\beta}p^s_1 * (t^{-n}(\rho)_{t^{-1}})$ via substitution and therefore have $(\fL{\beta}p^s_1)_t*\rho=0$ with the same argument as above. 
Therefore, we can always choose $\rho=0$ in the context of Theorem \ref{th:tt BuiCandy main result 1} since it does not affect the value of the left hand sides. 

Since $l_{\nabla_x}<1$, this argument works analogously for the kernels $\partial_{x_i}p^s_1$.

Regarding the remaining kernel $\partial_t p^s_1$, we need to cover the two cases $s\leq 1$ and $1<s<2$.
The former is again covered analogously to the kernels above since $l_{\partial_t}<s$.
In the latter case, $l_{\partial_t}\geq 1$ is possible, meaning that $\rho$ might be of degree one. 
Since we already know that the convolution of the kernel with constants disappears, we only need to cover the monomials of degree $1$.
For $\rho = x_j$, $j=1,\ldots,n$ we have 
\begin{equation*}
	(\partial_t p^s_1 * \rho)(x) = \int_\Rn \partial_t p^s_1(y)(x_j-y_j)\,dy = -\int_\Rn \partial_t p^s_1(y)y_j\,dy 
	\equiv \frac{1}{2\pi i} \partial_{x_j}\Fourier(\partial_t p^s_1)(0)
\end{equation*}
where we used that $\Fourier(-2\pi i x_j g(x))(\xi) = \partial_{\xi_j}\Fourier(g)(\xi)$ for any $g\in L^1(\Rn)$ satisfying $|x|g\in L^1(\Rn)$, see the proof of \cite[Proposition 2.2.11]{Gra14}.
With $s>1$ we have $(1+|\cdot|)\partial_t p^s_1 \in L^1(\Rn)$. 
Moreover, because of this, the first order derivatives of the Fourier transform of $\partial_t p^s_1$ are continuous. 
With $\partial_t p^s_1$ fulfilling \textbf{(C1)} for $r=s>1$,
\begin{equation*}
	\partial_t p^s_1 * \rho = \frac{1}{2\pi i} \partial_{x_j}\Fourier(\partial_t p^s_1)(0) 
	= \lim_{\xi\rightarrow 0} \frac{1}{2\pi i} \partial_{x_j}\Fourier(\partial_t p^s_1)(\xi)
	= 0.
\end{equation*}
Therefore, just as for the other kernels we can always choose $\rho = 0$ without affecting the estimates.

\underline{Step 4:}
First, we investigate how the kernels behave under dilation.
Recall that $(p^s_1)_t = p^s_t$. 
Therefore, we obtain $\partial_{x_i} p^s_t = \partial_{x_i} (p^s_1)_t = t^{-1} (\partial_{x_i} p^s_1)_t$ as well as
\begin{mathex}
	\partial_t p^s_t(x) = \partial_t (p^s_1)_t(x)
	&=& -n t^{-n-1} p^s_1(t^{-1}x) - t^{-n} t^{-2} x \cdot \nabla_x p^s_1(t^{-1}x) \\
	&=& -n t^{-1} (p^s_1)_t(x) - t^{-1} (x\cdot \nabla_x p^s_1)_t(x)\\
	&=& t^{-1} (\partial_t p^s_1)_t.
\end{mathex}
For $\fL{\beta}p^s_1$, we use that $\Fourier((g)_t)(\xi) = \Fourier(g)(t\xi)$ for $g\in L^1(\Rn)$, see \cite[Proposition 2.2.11]{Gra14}.
That proposition is formulated for Schwartz functions but the result can be extended to $L^1$ with the same arguments.
We then have
\begin{mathex}
	\Fourier(\fL{\beta} p^s_t)(\xi) 
	&=& (2\pi|\xi|)^\beta \Fourier( (p^s_1)_t )(\xi) 
	&=& (2\pi|\xi|)^\beta \Fourier( (p^s_1) )(t \xi) \\
	&=& t^{-\beta} (2\pi|t\xi|)^\beta \Fourier( (p^s_1) )(t\xi) 
	&=& t^{-\beta} \Fourier((\fL{\beta} p^s_1)_t)(\xi), \\
\end{mathex}
that is $\fL{\beta} p^s_t = t^{-\beta} (\fL{\beta} p^s_1)_t$.

We now conclude this proof by showing the Besov-Lipschitz estimates utilizing the interaction of the differential operators with the $s$-harmonic extension, in particular some integration-by-parts properties of these operators.
For the kernel $\fL{\beta}p^s_1$, with Lemma \ref{lem:fL integration by parts}
\begin{mathex}
	((\fL{\beta}p^s_1)_t * f)(x)	&=& t^{\beta} \int_{\Rn} \fL{\beta} p^s_t(x-y) f(y)\,dy  \\
								&=& t^{\beta} \int_{\Rn}  p^s_t(x-y) \fL{\beta}f(y)\,dy 	\\
								&=& t^{\beta} P^s_t(\fL{\beta}f)(x).
\end{mathex}
Since we have verified all the conditions for Theorems \ref{th:tt BuiCandy main result 1} and \ref{th:tt BuiCandy main result 2} in Step 2 and dealt with the appearing polynomials in Step 3, we now obtain (\ref{eq:tt beli with fracLapl}),
\begin{mathex}
	\|f\|_{\dot{B}^\alpha_{p,q}}
	&\approx &	\bigg(\int_0^\infty (t^{-\alpha}\|(\fL{\beta}p^s_1)_t * f\|_{L^p})^q\frac{dt}{t}\bigg)^\frac{1}{q}\\
	&=&			\bigg(\int_0^\infty (t^{-\frac{1}{q}-\alpha+\beta}\|P^s_t\fL{\beta}f\|_{L^p})^q\,dt\bigg)^\frac{1}{q}\\
	&=&			\left(\int_0^\infty	\left( \int_{\Rn} |t^{-\frac{1}{q}-\alpha+\beta}P^s_t\fL{\beta}f(x)|^p \,dx\right)^\frac{q}{p}\,dt	\right)^\frac{1}{q}
\end{mathex}
for $\alpha<\beta$, $0<p,q\leq\infty$.
For the next kernel, $\partial_t p^s_1$, via Lemma \ref{lem:Pst derivatives of Pst 1}
\begin{equation*}
	((\partial_t p^s_1)_t * f)(x)	= t (\partial_t p^s_t*f)(x)
									= t \partial_t P^s_t f (x).
\end{equation*}
We then obtain (\ref{eq:tt beli with partialt}) analogously to (\ref{eq:tt beli with fracLapl}).
Regarding the kernels $\partial_{x_i} p^s_1$, also via Lemma \ref{lem:Pst derivatives of Pst 1},
\begin{equation*}
	((\partial_{x_i} p^s_1)_t * f)(x)= t (\partial_{x_i} p^s_t*f)(x)
									= t \partial_{x_i} P^s_t f (x).
\end{equation*}
Assembling the estimates for these kernels and using that $\|\cdot\|_{L^{\tilde{p}}(\R^k)}$ is a quasi-norm for all $\tilde{p}>0$ and $k\in\N$, we obtain (\ref{eq:tt beli with nablax}),
\begin{mathex}
	\|f\|_{\dot{B}^\alpha_{p,q}}
	&\gtrsim &	\sum_{i=1}^n\bigg(\int_0^\infty(t^{-\frac{1}{q}-\alpha+1}\|\partial_{x_i}P^s_t f\|_{L^p})^q\,dt\bigg)^\frac{1}{q}\\
	&\gtrsim &	\bigg(\int_0^\infty\bigg(\sum_{i=1}^n t^{-\frac{1}{q}-\alpha+1}\|\partial_{x_i}P^s_t f\|_{L^p}\bigg)^q\,dt\bigg)^\frac{1}{q}\\
	&\gtrsim &	\bigg(\int_0^\infty\bigg(t^{-\frac{1}{q}-\alpha+1}\bigg\|\sum_{i=1}^n |\partial_{x_i}P^s_t f|\bigg\|_{L^p}\bigg)^q\,dt\bigg)^\frac{1}{q}\\
	&\approx &	\sum_{i=1}^n\bigg(\int_0^\infty(t^{-\frac{1}{q}-\alpha+1}\|\nabla_x P^s_t f\|_{L^p})^q\,dt\bigg)^\frac{1}{q}.
\end{mathex}
The Triebel-Lizorkin estimates are obtained analogously for all these kernels. 
\end{proof}

\begin{remark}\label{re:tt annotations regarding the trilibeli characterizations}
\begin{enumerate}[label=(\roman*)]
\item	While we we have (\ref{eq:tt beli with nablax}) and (\ref{eq:tt trili with nablax}) only as estimates, these are listed as another equivalency in \cite[Theorem 10.8]{LeSchi20}.
		We did not get the equivalency because the kernels $\partial_{x_i} p^s_1$	individually do not fulfill the Tauberian condition \textbf{(C2)}.
		While $\nabla_x p^s_1$ is not really a kernel in the sense of the results from \cite{BuiCan17}, it does sort of fulfill a Tauberian condition.
		For every $\xi\in S^{n-1}$ there exists $i\in\lbrace 1,\ldots,n\rbrace$ (depending on $\xi$) such that for every $t>0$
		\begin{equation*}
			|\Fourier(\partial_{x_i} p^s_1)(t\xi)| = 2\pi |\xi_i| G_s(\xi)>0.
		\end{equation*}
		This makes it seem very likely that there indeed is an equivalency. 

\item	We only needed $f\in\CicRn$ for two steps. 
		In the context of Theorem \ref{th:tt BuiCandy main result 2} this ensured that $f\in\dot{F}^s_{p,q}$ and $\dot{B}^s_{p,q}$ while in the context of Theorem \ref{th:tt BuiCandy main result 1} we used that $f$ as a Schwartz function is a distribution of growth $0$ in order to argue that we may choose $\rho=0$.
		Naturally, the estimates against the Triebel-Lizorkin and Besov-Lipschitz quasi-norms of Theorem \ref{th:tt triebel lizorkin besov char} hold for all functions $f$ that are distributions of growth $0$.
		
\item 	Based on Theorem \ref{th:tt triebel lizorkin besov char}, we can easily obtain a variety of further space characterizations via Theorem \ref{th:tt lifting property of trlibeli spaces}, the lifting property of Triebel-Lizorkin and Besov-Lipschitz spaces. 
		For example, (\ref{eq:tt beli with fracLapl}) and (\ref{eq:tt trili with fracLapl}) hold for arbitrary $\beta>0$.
		For $\beta > 1$, $\alpha<\beta$ and $0<p,q\leq\infty$ as well as $f\in\CicRn$, with (ii) we have
		\begin{mathex}
			\|f\|_{\dot{B}^\alpha_{p,q}} 
			&\approx & \|\fL{\beta-1}f\|_{\dot{B}^{\alpha-(\beta-1)}_{p,q}}\\
			&\approx &	\left(\int_0^\infty	\left( \int_{\Rn} 
								\left|t^{-\frac{1}{q}-(\alpha-(\beta-1))+1}P^s_t\fL{1}\left(\fL{\beta-1}f(x)\right)\right|^p 
						\,dx\right)^\frac{q}{p}\,dt	\right)^\frac{1}{q}\\
			&=&			\left(\int_0^\infty	\left( \int_{\Rn} 
								|t^{-\frac{1}{q}-\alpha+\beta}P^s_t\fL{\beta}f(x)|^p 
						\,dx\right)^\frac{q}{p}\,dt	\right)^\frac{1}{q}.
		\end{mathex}
		because $\fL{\beta-1}f$ is a distribution of growth $0$ due to Lemma \ref{lem:fL decay} and Lemma \ref{lem:fL Linfty estimate}.
		The respective extension for $\dot{F}^\alpha_{p,q}$ is obtained analogously.
		
\item 	The second example for applying the lifting property are double derivatives. 
		Based on (\ref{eq:tt beli with partialt}) and (\ref{eq:tt beli with partialt}), for $\alpha<s+1$, $0<p,q\leq\infty$ and $f\in\CicRn$ we obtain
		\begin{mathex}
			\|f\|_{\dot{B}^\alpha_{p,q}} 
			&\approx &	\sum_{j=1}^n\|\partial_{x_j}f\|_{\dot{B}^{\alpha-1}_{p,q}}\\
			&\approx & 	\sum_{j=1}^n\left(\int_0^\infty \bigg( \int_{\Rn} 
										|t^{-\frac{1}{q}-(\alpha-1)+1}\partial_t P^s_t (\partial_{x_i}f)(x)|^p 
						\,dx\bigg)^\frac{q}{p}\,dt	\right)^\frac{1}{q}.
		\end{mathex}
		 Redistributing the additional derivative with Lemma \ref{lem:Pst derivatives of Pst 2}, using the quasi-norm property of $\|\cdot\|_{L^{\tilde{p}}(\R^k)}$ for one direction of the estimate and by estimating the individual summands for the other direction, 
		\begin{mathex}
			\|f\|_{\dot{B}^\alpha_{p,q}}
			&\approx &	\left(\int_0^\infty \bigg( \int_{\Rn} 
										\sum_{j=1}^n|t^{-\frac{1}{q}-(\alpha-1)+1}\partial_t\partial_{x_i} P^s_t f(x)|^p 
						\,dx\bigg)^\frac{q}{p}\,dt	\right)^\frac{1}{q}\\
			&\approx &	\left(\int_0^\infty \left( \int_{\Rn} 
										|t^{-\frac{1}{q}-\alpha+2}\partial_t\nabla_x P^s_tf(x)|^p 
						\,dx\right)^\frac{q}{p}\,dt	\right)^\frac{1}{q}.
		\end{mathex}
		Analogously, for $\alpha<s+1$, $0<p<\infty$, $0<q\leq\infty$ and $f\in\CicRn$ we obtain
		\begin{equation*}
			\|f\|_{\dot{F}^\alpha_{p,q}} \approx 
			\left(\int_{\Rn}	\left( \int_0^\infty |t^{-\frac{1}{q}-\alpha+1}\partial_t\nabla_x P^s_tf(x)|^q \,dt\right)^\frac{p}{q}\,dx	\right)^\frac{1}{p}.
		\end{equation*}
		Apparently, we could reiterate this process, adding arbitrary many derivatives.
		This allows the characterization of $\dot{B}^\alpha_{p,q}$ and $\dot{F}^\alpha_{p,q}$ for arbitrary $\alpha\in\R$ as long as  enough derivatives are added.
		Of course, there are analogous versions for (\ref{eq:tt beli with nablax}) and (\ref{eq:tt trili with nablax}).
\end{enumerate}
\end{remark}

Due to the Tribel-Lizorkin characterization we now have some tools to deal with certain terms involving derivatives of the $s$-harmonic extension such as $\nabla_\Rnp P^s_t f$, $\nabla_x\nabla_\Rnp P^s_t f$ or $P^s_t(\fL{\beta}f)$.
We are still missing tools to deal with terms where the $s$-harmonic extension occurs unmodified.
Naturally, we wonder whether analogous results can be obtained for the unmodified kernel $p^s_1$.
Following the scheme of the proof for Theorem \ref{th:tt triebel lizorkin besov char}, one can easily see that the $p^s_1$ satisfies the conditions of Theorem \ref{th:tt BuiCandy main result 1} and Theorem \ref{th:tt BuiCandy main result 2} for $\alpha<0$ with $l=0$. 
But since $\int_{\Rn} p^s_1 \,dx \neq 0$ in contrast to the other kernels, we are not able to ignore the polynomial $\rho$ appearing in Theorem \ref{th:tt BuiCandy main result 1}.

But with a simple trick we are able to circumvent these limitations.
Next to the estimates for the derivatives of the $s$-harmonic extension, we therefore also obtain an analogous estimate in terms of the unmodified $s$-harmonic extension.

\begin{corollary}\label{cor:tt trilibeli spaces via unmodified Fs}
Let $0<p, q \leq \infty$ and $s\in(0,2)$.
Denote with $F^s(x,t) = P^s_tf(x)$ the $s$-harmonic extension of $f\in\CicRn$.
Then for $\beta\in(0,n)$
\begin{equation}\label{eq:tt beli via unmodified Fs}
	\|f\|_{\dot{B}^{-\beta}_{p,q}} \approx 
	\left(\int_0^\infty	\left( \int_{\Rn} |t^{-\frac{1}{q}+\beta}F^s(x,t)|^p \,dx\right)^\frac{q}{p}\,dt	\right)^\frac{1}{q},
\end{equation}
and for $p\neq\infty$
\begin{equation}\label{eq:tt trili via unmodified Fs}
	\|f\|_{\dot{F}^{-\beta}_{p,q}} \approx 
	\left(\int_{\Rn}	\left( \int_0^\infty |t^{-\frac{1}{q}+\beta}F^s(x,t)|^q \,dt\right)^\frac{p}{q}\,dx	\right)^\frac{1}{p}.
\end{equation}
\end{corollary}

\begin{proof}
We only show (\ref{eq:tt beli via unmodified Fs}), (\ref{eq:tt trili via unmodified Fs}) follows analogously.
$I^\beta f$ is a distribution of growth $0$ since $I^\beta f\in L^s( \Rn)$ for some $1<s<\infty$ due to Theorem \ref{th:Rp convergence and Lq boundedness}.
Therefore, with (\ref{eq:tt beli with fracLapl}), Remark \ref{re:tt annotations regarding the trilibeli characterizations}(ii),(iii) and the lifting property for homogeneous Besov-Lipschitz spaces we have
\begin{equation*}
	\|f\|_{\dot{B}^{-\beta}_{p,q}} \approx \|I^\beta f\|_{\dot{B}^{0}_{p,q}}\approx
	\left(\int_0^\infty	\left( \int_{\Rn} |t^{-\frac{1}{q}-0+\beta}P^s_t(\fL{\beta} I^\beta f)(x)|^p \,dx\right)^\frac{q}{p}\,dt	\right)^\frac{1}{q}.
\end{equation*}
The result follows since $\fL{\beta}$ and $I^\beta$ are inverse.
\end{proof}


\subsection{Lorentz spaces and an interpolation theorem}\label{subsec:tt - lorentz spaces}
Besides the Triebel-Lizorkin and Besov-Lipschitz spaces, there is a second type of spaces which we need to introduce in more detail,
the Lorentz spaces $\Lor{p}{q}{\Rn}$, $0<p<\infty, 0<q\leq\infty$.
They are a generalization of the Lebesgue spaces, offering a finer scale, which is not covered by the Triebel-Lizorkin and Besov-Lipschitz spaces. 
The Lebesgue spaces are included in the Lorentz scala with $\Lor{p}{p}{\Rn}=L^p(\Rn)$.
Despite not defining Lorentz spaces for $p=\infty$, we therefore still use the convention $\Lor{\infty}{\infty}{\Rn} = L^\infty(\Rn)$.
These Lorentz-spaces are defined as
\begin{equation*}
	\Lor{p}{q}{\Rn} \coloneqq 
	\lbrace f\colon \Rn\rightarrow\R \, ;\, f \text{ is measurable and } \|f\|_{\Lor{p}{q}{\Rn}}<\infty\rbrace.
\end{equation*}
For measurable functions $f\colon \Rn\rightarrow\R$ the Lorentz space (quasi-)norm $\|\cdot\|_\Lor{p}{q}{\Rn}$ is  
\begin{equation*}
	\|f\|_{\Lor{p}{q}{\Rn}}\coloneqq 
	\left\lbrace\begin{array}{ll}
		\left(\int_0^\infty\left(t^\frac{1}{p}f^*(t)\right)^q\frac{dt}{t}\right)^\frac{1}{q} \quad &\text{if } p,q\in(0,\infty),\\
		\sup_{t>0}t^\frac{1}{p}f^*(t)														\quad & \text{if }  q=\infty
	\end{array}\right.
\end{equation*}
where $f^*$ is the decreasing rearrangement of $f$ given by
\begin{equation*}
	f^\ast(t) \coloneqq \inf\lbrace s>0 \, ;\, \mathcal{L}^n\lbrace |f|>s\rbrace)\leq t \rbrace
\end{equation*}
with $\mathcal{L}^n$ denoting the $n$-dimensional Lebesgue-measure.

One key property to us will be that similar to Lebesgue spaces there is a Hölder inequality for the Lorentz spaces.

\begin{theorem}[Hölder's inequality in Lorentz spaces, \cite{ONe63}]\label{th:tt lorentz holder}
Let $ 0< p_1, p_2, p<\infty $ and $0<q_1,q_2,q\leq\infty$ with
\begin{equation*}
	\frac{1}{p}=\frac{1}{p_1}+\frac{1}{p_2} \quad\text{and}\quad \frac{1}{q}=\frac{1}{q_1}+\frac{1}{q_2}.
\end{equation*}
Then for measurable functions $f,g\colon \Rn\rightarrow\R$ we have
\begin{equation}
	\|fg\|_{\Lor{p}{q}{\Rn}} \lesssim \|f\|_{\Lor{p_1}{q_1}{\Rn}} \|g\|_{\Lor{p_2}{q_2}{\Rn}}.
\end{equation}
In particular, for $1=\frac{1}{p_1}+\frac{1}{p_2}$ and $1=\frac{1}{q_1}+\frac{1}{q_2}$,
\begin{equation}
	\int_\Rn |fg| \lesssim \|f\|_{\Lor{p_1}{q_1}{\Rn}} \|g\|_{\Lor{p_2}{q_2}{\Rn}}.
\end{equation}
\end{theorem}

The other key property to us is the weak type interpolation theorem for Lorentz spaces as given in \cite[Section 3]{Hun66}.
We use this theorem to obtain Lorentz estimates from $L^p$-estimates, adding another scale to our estimates without much effort.

\begin{theorem}[Weak type interpolation theorem for Lorentz spaces, \cite{Hun66}]\label{th:tt lorentz interpolation}
Let $T\colon \Lor{p_i}{q_i}{\Rn}\rightarrow \Lor{\tilde{p}_i}{\tilde{q}_i}{\Rn}$ with $0<p_i,\tilde{p}_i<\infty$, $0<q_i,\tilde{q}_i\leq\infty$ for $i=0,1$ and $p_0<p_1$, $\tilde{p}_0\neq \tilde{p}_1$.
Assume that $T$ is quasi-linear, that is that $T(f+g)$ is defined whenever $Tf$ and $Tg$ are defined and 
	$|T(f+g)|\lesssim |Tf|+|Tg|$ 
almost everywhere.
Further assume that  
\begin{equation*}
	\|Tf\|_{\Lor{\tilde{p}_i}{\tilde{q}_i}{\Rn}} \lesssim \|f\|_{\Lor{p_i}{q_i}{\Rn}}, \quad i=0,1.
\end{equation*}
Then
\begin{equation*}
	\|Tf\|_{\Lor{\tilde{p}_\theta}{s}{\Rn}} \lesssim \|f\|_{\Lor{p_\theta}{q}{\Rn}},
\end{equation*}
where $q\leq s$ and, for $0<\theta<1$, 
	$\frac{1}{p_\theta}=\frac{1-\theta}{p_0}+\frac{\theta}{p_1}$ as well as 
	$\frac{1}{\tilde{p}_\theta}=\frac{1-\theta}{\tilde{p}_0}+\frac{\theta}{\tilde{p}_1}$.

The theorem also holds for $(p_1,q_1)=(\infty,\infty)$ or $(\tilde{p}_i,\tilde{q}_i)=(\infty,\infty)$.
\end{theorem}

\begin{remark}
Note that the image of $(0,1)$ under $\theta \mapsto p_\theta$ is $(p_0,p_1)$. 
Therefore, the interpolation theorem captures all $p$ between $p_0$ and $p_1$ as well as all $\tilde{p}$ between $\tilde{p}_0$ and $\tilde{p_1}$ but only in the specific relation dictated by the respective $\theta$.
\end{remark}

%
%

\subsection{Building blocks: $BMO$, fractional Sobolev \& H\"older space estimates}\label{subsec:tt - building blocks}

At the start of Subsection \ref{subsec:tt - triebel lizorkin besov}, we listed identifications of Triebel-Lizorkin and Besov-Lipschitz spaces with a variety of well known and important function spaces.
Thanks to these, based on Theorem \ref{th:tt triebel lizorkin besov char} we now easily obtain a multitude of more specific space characterizations via the generalized Poisson operator.

The first couple of more specific space characterizations are related to the two common variants of fractional Sobolev spaces.
On the one hand, there are the Slobodeckij spaces $W^s_p(\Rn)$.
For $1\leq p <\infty$ and $0<s\notin\N$ with $s=[s] +\lbrace s \rbrace$, where  $[s]$ is the integer part of $s$, they are defined as 
\begin{equation}\label{eq:tt def Slobodeckij}
	W^s_p(\Rn) \coloneqq \left\lbrace f\in W^{[s]}_p(\Rn) \, ;\, 
						\|f\|_{W^s_p(\Rn)} \coloneqq	\|f\|_{W^{[s]}_p(\Rn)} 
													+ \sum_{|\kappa|=[s]} [\partial_{\kappa}f]_{W^{\lbrace s\rbrace,p}(\Rn)}<\infty
				\right\rbrace,
\end{equation}
where the seminorm $[\cdot]_{W^{\nu,p}(\Rn)}$, $\nu\in(0,1)$, is given by 
\begin{equation*}
	[f]_{W^{\nu,p}(\Rn)}\coloneqq \left( \int_\Rn \int_\Rn \frac{|f(x)-f(y)|^p}{|x-y|^{n+\nu p}}\,dx\,dy\right)^\frac{1}{p}.
\end{equation*}
Recalling the list from Subsection \ref{subsec:tt - triebel lizorkin besov}, the Slobodeckij spaces can be identified with $F^s_{p,p}$.
On the other hand, there are the Bessel-potential spaces $H^s_p(\Rn)$ defined as
\begin{equation*}
	H^s_p(\Rn) \coloneqq \left\lbrace f\in\mathcal{S}'(\Rn) \, ;\, 
					\|f\|_{H^s_p(\Rn)} 
					\coloneqq \left\|\Fourier^{-1}\left((1+|\xi|^2)^\frac{s}{2}\Fourier(f)\right)\right\|_{L^p(\Rn)}<\infty
				 \right\rbrace
\end{equation*}
for $s\in\R$ and $1<p<\infty$. These spaces can be identified with $F^s_{p,2}(\Rn)$.

Note that for $m\in\Nz$, $p\neq 2$ we have $W^{m,p}(\Rn) =  H^m_p(\Rn) = F^m_{p,2}(\Rn) \neq F^m_{p,p}(\Rn)$.
Therefore, whenever we obtain results for the Slobodeckij spaces via the Triebel-Lizorkin identification, we need to somewhat counter-intuitively exclude the integer case, i.e. the classical Sobolev spaces. 

Further note that both of these spaces are inhomogeneous while we only have the space characterization for homogeneous Triebel-Lizorkin spaces. 
The corresponding $\dot{F}^s_{p,p}$-norm or $\dot{F}^s_{p,2}$-norm, which we obtain via the Poisson kernels in the following two propositions, might intuitively be described as the highest order seminorm contributing to the norm of the corresponding inhomogeneous space. 

\begin{proposition}[Slobodeckij space characterizations]\label{prop:tt slobodeckij space char}
Let  $1<p<\infty$ and $s\in(0,2)$.
Denote with $F^s(x,t)=P^s_tf(x)$ the $s$-harmonic extension of $f\in\CicRn$.\\
Regarding the $x$-derivatives, for $\nu\in(0,1)$ we have
\begin{equation}\label{eq:tt Slobodeckij via nabla_xFs}
	\left(\int_\Rn \int_0^\infty |t^{1-\frac{1}{p}-\nu}\nabla_x F^s(x,t)|^p \,dt\,dx \right)^\frac{1}{p}
	\lesssim [f]_{W^{\nu,p}(\Rn)}.
\end{equation}
Regarding the $t$-derivative, for $\nu\in(0,s)$, $\nu\neq 1$ we have
\begin{equation}\label{eq:tt Slobodeckij via partial_tFs}
	\left(\int_\Rn \int_0^\infty |t^{1-\frac{1}{p}-\nu}\partial_t F^s(x,t)|^p \,dt\,dx \right)^\frac{1}{p}
	\approx [f]_{W^{\nu,p}(\Rn)}.
\end{equation}
Regarding double derivatives, for $\nu\in(0,2)$, $\nu\neq 1$,
\begin{equation}\label{eq:tt Slobodeckij via nabla_x^2Fs}
	\left(\int_\Rn \int_0^\infty |t^{2-\frac{1}{p}-\nu}\nabla_x^2 F^s(x,t)|^p \,dt\,dx \right)^\frac{1}{p}
	\lesssim [f]_{W^{\nu,p}(\Rn)},
\end{equation}
and for $\nu\in(0,1+s)$, $\nu\notin\N$,
\begin{equation}\label{eq:tt Slobodeckij via partial_tnabla_xFs}
	\left(\int_\Rn \int_0^\infty |t^{2-\frac{1}{p}-\nu}\partial_t\nabla_x F^s(x,t)|^p \,dt\,dx \right)^\frac{1}{p}
	\lesssim [f]_{W^{\nu,p}(\Rn)}.
\end{equation}
\end{proposition}

\begin{proof}[Estimates (\ref{eq:tt Slobodeckij via nabla_xFs})-(\ref{eq:tt Slobodeckij via partial_tnabla_xFs})]
According to \cite[Theorem 2.3.3(ii)]{Tri92} we have
\begin{equation*}
	\|f\|_{F^\nu_{p,p}(\Rn)} \approx \|f\|_{L^p(\Rn)} + \|f\|_{\dot{F}^\nu_{p,p}(\Rn)}
\end{equation*}
Therefore, since $F^\nu_{p,p}(\Rn) = W^{\nu,p}(\Rn)$ for $\nu\notin\Nz$,
\begin{equation}\label{eq:tt frac sobolev char homogeneous trili vs slobodeckij}
	\|f\|_{L^p(\Rn)} + \|f\|_{\dot{F}^\nu_{p,p}(\Rn)} \approx \|f\|_{L^p(\Rn)} + [f]_{W^{\nu,p}(\Rn)}.
\end{equation}
All of the occurring (semi-)norms are homogeneous, but for different orders.
For all $\lambda >0$ we have 
\begin{mathex}
	\|f(\lambda\,\cdot)\|_{L^p(\Rn)} &=& \lambda^{-\frac{n}{p}}\|f\|_{L^p(\Rn)},\\
	{[}f(\lambda\,\cdot){]}_{W^{\nu,p}(\Rn)}&=& \lambda^{\nu-\frac{n}{p}} [f]_{W^{\nu,p}(\Rn)},\\
	\|f(\lambda\,\cdot) \|_{\dot{F}^\nu_{p,p}(\Rn)} &\approx& \lambda^{\nu-\frac{n}{p}}\|f\|_{\dot{F}^\nu_{p,p}(\Rn)},
\end{mathex}
where the first equations are obtained easily via substitution and the last estimate is due to \cite[Remark 5.1.3.4]{Tri83} with the constants independent of $\lambda$.
These observations in combination with (\ref{eq:tt frac sobolev char homogeneous trili vs slobodeckij}) yield
\begin{equation*}
	\lambda^{-\frac{n}{p}}\|f\|_{L^p(\Rn)} + \lambda^{\nu-\frac{n}{p}}\|f\|_{\dot{F}^\nu_{p,p}(\Rn)} 
	\leq C\left(\lambda^{-\frac{n}{p}}\|f\|_{L^p(\Rn)} + \lambda^{\nu-\frac{n}{p}}[f]_{W^{\nu,p}(\Rn)}\right)
\end{equation*}
for some constant $C>0$ independent of $\lambda$ and $f$. 
Thus,
\begin{equation*}
	\|f\|_{\dot{F}^\nu_{p,p}(\Rn)} 
	\leq \lambda^{-\nu}(C-1)\|f\|_{L^p(\Rn)} + C[f]_{W^{\nu,p}(\Rn)}.
\end{equation*}
Since $\nu>0$, for $\lambda \rightarrow \infty$ we obtain $\|f\|_{\dot{F}^\nu_{p,p}(\Rn)}\lesssim [f]_{W^{\nu,p}(\Rn)}$. 
Analogously, we get the estimate $[f]_{W^{\nu,p}(\Rn)}\lesssim \|f\|_{\dot{F}^\nu_{p,p}(\Rn)}$ and with that $[f]_{W^{\nu,p}(\Rn)}\approx \|f\|_{\dot{F}^\nu_{p,p}(\Rn)}$.

Therefore, we obtain (\ref{eq:tt Slobodeckij via nabla_xFs}) as a special case of (\ref{eq:tt trili with nablax}), (\ref{eq:tt Slobodeckij via partial_tFs}) from (\ref{eq:tt trili with partialt}) and (\ref{eq:tt Slobodeckij via nabla_x^2Fs}) as well as (\ref{eq:tt Slobodeckij via partial_tnabla_xFs}) with Remark \ref{re:tt annotations regarding the trilibeli characterizations}(iv) applied to (\ref{eq:tt trili with nablax}) and (\ref{eq:tt trili with partialt}) respectively.
\end{proof}

\begin{proposition}[Bessel-potential space characterizations]\label{prop:tt bessel potential spaces char}
Let  $1<p<\infty$, $0< q \leq \infty$ and $s\in(0,2)$.
Denote with $F^s(x,t)=P^s_tf(x)$ the $s$-harmonic extension of $f\in\CicRn$.\\
Regarding the $x$-derivatives, for $\nu\in[0,1)$ we have
\begin{equation}\label{eq:tt Bessel potential spaces via  nabla_xFs}
	\left\| x\mapsto \left(\int_0^\infty |t^{\frac{1}{2}-\nu}\nabla_x F^s(x,t)|^2\,dt\right)^\frac{1}{2}\right\|_{\Lor{p}{q}{\Rn}}
	\lesssim	\|\fL{\nu}f\|_{\Lor{p}{q}{\Rn}},
\end{equation}
Regarding the $t$-derivative, for $\nu\in(0,s)$ we have
\begin{equation}\label{eq:tt Bessel potential spaces via partial_tFs}
	\left\| x\mapsto \left(\int_0^\infty |t^{\frac{1}{2}-\nu}\partial_t F^s(x,t)|^2\,dt\right)^\frac{1}{2}\right\|_{\Lor{p}{q}{\Rn}}
	\lesssim	\|\fL{\nu}f\|_{\Lor{p}{q}{\Rn}}.
\end{equation}
Regarding double derivatives, for $\nu\in[0,2)$ we have
\begin{equation}\label{eq:tt Bessel potential spaces via nabla_x^2Fs}
	\left\| x\mapsto \left(\int_0^\infty |t^{\frac{3}{2}-\nu}\nabla_x^2 F^s(x,t)|^2\,dt\right)^\frac{1}{2}\right\|_{\Lor{p}{q}{\Rn}}
	\lesssim	\|\fL{\nu}f\|_{\Lor{p}{q}{\Rn}},
\end{equation}
and for $\nu\in[0,1+s)$,
\begin{equation}\label{eq:tt Bessel potential spaces via partial_tnabla_xFs}
	\left\| x\mapsto \left(\int_0^\infty |t^{\frac{3}{2}-\nu}\partial_t\nabla_x F^s(x,t)|^2\,dt\right)^\frac{1}{2}\right\|_{\Lor{p}{q}{\Rn}}
	\lesssim	\|\fL{\nu}f\|_{\Lor{p}{q}{\Rn}}.
\end{equation}
All these estimates also hold for $-n<\nu<0$ with $\fL{\nu}f$ replaced by the Riesz potential $I^{-\nu}f$. 
In fact, for $0<\nu<n$ we also have
\begin{equation}\label{eq:tt Bessel potential spaces via unmodified Fs}
	\left\| x\mapsto \left(\int_0^\infty |t^{-\frac{1}{2}+\nu}F^s(x,t)|^2 dt \right)^\frac{1}{2}   \right\|_{\Lor{p}{q}{\Rn}}
	\lesssim \|I^{\nu}f\|_{\Lor{p}{q}{\Rn}}.
\end{equation}
\end{proposition}

\begin{proof}[Estimates (\ref{eq:tt Bessel potential spaces via  nabla_xFs})-(\ref{eq:tt Bessel potential spaces via unmodified Fs})]
We first show the estimates for $p=q$, i.e. the $L^p$-norm.
Then, we obtain (\ref{eq:tt Bessel potential spaces via  nabla_xFs})-(\ref{eq:tt Bessel potential spaces via unmodified Fs}) via interpolation.

\underline{Step 1:}
Due to Theorem \ref{th:tt lifting property of trlibeli spaces} and since $\dot{F}^0_{p,2}(\Rn)=L^p(\Rn)$ for $\nu\geq 0$, we have
\begin{equation*}
	\|f\|_{\dot{F}^\nu_{p,2}(\Rn)} \approx \|\fL{\nu}f\|_{L^p(\Rn)}
\end{equation*}
while for $-n<\nu<0$
\begin{equation*}
	\|f\|_{\dot{F}^\nu_{p,2}(\Rn)} \approx \|I^{-\nu}f\|_{L^p(\Rn)}.
\end{equation*}
Therefore, we obtain (\ref{eq:tt Bessel potential spaces via  nabla_xFs}) and (\ref{eq:tt Bessel potential spaces via partial_tFs}) for $q=p$ from (\ref{eq:tt trili with nablax}) and (\ref{eq:tt trili with partialt}) respectively, (\ref{eq:tt Bessel potential spaces via unmodified Fs}) from Corollary \ref{cor:tt trilibeli spaces via unmodified Fs}. The $L^p$-versions of the other two estimates can be obtained with Remark \ref{re:tt annotations regarding the trilibeli characterizations}(iv).

\underline{Step 2:}
Based on these $L^p$-estimates we now obtain the $L^{p,q}$-estiamtes via the weak type interpolation theorem for Lorentz spaces, Theorem \ref{th:tt lorentz interpolation}, see \cite[Section 3]{Hun66}.
The details are only provided for (\ref{eq:tt Bessel potential spaces via  nabla_xFs}), the same arguments work for the other equations.
We want to apply the interpolation theorem to the operator 
\begin{equation*}
	g\rightarrow Tg\coloneqq\left(x\mapsto \left(\int_0^\infty |t^{\frac{1}{2}-\nu}\nabla_x P^s_t \dot{J}^{-\nu} g(x)|^2\,dt\right)^\frac{1}{2} \right)
\end{equation*}
that is defined for $g\in \dot{J}^{\nu}(\CicRn)$.
Here, $\dot{J}^\nu$ is the lifting operator from Theorem \ref{th:tt lifting property of trlibeli spaces} that up to constants equals the fractional Laplacians for $\nu\geq 0$ and the Riesz potentials for $-n<\nu<0$.
For any $1<p_0<p_1<\infty$ we have
\begin{equation*}
	\|Tg\|_{\Lor{p_i}{p_i}{\Rn}} = \|Tg\|_{L^{p_i}(\Rn)}\lesssim \|\dot{J}^\nu \dot{J}^{-\nu} g\|_{L{p_i}(\Rn)} = \| g\|_{\Lor{p_i}{p_i}{\Rn}}
\end{equation*}
due to the $L^p$-estimates from above. 
Furthermore, we have
\begin{equation*}
	|T(g_1+g_2)|\leq |Tg_1|+|Tg_2|
\end{equation*}
due to the Minkowski inequality.
Usually, in order to apply Theorem \ref{th:tt lorentz interpolation} the operator $T$ has to be defined on all of $L^{p_0}(\Rn)+L^{p_1}(\Rn)$.
We are only interested in interpolation estimates for functions $g\in \dot{J}^{\nu}(\CicRn)\subset L^{p_0}(\Rn)\cap L^{p_1}(\Rn)$ though.
Investigating the proof of Theorem \ref{th:tt lorentz interpolation} in \cite[Section 3]{Hun66}, for this weaker result it suffices to have $T$ defined on $L^{p_0}(\Rn)\cup L^{p_1}(\Rn)$.

Therefore, let us define $T$ on $L^{p_i}(\Rn)$ for $i=0,1$.
This would be very easy, if all elements in $\dot{J}^{-\nu}(L^{p_i}(\Rn))$ were distributions of growth $0$ since then $T$ would already be defined on $L^{p_i}(\Rn)$ without any further modifications due to Remark \ref{re:tt annotations regarding the trilibeli characterizations}. 
Sadly, this is not always the case.

According to \cite[Theorem 5.1.5]{Tri83}, $\mathcal{S}(\Rn)$ and therefore $\CicRn$ are dense in $\dot{F}^{-\nu}_{p_i,2}$.
Thus, due to Theorem \ref{th:tt lifting property of trlibeli spaces} $\dot{J}^{\nu}(\CicRn)$ is dense in $L^{p_i}(\Rn) = \dot{F}^0_{p_i,2}(\Rn)$.
Let $g\in L^{p_i}(\Rn)$.
Then there exists an approximating sequence $(g_n)_{n\in\N}\subset \dot{J}^{\nu}(\CicRn)$ with $g_n\overset{n\rightarrow\infty}{\longrightarrow} g$ in $L^{p_i}(\Rn)$.
We have 
\begin{equation*}
	Tg_n-Tg_m \leq Tg_m + T(g_n-g_m) - Tg_m = T(g_n-g_m) = T(g_m-g_n)
\end{equation*}
and therefore via symmetry $|Tg_n-Tg_m|\leq T(g_n-g_m)$.
Then $Tg_n$ is a Cauchy sequence in $\Lor{p_i}{p_i}{\Rn}$ since
\begin{equation*}
		\|Tg_n - Tg_m\|_{L^{p_i}(\Rn)}
		\leq	\|T(g_n-g_m)\|_{L^{p_i}(\Rn)}
		\lesssim	\|g_n-g_m\|_{L^{p_i}(\Rn)}.
\end{equation*}
Because the Lebesgue spaces are complete we can then set $Tg\coloneqq \lim_{n\rightarrow\infty} Tg_n$.
This definition is independent of the chosen sequence as we obtain $\|T\tilde{g}_n-Tg_n\|_{L^{p_i}(\Rn)} \lesssim \|\tilde{g}_n-g_n\|_{L^{p_i}(\Rn)}$ with the same argument as above for any other sequence $(\tilde{g}_n)_{n\in\N}\dot{J}^{\nu}(\CicRn)$.
Moreover, this choice is also well-defined when considering $g\in L^{p_0}(\Rn)\cap L^{p_1}(\Rn)$;
$\dot{J}^{-\nu}g$ can be approximated with the same sequence of compactly supported smooth functions both in $\dot{F}^{-\nu}_{p_0,2}$ and $\dot{F}^{-\nu}_{p_1,2}$, see the proof of \cite[Theorem 2.3.3]{Tri83}.
Applying $\dot{J}^{\nu}$ to this sequence we obtain a sequence in $\dot{J}^{\nu}(\CicRn)$ approximating $g$ both in $L^{p_0}(\Rn)$ as well as $L^{p_1}(\Rn)$.
Therefore, both approximations in $L^{p_0}(\Rn)$ as well as $L^{p_1}(\Rn)$ can be compared to this specific sequence, which yields the well-definedness with the above arguments.

It remains to show that the inequaltities for the quasi-linearity and boundedness are fulfilled on all of $L^{p_0}(\Rn)\cup L^{p_1}(\Rn)$.
Let  $g\in L^{p_{i_1}}(\Rn)$ and $\tilde{g}\in L^{p_{i_2}}(\Rn)$ such that $g+\tilde{g}\in L^{p_i}(\Rn)$ for some choice of $i,i_1,i_2\in\lbrace 0,1\rbrace$.
When choosing $(g_n)_{n\in\N}$ and $\tilde{g}_{n\in\N}$ as the specific approximation sequences for $g$ and $g_n$ as described above, then by linearity $(g_n + \tilde{g}_n)_{n\in\N}\subset \dot{J}^{\nu}(\CicRn)$ is that specific sequence for $g_n + \tilde{g}_n$. 
Because of the $L^{p_{i_0}}$, $L^{p_{i_1}}$ and $L^{p_i}$ convergence, there exists a shared subsequence such that $Tg_{n_k}$, $T\tilde{g}_{n_k}$ and $T(g_{n_k}+\tilde{g}_{n_k})$ converge against $Tg$, $T\tilde{g}$ and $T(g+\tilde{g})$ respectively almost everywhere. 
Hence,
\begin{equation*}
	|T(g+\tilde{g})| 
	= 		\lim_{k\rightarrow\infty} |T(g_{n_k}+\tilde{g}_{n_k})| 
	\leq		\lim_{k\rightarrow\infty} \left(|Tg_{n_k}|+|T\tilde{g}_{n_k}| \right)
	=		|Tg|+|T\tilde{g}|
\end{equation*}
almost everywhere.
Regarding the boundedness of $T$, we have
\begin{equation*}
	\|Tg\|_{L^{p_{i_1}}(\Rn)} = \lim_{n\rightarrow\infty} \|Tg_n\|_{L^{p_{i_1}}(\Rn)} 
		\lesssim	\lim_{n\rightarrow\infty} \|g_n\|_{L^{p_{i_1}}(\Rn)} = \|g\|_{L^{p_{i_1}}(\Rn)}.
\end{equation*}

In total, the slightly modified requirements for Theorem \ref{th:tt lorentz interpolation} are met and for any $p\in(p_0,p_1)$, $0<q\leq\infty$ and $f\in\CicRn$ we obtain
\begin{equation*}
	\left\| x\mapsto \left(\int_0^\infty |t^{\frac{1}{2}-\nu}\nabla_x F^s(x,t)|^2\,dt\right)^\frac{1}{2}\right\|_{\Lor{p}{q}{\Rn}}
	=\| T(J^\nu f) \|_{\Lor{p}{q}{\Rn}} 
	\lesssim \|J^\nu f\|_{\Lor{p}{q}{\Rn}}.
\end{equation*}
Since $1<p_0,p_1<\infty$ can be chosen arbitrarily, this is (\ref{eq:tt Bessel potential spaces via  nabla_xFs}).
\end{proof}

\begin{remark}\label{re:tt bessel potential space estimates for larger parameters}
Since the estimates in Proposition \ref{prop:tt bessel potential spaces char} are derived from Triebel-Lizorkin space estimates, there is a variety of additional weaker estimates which can analogously be obtained from some Triebel-Lizorkin space properties. 
Most important to us will be the following example.
Let $k\geq 2$. 
Due to Proposition \ref{prop:tt elementary embeddings for trilibeli}, under the assumptions for (\ref{eq:tt Bessel potential spaces via  nabla_xFs}) we have
\begin{equation*}
	\left\|x\rightarrow	\left( \int_0^\infty |t^{1-\frac{1}{k}-\nu}\nabla_x P^s_tf(x)|^k \,dt\right)^\frac{1}{k}\right\|_{L^p(\Rn)}
	\lesssim\|f\|_{\dot{F}^\nu_{p,k}(\Rn)} \lesssim \|f\|_{\dot{F}^\nu_{p,2}(\Rn)}.
\end{equation*}
Following the proof of Proposition \ref{prop:tt bessel potential spaces char}, we then obtain
\begin{equation*}
	\left\|x\rightarrow	\left( \int_0^\infty |t^{1-\frac{1}{k}-\nu}\nabla_x P^s_tf(x)|^k \,dt\right)^\frac{1}{k}\right\|_{\Lor{p}{q}{\Rn}}
	\lesssim  \|\fL{\nu}f\|_{\Lor{p}{q}{\Rn}}.
\end{equation*}
Of course, analogous versions hold for the other estimates of the proposition.
\end{remark}

Next, we cover the space of all functions with bounded mean oscillation, known as the $BMO$-space. 
The $BMO$-seminorm is given by 
\begin{equation}\label{eq:tt def BMO}
	[f]_{BMO} = \sup_{B\subset\Rn} \frac{1}{|B|}\int_B|f(y)-(f)_B|\,dy,
\end{equation}
where the supremum is over all balls $B\subset\Rn$ and $(f)_B\coloneqq |B|^{-1}\int_B f$.
We can then define the space $BMO = \lbrace f\in L^1_{\text{loc}}(\Rn)\,;\, [f]_{BMO}<\infty \rbrace$.
The following result is a special case of \cite[Chapter IV, §4.3, p.159, Theorem 3]{Ste93}, a well known result covering the relation between $BMO$ and certain sets of so called Carleson measures on $\Rnpp$. 
Stein already remarked that the result applies to certain Poisson-like kernels, see \cite[Chapter IV, §4.4.3, p.165]{Ste93}.
For us, this special case is another easy consequence from the Triebel-Lizorkin space characterizations.

\begin{proposition}[BMO characterization]\label{prop:tt bmo char}
Let $s\in(0,2)$.
Denote with $F^s(x,t)=P^s_tf(x)$ the $s$-harmonic extension of $f\in\CicRn$. Then
\begin{equation}\label{eq:tt BMO characterization with nabla_Rnp}
	[f]_{BMO(\Rn)} \gtrsim
	\sup_{x\in\Rn, r>0}\left(\frac{1}{|B(x,r)|}\int_{T(B(x,r))} t|\nabla_\Rnp F^s(y,t)|^2 \,dy\,dt\right)^\frac{1}{2}.
\end{equation}
Additionally, for $\beta>0$
\begin{equation}\label{eq:tt BMO characterization with fracLapl(beta)}
	[f]_{BMO(\Rn)} \approx
	\sup_{x\in\Rn, r>0}\left(\frac{1}{|B(x,r)|}\int_{T(B(x,r))} t^{2\beta-1}|P^s_t(\fL{\beta}f)(y)|^2 \,dy\,dt\right)^\frac{1}{2}.
\end{equation}
and
\begin{equation}\label{eq:tt BMO characterization with nabla_Rnpp fracLapl(beta)}
	[f]_{BMO(\Rn)} \gtrsim
	\sup_{x\in\Rn, r>0}\left(\frac{1}{|B(x,r)|}\int_{T(B(x,r))} t^{2\beta+1}|\nabla_\Rnp P^s_t(\fL{\beta}f)(y)|^2 \,dy\,dt\right)^\frac{1}{2}.
\end{equation}
Here, $T(B(x,r))=\lbrace (y,t)\in\Rnpp\, :\,|y-x|<r-t\rbrace$ is the ``tent'' over $B(x,r)$.
\end{proposition}

\begin{proof}
Recalling the list from the beginning of subsection \ref{subsec:tt - triebel lizorkin besov}, we have $BMO=\dot{F}^0_{\infty,2}$.
In the proof of Theorem \ref{th:tt triebel lizorkin besov char} we have already seen that the kernels $\partial_t p^s_1$, $\partial_{x_i}p^s_1$ and $\fL{\beta}p^s_1$, $\beta\in(0,1]$ satisfy the conditions of Theorem \ref{th:tt BuiCandy main result 1} with $\alpha=0$ and allow the choice of $\rho = 0$ for the appearing polynomial. 
Thus, we have (\ref{eq:tt BMO characterization with nabla_Rnp}) and (\ref{eq:tt BMO characterization with fracLapl(beta)}) for $\beta\in(0,1]$ as an immediate consequence of Theorem \ref{th:tt BuiCandy main result 1}.

Sadly, Theorem \ref{th:tt lifting property of trlibeli spaces}, the lifting property, is only applicable to Triebel-Lizorkin spaces with $p<\infty$. 
Therefore, we can not argue as in Remark \ref{re:tt annotations regarding the trilibeli characterizations} to obtain estimates for higher order kernels. 
Following the same line of argument as for the proof Theorem \ref{th:tt triebel lizorkin besov char} though, we get that the kernels $\fL{\beta}p^s_1$, $\beta>1$, as well as $\fL{\beta}\partial_t p^s_1$ and $\fL{\beta}\partial_{x_i} p^s_1$ also satisfy the conditions of Theorem \ref{th:tt BuiCandy main result 1} for $\alpha=0$, $p=\infty$ and $q=2$.
This results in the rest of the estimates.
\end{proof}

By plugging suitable functions in these $BMO$-estimates, we obtain some additional estimates.

\begin{corollary}\label{cor:tt BMO(fLnu) characterization}
Let $s\in(0,2)$.
Denote with $F^s(x,t)=P^s_tf(x)$ the $s$-harmonic extension of $f\in\CicRn$. 
Then, for $0<\nu<1$ ,
\begin{equation}\label{eq:tt BMO(fLnu) characterization with nabla_x}
	[\fL{\nu}f]_{BMO(\Rn)} \gtrsim
	\sup_{x\in\Rn, r>0}\left(\frac{1}{|B(x,r)|}\int_{T(B(x,r))} t^{1-2\nu}|\nabla_y F^s(y,t)|^2 \,dy\,dt\right)^\frac{1}{2},
\end{equation}
and for $0<\nu<s$,
\begin{equation}\label{eq:tt BMO(fLnu) characterization with partial_t}
	[\fL{\nu}f]_{BMO(\Rn)} \approx
	\sup_{x\in\Rn, r>0}\left(\frac{1}{|B(x,r)|}\int_{T(B(x,r))} t^{1-2\nu}|\partial_tF^s(y,t)|^2 \,dy\,dt\right)^\frac{1}{2}.
\end{equation}
Additionally, for $0<\nu<\min\lbrace 1+s,2\rbrace$,
\begin{equation}\label{eq:tt BMO(fLnu) characterization with double deriv}
	[\fL{\nu}f]_{BMO(\Rn)} \gtrsim
	\sup_{x\in\Rn, r>0}\left(\frac{1}{|B(x,r)|}\int_{T(B(x,r))} t^{3-2\nu}|\nabla_x\nabla_\Rnp F^s(y,t)|^2 \,dy\,dt\right)^\frac{1}{2}.
\end{equation}
\end{corollary}

\begin{proof}
We obtain (\ref{eq:tt BMO(fLnu) characterization with nabla_x}) by plugging $\mathcal{R}_i\fL{\nu}f$ into (\ref{eq:tt BMO characterization with fracLapl(beta)}) with $\beta=1-\nu$ since
\begin{equation*}
	P^s_t(\fL{1-\nu}\mathcal{R}_i\fL{\nu}f)(x) = P^s_t(\mathcal{R}_i\fL{1}f)(x) 
	= P^s_t(\partial_{x_i}f)(x) = \partial_{x_i}F^s(x,t).
\end{equation*}
This can be easily seen on the side of the Fourier transforms.
Thus, we have (\ref{eq:tt BMO(fLnu) characterization with nabla_x}) since the Riesz transform maps $BMO$ continuously onto itself, see \cite[Chapter IV, §6.3a(b), p.179]{Ste93}.

Regarding (\ref{eq:tt BMO(fLnu) characterization with partial_t}), we have
\begin{equation*}
	t^{1-\nu}\partial_t P^s_tf(x) = Ct^{s-\nu} \left(\fL{s-\nu}p^{2-s}_t\ast\fL{\nu}f\right)(x)
\end{equation*}
For a detailed proof of this fact, see the proof of Estimate (\ref{eq:tt square function with partial_t}), for the interaction of dilations and fractional Laplacians see Step 4 in the proof of Theorem \ref{th:tt triebel lizorkin besov char}.
Therefore, plugging $\fL{\nu}f$ into (\ref{eq:tt BMO characterization with fracLapl(beta)}) with $\beta = s-\nu$ and the $(2-s)$-harmonic extension, for the integrand we have
\begin{mathex}[LCL]
	t^{2(s-\nu)-1}|P^{2-s}_t(\fL{s-\nu}\fL{\nu}f)(y)|^2 =  Ct^{2s-2\nu-1}|t^{1-s}\partial_tP^s_tf(y)|^2.
\end{mathex}
Thus, we obtain (\ref{eq:tt BMO(fLnu) characterization with partial_t}).

In the same manner, we obtain (\ref{eq:tt BMO(fLnu) characterization with double deriv}) from (\ref{eq:tt BMO characterization with fracLapl(beta)}). 
By plugging in $\mathcal{R}_i\Rj\fL{\nu}f$ for $\beta = 2-\nu$, we obtain the estimate for $\partial_{x_i}\partial_{x_j}F^s$.
By plugging in $\mathcal{R}_i\fL{\nu}f$ for $\beta = 1+s-\nu$ and the $(2-s)$-harmonic extension, we obtain the estimate for $\partial_{x_i}\partial_tF^s$.
\end{proof}

\begin{remark}\label{re:tt equivalencies vs estimates}
All estimates in Proposition \ref{prop:tt slobodeckij space char}, Proposition \ref{prop:tt bessel potential spaces char}, and Proposition \ref{prop:tt bmo char} are actually listed as (norm-)equivalencies in \cite{LeSchi20}.
We only obtain the one sided estimates for two reasons.
First, for the kernels involving $x$-derivatives, this is due to the same reason we explained in Remark \ref{re:tt annotations regarding the trilibeli characterizations}(i).
Second, for the Lorentz scale estimates, results in the style of the classical interpolation theorems such as listed in \cite[Section 3]{Hun66} only yield one direction of the estimate.
Since proving the blackbox estimates in Subsection \ref{subsec:tt - blackbox estimates} and therefore also the main part only requires this one direction of the space characterizations, we do not further investigate the reverse estimates.
\end{remark}

We finish this subsection by covering the Hölder spaces.
For $\nu>0$, the corresponding Hölder-seminorm is 
\begin{equation*}
	[f]_{\C^\nu(\Rn)} \coloneqq
	\left\lbrace \begin{array}{llll}
		 \sup_{x\neq y\in\Rn}\frac{|\nabla^{\lfloor\nu\rfloor}f(x)-\nabla^{\lfloor\nu\rfloor}f(y)|}{|x-y|^{\nu-\lfloor\nu\rfloor}}
		 															\quad&	\text{if } \nu\notin\N, \\
		 \|\nabla^\nu f\|_{L^\infty(\Rn)}							\quad&	\text{if }\nu\in\N,
		
	\end{array}\right.
\end{equation*}
where $\lfloor \cdot \rfloor$ is the floor function.
We may denote $[f]_{\operatorname{Lip}}=[f]_{\C^1}$.
The Hölder spaces are defined as 
\begin{equation*}
	\C^\nu(\Rn) = 	\lbrace f\in C^{\lfloor\nu\rfloor}(\Rn) \,;\, 
						\|f\|_{\C^\nu(\Rn)} = \|f\|_{\C^{\lfloor\nu\rfloor}(\Rn)} + [f]_{\C^\nu(\Rn)}<\infty 
					\rbrace.
\end{equation*}

\begin{proposition}[Hölder space characterizations]\label{prop:tt hölder char}
Let $s\in(0,2)$. Denote with $F^s(x,t)=P^s_tf(x)$ the $s$-harmonic extension of $f\in\CicRn$.
Then for $\nu\in(0,s)\setminus\N$ we have
\begin{equation}\label{eq:tt Hölder characterization with partial_t}
	\sup_{(x,t)\in\Rnpp} t^{1-\nu}|\partial_t F^s(x,t)| \approx [f]_{\C^\nu(\Rn)}.
\end{equation}
Additionally, for any $\nu\in(0,1]$,
\begin{equation}\label{eq:tt Hölder characterization with nabla_x}
	\sup_{(x,t)\in\Rnpp} t^{1-\nu}|\nabla_x F^s(x,t)| \lesssim [f]_{\C^\nu(\Rn)}.
\end{equation}
Regarding double derivatives, for $\nu\in (0,\min\lbrace 1+s, 2\rbrace)\setminus \N$ we have
\begin{equation}\label{eq:tt Hölder characterization with nabla_Rnp nabla_x}
	\sup_{(x,t)\in\Rnpp} t^{2-\nu}|\nabla_x \nabla_{\Rnp}F^s(x,t)| \lesssim [f]_{\C^\nu(\Rn)}.
\end{equation}
\end{proposition}

\begin{proof}
According to the Besov-Lipschitz space identifications listed at the start of subsection \ref{subsec:tt - triebel lizorkin besov}, we have $\C^\nu(\Rn)= B^\nu_{\infty,\infty}(\Rn)$ for $\nu\notin\N$.
According to \cite[2.3.3(i)]{Tri92} we also have 
\begin{equation*}
	\|f\|_{B^\nu_{\infty,\infty}} \approx \|f\|_{L^\infty(\Rn)} + \|f\|_{\dot{B}^\nu_{\infty,\infty}}.
\end{equation*}
Therefore,
\begin{equation}
	\|f\|_{L^\infty(\Rn)} + \|f\|_{\dot{B}^\nu_{\infty,\infty}} \approx \|f\|_{\C^0(\Rn)} + [f]_{\C^\nu(\Rn)}.
\end{equation}
For $\lambda>0$ directly via the definition of the (semi-)norms we obtain
\begin{mathex}
	\|f(\lambda \cdot)\|_{L^\infty(\Rn)} &=& \|f\|_{L^\infty(\Rn)}\\
	\|f(\lambda\cdot)\|_{\C^0(\Rn)} &=& \|f\|_{\C^0(\Rn)}\\
	{[}f(\lambda\cdot){]}_{\C^\nu(\Rn)} &=& \lambda^\nu[f]_{\C^\nu(\Rn)}.
\end{mathex}
For the homogeneous Besov-Lipschitz norm, due to \cite[Remark 5.1.3.4]{Tri83}, we have 
\begin{equation*}
	\|f(\lambda\cdot)\|_{\dot{B}^\nu_{\infty,\infty}} \approx \lambda^\nu \|f\|_{\dot{B}^\nu_{\infty,\infty}}.
\end{equation*}
With the same argument as in the proof of the estimates (\ref{eq:tt Slobodeckij via nabla_xFs})-(\ref{eq:tt Slobodeckij via partial_tFs}) we get the identification $\|f\|_{\dot{B}^\nu_{\infty,\infty}}\approx [f]_{\C^\nu(\Rn)}$.
Thus, for $\nu\notin\N$ the estimates (\ref{eq:tt Hölder characterization with partial_t}) and (\ref{eq:tt Hölder characterization with nabla_x}) follow immediately from Theorem \ref{th:tt triebel lizorkin besov char} while (\ref{eq:tt Hölder characterization with nabla_Rnp nabla_x}) follows with Remark \ref{re:tt annotations regarding the trilibeli characterizations}(iv).

It remains to show (\ref{eq:tt Hölder characterization with nabla_x}) for $\nu=1$.
The estimate follows immediately from Corollary \ref{cor:Pst boundedness of Pst derivatives} since $[f]_{\C^1(\Rn)}= [f]_{\operatorname{Lip}} = \|\nabla f \|_{L^\infty(\Rn)}$.
\end{proof}

\subsection{Building blocks: Square function estimates}\label{subsec:tt - square function estimates}
Next to all the estimates and characterizations from the previous subsection, which are based on the $\dot{F}^\alpha_{p,q}$ and $\dot{B}^\alpha_{p,q}$ characterizations, we also need some estimates that originate from the theory of so-called nontangential square functions.
These are related to singular integrals and maximal functions, confer for example \cite[Chapter I, §6.3]{Ste93}.

Let $\phi$ on $\Rn$ be a kernel sufficiently small at infinity (e.g., $\phi \in\mathcal{S}(\Rn)$) with $\int_\Rn \phi \,dx= 0$.
Then the (regular) square function $s_\phi$ is defined as 
\begin{equation*}
	f\mapsto s_\phi(f),\qquad (s_\phi f)(x) = \left( \int_0^\infty |(f\ast \phi_t)(x)|^2\frac{dt}{t}\right)^\frac{1}{2},
\end{equation*}
while the so-called ``nontangential'' square function $S_\phi$ is given by
\begin{equation*}
	f\mapsto S_\phi(f),\qquad (S_\phi f)(x) = \left( \int_{\lbrace(y,t);|x-y|<t\rbrace} |(f\ast \phi_t)(y)|^2\frac{dt\,dy}{t^{n+1}}\right)^\frac{1}{2}.
\end{equation*}
Again, $\phi_t$ denotes the standard dilation.
For these square functions we find an $L^p$-estimate in \cite[Chapter I, §8.23]{Ste93}, which we can extend to a Lorentz estimate.

\begin{theorem}\label{th:sqfunc Lpq estimate}
Let $1<p<\infty$, $0<q\leq\infty$.
Assume $\phi\in\C^1(\Rn)$ with $\int_\Rn \phi \,dx = 0$ and $\sigma_1,\sigma_2>0$, such that for all $x\in\Rn$
\begin{equation*}
	|\phi(x)|\leq C(1+|x|)^{-n-\sigma_1}, \qquad |\nabla \phi(x)|\leq C(1+|x|)^{-n-\sigma_2}.
 \end{equation*} 
Then for all $f\in \Lor{p}{q}{\Rn}$ we have
\begin{equation*}
	\|s_\phi(f)\|_{\Lor{p}{q}{\Rn}} \lesssim \|f\|_{\Lor{p}{q}{\Rn}},
\end{equation*}
and
\begin{equation*}
	\|S_\phi(f)\|_{\Lor{p}{q}{\Rn}} \lesssim \|f\|_{\Lor{p}{q}{\Rn}}.
\end{equation*}
\end{theorem}

\begin{proof}
The result from \cite{Ste93} is originally formulated for kernels $\phi\in\mathcal{S}$, the proof in §8.23 bulding upon observations from §6.3 and §6.4.
According to \cite[Chapter I, §8.23(c)]{Ste93}, the decay condition can be relaxed for example to $|\phi(x)|\leq C(1+|x|)^{-n-1}$ and $|\nabla \phi(x)|\leq C(1+|x|)^{-n-1}$.
In fact, the conditions we stated in the Theorem suffice. 
First, we list the necessary adjustments to the proof in \cite{Ste93}. 
Second, we apply Theorem \ref{th:tt lorentz interpolation} to obtain the full Lorentz scale. 

\underline{Step 1:}
The part of the proof given in §8.23 is unchanged, while in §6.3 and §6.4 some small adjustments are necessary. 
In §6.4, instead of $\phi_t(x)\leq At |x|^{-n-1}$ we obtain $\phi_t(x)\leq At^{\sigma_1} |x|^{-n-\sigma_1}$.
In §6.3, instead of $\Fourier(\phi)(\xi)\leq A|\xi|$ we have  $\Fourier(\phi)(\xi)\leq A|\xi|^{\tilde{\sigma}}$ for $0<\tilde{\sigma}<\min\lbrace 1,\sigma_1\rbrace$ because $\Fourier(\phi)(0) = 0$ and $\Fourier(\phi)\in\C^{\tilde{\sigma}}$.
The former is a direct consequence of the integral over $\phi$ vanishing.
We obtain the latter via Triebel-Lizorkin and Besov-Lipschitz space identifications. 
For $\sigma<\sigma_1$ we have
\begin{equation*}
	J^\sigma\Fourier(\phi) = C\Fourier((1+|x|)^\sigma \phi)\in L^p \quad \text{for all } 2\leq p<\infty,
\end{equation*}
because $(1+|x|)^\sigma \phi\in L^1\cap L^2$ if $|\phi(x)|\leq C(1+|x|)^{-n-\sigma_1}$.
Since $L^p=F^0_{p,2}$, with Theorem \ref{th:tt lifting property of trlibeli spaces} we obtain $\Fourier(\phi)\in F^\sigma_{p,2}$ for all $2\leq p<\infty$.
We then have $\Fourier(\phi)\in F^\sigma_{p,2}\subset B^{\sigma}_{p,p}$ due to Proposition \ref{prop:tt elementary embeddings for trilibeli}(iii). 
In order to manipulate the other parameters, we use a generalization of the classical Sobolev embedding theorem for Triebel-Lizorkin and Besov-Lipschitz spaces.
With \cite[Theorem 2.7.1]{Tri83} we get $\Fourier(\phi)\in B^{\tilde{\sigma}}_{\infty,p}$ for all $\tilde{\sigma}<\sigma$. 
Finally, Proposition \ref{prop:tt elementary embeddings for trilibeli}(i) yields $\Fourier(\phi)\in B^{\tilde{\sigma}}_{\infty,p}\subset B^{\sigma}_{\infty,\infty}=\C^\sigma$.

\underline{Step 2:}
The operators $s_\phi$ and $S_\phi$ are obviously quasi-linear due to the Minkowski inequality. 
With Theorem \ref{th:tt lorentz interpolation} we then obtain the above result, since we already have the boundedness regarding $L^p=L^{(p,p)}$ for any $1<p<\infty$. 
\end{proof}

Applying the non-tangential square function estimate to various modifications of our Poisson kernel, we now obtain a variety of estimates. 
These are interesting, since they occur as counterpart to the integrals equivalent to the $BMO$-space norm, see Lemma \ref{lem:tt carleson measure and square function estimate}.
The regular square function estimate also provides some estimates, but these are already covered via the Triebel-Lizorkin space characterizations.

\begin{proposition}[Square function estimates]\label{prop:tt square functions to Bessel potential spaces}
Let $1<p<\infty$, $0< q \leq \infty$ and $s\in(0,2)$.
Denote with $F^s(x,t)=P^s_tf(x)$ the $s$-harmonic extension of $f\in\CicRn$.\\
Regarding classical derivatives, for $\nu \in [0,s)$ we have
\begin{equation}\label{eq:tt square function with partial_t}
	\left\| x \mapsto 
		\left( \int_{\lbrace (y,t);|y-x|<t\rbrace}	t^{1-2\nu-n}|\partial_t F^s(y,t)|^2 \,dy\,dt\right)^\frac{1}{2}
	\right\|_{\Lor{p}{q}{Rn}}	
	\lesssim \|\fL{\nu}f\|_{\Lor{p}{q}{\Rn}},
\end{equation}
and for $\nu\in[0,1)$
\begin{equation}\label{eq:tt square function with nabla_x}
	\left\| x \mapsto 
		\left( \int_{\lbrace (y,t);|y-x|<t\rbrace}	t^{1-2\nu-n}|\nabla_x F^s(y,t)|^2 \,dy\,dt\right)^\frac{1}{2}
	\right\|_{\Lor{p}{q}{Rn}}	
	\lesssim \|\fL{\nu}f\|_{\Lor{p}{q}{\Rn}}.
\end{equation}
Moreover, for any $\nu\in[0,\min\lbrace 1+s,2\rbrace)$ we have 
\begin{equation}\label{eq:tt square function with double derivative}
	\left\| x \mapsto 
		\left( \int_{\lbrace (y,t);|y-x|<t\rbrace}	t^{3-2\nu-n}|\nabla_x\nabla_{\Rnp} F^s(y,t)|^2 \,dy\,dt\right)^\frac{1}{2}
	\right\|_{\Lor{p}{q}{Rn}}	
	\lesssim \|\fL{\nu}f\|_{\Lor{p}{q}{\Rn}}.
\end{equation}
Regarding fractional derivatives, for $\beta>0$ and $\nu\in[0,\beta)$ we have
\begin{equation}\label{eq:tt square function with fracLapl}
	\left\| x \mapsto 
		\left( \int_{\lbrace (y,t);|y-x|<t\rbrace}	t^{2\beta-2\nu-n-1}| P^s_t(\fL{\beta} f)(y)|^2 \,dy\,dt\right)^\frac{1}{2}
	\right\|_{\Lor{p}{q}{Rn}}	
	\lesssim \|\fL{\nu}f\|_{\Lor{p}{q}{\Rn}}.
\end{equation}
All these estimates also hold for $-n<\nu<0$ when replacing $\|\fL{\nu}f\|_{\Lor{p}{q}{\Rn}}$ with $\|I^{-\nu}f\|_{\Lor{p}{q}{\Rn}}$.
\end{proposition}

\begin{proof}[Estimate (\ref{eq:tt square function with partial_t})]
In order to show (\ref{eq:tt square function with partial_t}), we only need to bring the left side of the estimate in the form of a nontangential square function and verify that the respective kernels fulfill the conditions of Theorem \ref{th:sqfunc Lpq estimate}. 
We have 
\begin{equation*}
	\left( \int_{\Gamma(x)}	t^{1-2\nu-n}|\partial_t F^s(y,t)|^2 \,dy\,dt\right)^\frac{1}{2}
	= 
	\left( \int_{\Gamma(x)}	t^{-1-n}|t^{1-\nu}(\partial_t p^s_t*f)(y)|^2 \,dy\,dt\right)^\frac{1}{2}
\end{equation*}
with $\Gamma(x)=\lbrace (y,t)\,;\,|y-x|<t\rbrace$.
In order to obtain the estimate, we need to find a kernel $q$ that satisfies the conditions of Theorem \ref{th:sqfunc Lpq estimate} and fulfills
\begin{equation*}
	t^{1-\nu}\partial_t p^s_t*f = q_t*\fL{\nu}f.
\end{equation*} 
Applying the Fourier transform to the left side and recalling (\ref{eq:tt partial_t pst Fourier}), we obtain
\begin{mathex}
	\Fourier(t^{1-\nu}\partial_t p^s_t*f)(\xi)
	&=&	C t^{2-\nu}|\xi|^2 \left(\int_0^\infty \lambda^{\frac{s-2}{2}} e^{-\lambda-\frac{|t\xi|^2}{c\lambda}} \frac{d\lambda} {\lambda}\right)\Fourier(f)(\xi)\\
	&=&	C|t\xi|^{2-\nu}\left(\int_0^\infty \lambda^{\frac{s-2}{2}} e^{-\lambda-\frac{|t\xi|^2}{c\lambda}} \frac{d\lambda} {\lambda}\right)\Fourier(\fL{\nu}f)(\xi).
\end{mathex}
Taking into consideration the interaction between Fourier transform and dilation, the candidate for a fitting kernel is given by
\begin{equation}\label{eq:sqfunc kernel candidate for partial_t}
	q(x) = \Fourier^{-1}\left(\xi\mapsto C|\xi|^{2-\nu}\left(\int_0^\infty \lambda^{\frac{s-2}{2}} e^{-\lambda-\frac{|\xi|^2}{c\lambda}} \frac{d\lambda} {\lambda}\right)\right)
\end{equation}
for some constant $C\neq 0$. 
Substituting $\tilde{\lambda}= \frac{|\xi|^2}{c\lambda}$, we obtain
\begin{mathex}
	C|\xi|^{2-\nu}\left(\int_0^\infty \lambda^{\frac{s-2}{2}} e^{-\lambda-\frac{|\xi|^2}{c\lambda}} \frac{d\lambda} {\lambda}\right)
	&=&	C|\xi|^{s-\nu} \left(\int_0^\infty \left(\frac{|\xi|^2}{c\lambda}\right)^{\frac{2-s}{2}} e^{-\lambda-\frac{|\xi|^2}{c\lambda}} \frac{d\lambda} {\lambda}\right)\\
	&=& C|\xi|^{s-\nu} \left(\int_0^\infty \tilde{\lambda}^{\frac{2-s}{2}} e^{-\frac{|\xi|^2}{c\tilde{\lambda}}-\tilde{\lambda}} \frac{d\tilde{\lambda}} {\tilde{\lambda}}\right).
\end{mathex}
Therefore, due to \cite[Proposition 7.6(a)]{Hao16} and the definition of the fractional Laplacian, an explicit formula for the kernel is
\begin{equation}\label{eq:tt fitting kernel for partial_t ps1 square function estimate}
	q(x) = C\fL{s-\nu}\left((|x|^2+1)^{-\frac{s-2-n}{2}}\right).
\end{equation}
Combining Lemma \ref{lem:fL Linfty estimate} and Lemma \ref{lem:fL decay}, the $L^\infty$-estimate and the decay estimate for the fractional Laplacian,
\begin{equation*}
	|q(x)| \leq C \frac{1}{(1+|x|)^{-n-s+\nu}},
	\qquad
	|\nabla_x q(x)| \leq C\frac{1}{(1+|x|)^{-n-s+\nu}}.
\end{equation*}
Because of $\lim_{\xi\rightarrow 0}\Fourier(q)(\xi)=0$ and the continuity of $\Fourier(q)$, we further have $\int_\Rn q\,dx= 0$, compare Step 3 of the proof for Theorem \ref{th:tt triebel lizorkin besov char}.
Since $\nu<s$ and therefore $s-\nu>0$, the conditions for the square function estimate are fulfilled.
Thus, we have shown (\ref{eq:tt square function with partial_t}).

The same arguments hold for $-n<\nu<0$, where we use that $\Fourier(I^{-\nu}f)(\xi)=C|\xi|^{\nu}$.
\end{proof}

\begin{remark}
Based on (\ref{eq:sqfunc kernel candidate for partial_t}), one might be tempted to immediately assume 
\begin{equation*}
	q(x) = C\fL{2-\nu}\left((1+|x|^2)^{-\frac{n+s-2}{2}}\right).
\end{equation*}
This would only work for $n\geq 2$ though, because for \cite[Proposition 7.6(a)]{Hao16} we need that $n+s-2>0$
Instead, we basically use the following interesting interaction between fractional Laplacian and Poisson kernel on the side of the Fourier transforms.
For $0<\gamma<n$ we have
\begin{equation*}
	\fL{\gamma}(|x|^2+1)^\frac{\gamma-n}{2} = C (|x|^2+1)^\frac{-\gamma-n}{2}.
\end{equation*}
This interaction can easily be proven on the side of the Fourier transforms with the same transformations we used to obtain the explicit form of $q$. 
Also, on the Fourier side, we can choose $\gamma\in\R$.
\end{remark}

\begin{proof}[Estimate (\ref{eq:tt square function with nabla_x})]
This estimate is obtained analogously to (\ref{eq:tt square function with partial_t}), starting from
\begin{mathex}
	\Fourier(t^{1-\nu}\partial_{x_i} p^s_t*f)(\xi)
	&=&	C t^{1-\nu}\xi_i\left(\int_0^\infty \lambda^{\frac{s}{2}} e^{-\lambda-\frac{|t\xi|^2}{c\lambda}} \frac{d\lambda} {\lambda}\right)\Fourier(f)(\xi)\\
	&=&	C |t\xi|^{1-\nu}\frac{t\xi_i}{|t\xi|}\left(\int_0^\infty \lambda^{\frac{s}{2}} e^{-\lambda-\frac{|t\xi|^2}{c\lambda}} \frac{d\lambda} {\lambda}\right)\Fourier(\fL{\nu}f)(\xi).\\
\end{mathex}
Recalling the Fourier symbol of the Riesz transforms, we then get the kernel
\begin{equation}\label{eq:tt fitting kernel for partial_xi ps1 square function estimate}
	q(x)= C \fL{1-\nu}\mathcal{R}_i  \left((1+|x|^2)^{-\frac{n+s}{2}}\right).
\end{equation}
The conditions of Theorem \ref{th:sqfunc Lpq estimate} are verified via the same arguments as above when additionally taking into account Lemma \ref{lem:Rt decay} and Lemma \ref{lem:Rt Linfty estimate}, the decay estimate and the $L^\infty$-estimate for the Riesz transforms.
\end{proof}

\begin{proof}[Estimate (\ref{eq:tt square function with double derivative})]
For this estimate we have two different types of kernels.
Analogous to the proofs of the other two estimates, we start with
\begin{equation*}
	\Fourier(t^{2-\nu}\partial_{x_i}\partial_{x_j} p^s_t*f)(\xi)
	=	C |t\xi|^{2-\nu}\frac{t\xi_i}{|t\xi|}\frac{t\xi_j}{|t\xi|}\left(\int_0^\infty \lambda^{\frac{s}{2}} e^{-\lambda-\frac{|t\xi|^2}{c\lambda}} \frac{d\lambda} {\lambda}\right)\Fourier(\fL{\nu}f)(\xi)
\end{equation*}
and
\begin{equation*}
	\Fourier(t^{2-\nu}\partial_{x_i}\partial_{t} p^s_t*f)(\xi)
	=	C \frac{t\xi_i}{|t\xi|} |t\xi|^{1+s-\nu}|t\xi|^{2-s}\left(\int_0^\infty \lambda^{\frac{s-2}{2}} e^{-\lambda-\frac{|t\xi|^2}{c\lambda}} \frac{d\lambda} {\lambda}\right)\Fourier(\fL{\nu}f)(\xi).
\end{equation*}
We then obtain 
\begin{equation*}
	q(x)= C \fL{2-\nu}\mathcal{R}_i\Rj  \left((1+|x|^2)^{-\frac{n+s}{2}}\right),
\end{equation*}
and
\begin{equation*}
	q(x)= C \fL{1+s-\nu}\mathcal{R}_i  \left((1+|x|^2)^{-\frac{n+2-s}{2}}\right)
\end{equation*}
respectively.
Following the same line of argument as for the other two estimates, both types of kernels fulfill the conditions of Theorem \ref{th:sqfunc Lpq estimate}.
\end{proof}

\begin{proof}[Estimate (\ref{eq:tt square function with fracLapl})]
Again, we first need to find a suitable kernel, so that the left side of the estimate is in the form of a nontangential square function.
We have 
\begin{equation*}
	\Fourier(t^{\beta-\nu}P^s_t(\fL{\beta} f))(\xi) = C |t\xi|^{\beta-\nu}\, \Fourier((p^s_1)_t)\,\Fourier(\fL{\nu}f).
\end{equation*}
Therefore, a suitable kernel is 
\begin{equation*}
	q(x) = C \fL{\beta-\nu} \left((1+|x|^2)^{-\frac{n+s}{2}}\right).
\end{equation*}
Following the arguments from the proof of (\ref{eq:tt square function with partial_t}), this kernel satisfies the conditions of Theorem \ref{th:sqfunc Lpq estimate}.
\end{proof}

\begin{corollary}
Let $1<p<\infty$, $0< q \leq \infty$ and $s\in(0,2)$.
Denote with $F^s(x,t)=P^s_tf(x)$ the $s$-harmonic extension of $f\in\CicRn$.
Then for $0<\nu<n$ we have 
\begin{equation}\label{eq:tt square function with unmodified kernel}
	\left\| x \mapsto 
		\left( \int_{\lbrace (y,t);|y-x|<t\rbrace}	t^{2\nu-n-1}|F^s(y,t)|^2 \,dy\,dt\right)^\frac{1}{2}
	\right\|_{\Lor{p}{q}{Rn}}	
	\lesssim \|I^{\nu}f\|_{\Lor{p}{q}{\Rn}}.
\end{equation}
\end{corollary}

\begin{proof}
Insert $I^{\nu}f$ into (\ref{eq:tt square function with fracLapl}) with $\beta=\nu$.
\end{proof}

\begin{remark}
Same as in the subsection above, the estimates in Proposition \ref{prop:tt square functions to Bessel potential spaces} are listed as (norm-)equivalencies in \cite{LeSchi20}.
Here, in order to obtain equivalency for the square-function estimates, we would need to verify that the kernels are nondegenerate as described in \cite[Chapter I, §8.23(b)]{Ste93}, which is basically the Tauberian condition \textbf{(C2)}, see \cite[Chapter IV, §4.3, p.159]{Ste93}. 
And again, for the Lorentz scale estimates, our interpolation theorem only yields one direction of the estimate.
\end{remark}

\subsection{Building blocks: Maximal function and pointwise estimates of $P^s_tf$}\label{subsec:tt - maximal function estimates}

In addition to the Triebel-Lizorkin and Besov-Lipschitz space characterizations as well as the square function estimates, we further need pointwise estimates, which offer more than for example Corollary \ref{cor:Pst decay of Pst and derviatives} and Corollary \ref{cor:Pst boundedness of Pst derivatives}. 
The central vehicle for these pointwise estimates will be the maximal function, which we already briefly mentioned  at the beginning of the previous subsection.
For $f\in L^1_{\text{loc}}(\Rn)$ the maximal function $\mathcal{M}f$ is defined by 
\begin{equation}\label{eq:maxfunc - definition}
	(\mathcal{M}f)(x) \coloneqq \sup_{\delta>0} \frac{1}{|B(x,\delta)|}\int_{B(x,\delta)} |f(y)| \,dy.
\end{equation}

The idea is to obtain pointwise estimates of the $s$-harmonic extension against the maximal function of $f$ or derivatives of $f$ with the following proposition from \cite[Chapter II, §2.1]{Ste93}.

\begin{proposition}\label{prop:maxfunc - convolution estimate}
Assume that a measurable function $\psi\colon\Rn\rightarrow\R$ has a radial majorant which is non-increasing, bounded and integrable. 
Let $f\in L^1_{\text{loc}}(\Rn)$.
Then, for $F\colon\Rnpp\rightarrow \R$, $F(x,t)=(f\ast \psi_t)(x)$, where $\psi_t$ is the standard dilation of $\psi$, we have 
\begin{equation*}
	\sup_{(y,t)\,;\, |y-x|<t} |F(y,t)| \lesssim \mathcal{M}f(x).
\end{equation*}
\end{proposition}

These estimates are useful because in turn the Lorentz space norm of the maximal function can be estimated against the respective Lorentz space norm of the original function with the following theorem from \cite[Chapter I, §3.1]{Ste93}. 

\begin{theorem}\label{th:maxfunc - Lp estimate}
Let $f$ be a function defined on $\Rn$.
\begin{enumerate}[label=(\roman*)]
\item	If $f\in L^p(\Rn)$, $1\leq p\leq \infty$, then $\mathcal{M}f$ is finite almost everywhere.
\item	If $f\in\Lor{p}{q}{\Rn}$, $1< p< \infty$, $0< q\leq \infty$ or $p=q=\infty$, then $\mathcal{M}f\in \Lor{p}{q}{\Rn}$ and 
		\begin{equation*}
			\|\mathcal{M}f\|_{\Lor{p}{q}{\Rn}}\lesssim \|f\|_{\Lor{p}{q}{\Rn}}.
		\end{equation*}
\end{enumerate}
\end{theorem}

\begin{proof}
From \cite[Chapter I, §3.1]{Ste93} we obtain (i) and (ii) for $q=p$.
We then obtain the full Lorentz version of (ii) via interpolation with Theorem \ref{th:tt lorentz interpolation}. 
The maximal function is defined for all $f\in L^1_{\text{loc}} \supset L^{p_0}+L^{p_1}$, $1<p_0,p_1\leq\infty$. 
Obviously, the maximal function is quasi linear. 
Therefore, Theorem \ref{th:tt lorentz interpolation} is applicable.
\end{proof}

We will use this theorem in combination with the following pointwise estimates.
Besides the promised estimates for the unmodified $s$-harmonic extension, we obtain some additional pointwise estimates for $\partial_t P^s_t f$ and $\nabla_x P^s_t f$.
Note that we take the supremum over the same set, which already occured in the nontangential square functions, that is $\lbrace (y,t)\, ;\,|y-x|<t\rbrace$.

\begin{proposition}[Pointwise estimates with the maximal function]\label{prop:maxfunc pointwise estimates}
Let $s\in(0,2)$ and $f\in\CicRn$.
Denote with $F^s(x,t)=P^s_t f(x)$.
Then
\begin{equation}\label{eq:maxfunc pw est of Fs}
	\sup_{(y,t)\, ;\,|y-x|<t} |F^s(y,t)|\lesssim \mathcal{M}f(x).
\end{equation}
For $\sigma\in [0,s]$, regarding the derivative in $t$ direction,
\begin{equation}\label{eq:maxfunc pw est of t^(1-s)partial_t Fs}
	\sup_{(y,t)\, ;\,|y-x|<t} |t^{1-\sigma}\partial_t F^s(y,t)|\lesssim \mathcal{M}(\fL{\sigma}f)(x).
\end{equation}
For $\sigma\in[0,1]$, denoting $\nabla^\sigma = \nabla I^{1-\sigma} = \mathcal{R}\fL{\sigma}= (\mathcal{R}_i\fL{\sigma})_i$,
\begin{equation}\label{eq:maxfunc pw est of t^(1-s)nabla_x Fs}
	\sup_{(y,t)\, ;\,|y-x|<t} |t^{1-\sigma}\nabla_x F^s(y,t)|\lesssim \mathcal{M}(\nabla^\sigma f)(x).
\end{equation}
Finally, for any $0<\sigma<n$,
\begin{equation}\label{eq:maxfunc pw est of t^sigma Fs}
	\sup_{(y,t)\, ;\,|y-x|<t} t^\sigma|F^s(y,t)|\lesssim \mathcal{M}(I^\sigma f)(x).
\end{equation}
\end{proposition}

\begin{proof}[Estimate (\ref{eq:maxfunc pw est of Fs})]
Estimate (\ref{eq:maxfunc pw est of Fs}) is a direct consequence of Proposition \ref{prop:maxfunc - convolution estimate} since $P^s_t f = C_{n,s}(p^s_t*f) = C_{n,s} (p^s_1)_t*f$ and 
\begin{equation*}
	p^s_1(x) = \frac{1}{(1+|x|^2)^\frac{n+s}{2}}
\end{equation*}
is radial, non-increasing, bounded and integrable.
\end{proof}

\begin{proof}[Estimate (\ref{eq:maxfunc pw est of t^(1-s)partial_t Fs})]
Following the proof of Proposition \ref{prop:tt square functions to Bessel potential spaces}, more specifically (\ref{eq:tt fitting kernel for partial_t ps1 square function estimate}), we have
\begin{equation*}
	t^{1-\sigma}\partial_t F^s(x,t) = (q_t*\fL{\sigma}f)(x)
\end{equation*}
with
\begin{equation*}
	q(x) = C\fL{s-\sigma}\left((|x|^2+1)^\frac{s-2-n}{2}\right).
\end{equation*}
For $\sigma<s$, just as in the proof of Proposition \ref{prop:tt square functions to Bessel potential spaces}, we obtain
\begin{equation*}
	|q(x)| \leq C \frac{1}{(1+|x|)^{n+s-\sigma}}
\end{equation*}
with Lemma \ref{lem:fL Linfty estimate} and Lemma \ref{lem:fL decay}.
For $\sigma=s$ we have 
\begin{equation*}
	|q(x)| = C \frac{1}{(1+|x|)^{n+2-s}}.
\end{equation*}
Either way, we have a radial, non-increasing, bounded and integrable majorant.
Therefore, Proposition \ref{prop:maxfunc - convolution estimate} yields (\ref{eq:maxfunc pw est of t^(1-s)partial_t Fs}) because $\fL{\sigma}f\in L^1_{\text{loc}}(\Rn)$ since $f\in\CicRn$. 
\end{proof}

\begin{proof}[Estimate (\ref{eq:maxfunc pw est of t^(1-s)nabla_x Fs})]
From the proof of Proposition \ref{prop:tt square functions to Bessel potential spaces}, more specifically (\ref{eq:tt fitting kernel for partial_xi ps1 square function estimate}), we know that
\begin{equation*}
	t^{1-\sigma}\partial_{x_i} F^s(x,t) = (q_t*\mathcal{R}_i\fL{\sigma}f)(x)
\end{equation*}
where 
\begin{equation*}
	q(x) = C \fL{1-\sigma}\left((1+|x|^2)^{-\frac{n+s}{2}}\right).
\end{equation*}
For $\sigma<1$, combining Lemma \ref{lem:fL Linfty estimate} and Lemma \ref{lem:fL decay}, we have
\begin{equation*}
	q(x) \leq C(1+|x|^2)^{-n-(1-\sigma)},
\end{equation*}
which is a radial, non-increasing, bounded and integrable majorant.
For $\sigma = 1$, $q(x)$ itself is radial, non-increasing, bounded and integrable.
Therefore, with Proposition \ref{prop:maxfunc - convolution estimate} we have
\begin{equation*}
	\sup_{(y,t)\, ;\,|y-x|<t} |t^{1-\sigma}\partial_{x_i} F^s(y,t)| \lesssim \mathcal{M}(\mathcal{R}_i\fL{\sigma}f)(x).
\end{equation*}
Combining these equations for all $i\in\lbrace 1,\ldots,n\rbrace$ leads to (\ref{eq:maxfunc pw est of t^(1-s)nabla_x Fs}).
\end{proof}

\begin{proof}[Estimate (\ref{eq:maxfunc pw est of t^sigma Fs})]
Via the Fourier symbols of fractional Laplacian and Riesz potential we obtain
\begin{mathex}
	t^\sigma F^s(x,t)	&=& \int_{\Rn} t^\sigma p^s_t(x-y)\, f(y) \,dy\\ 
						&=&	\int_{\Rn} t^\sigma (2\pi\xi)^\sigma\Fourier(p^s_t)(\xi)\, (2\pi\xi)^{-\sigma}\Fourier(f)(\xi)\\
						&=& \int_{\Rn} t^\sigma \fL{\sigma}p^s_t(x-y)\, I^\sigma f(y) \,dy.
\end{mathex}
At the beginning of step 4 in the proof of \ref{th:tt triebel lizorkin besov char}, we saw that $\fL{\sigma} p^s_t = t^{-\sigma}(\fL{\sigma}p^s_1)_t$.
Therefore,
\begin{equation*}
	t^\sigma F^s(x,t) = (q_t*I^\sigma f)(x)
\end{equation*}
with 
\begin{equation*}
	q(x) =  \fL{\sigma} p^s_1 = \fL{\sigma} \frac{1}{(|x|^2+1)^\frac{n+s}{2}}.
\end{equation*}
As in the proof of the other estimates, combining Lemma \ref{lem:fL Linfty estimate} and Lemma \ref{lem:fL decay}, we obtain a radial, non-increasing, bounded and integrable majorant.
Thus, (\ref{eq:maxfunc pw est of t^(1-s)nabla_x Fs}) follows with Proposition \ref{prop:maxfunc - convolution estimate}.
\end{proof}

\begin{remark}\label{re:maxfunc extended range for parameters in pw estiamtes}
In the same way that (\ref{eq:maxfunc pw est of t^sigma Fs}) is the extension of (\ref{eq:maxfunc pw est of Fs}), the estimates (\ref{eq:maxfunc pw est of t^(1-s)partial_t Fs}) and (\ref{eq:maxfunc pw est of t^(1-s)nabla_x Fs}) can also be extended to the range of $-n<\sigma<0$ as
\begin{equation*}
	\sup_{(y,t)\, ;\,|y-x|<t} |t^{1-\sigma}\partial_t F^s(y,t)|\lesssim \mathcal{M}(I^{-\sigma}f)(x)
\end{equation*}
and
\begin{equation*}
	\sup_{(y,t)\, ;\,|y-x|<t} |t^{1-\sigma}\nabla_x F^s(y,t)|\lesssim \mathcal{M}((\mathcal{R}_i I^{-\sigma} f)_i)(x).
\end{equation*}
\end{remark}

\begin{remark}\label{re:maxfunc pw estimates for double derivatives}
Plugging in $\partial_{x_i}F^s(x,t) = P^s_t(\partial_{x_i}f)(x)$ for $F^s(x,t)$ and $\tilde{\sigma} = 1+\sigma$ in (\ref{eq:maxfunc pw est of t^(1-s)partial_t Fs}), we obtain
\begin{equation}\label{eq:maxfunc pw est of t^(2-s)partial_tnabla_x Fs}
	\sup_{(y,t)\, ;\,|y-x|<t} |t^{2-\tilde{\sigma}}\partial_t\nabla_x F^s(y,t)|\lesssim \mathcal{M}((\mathcal{R}_i\fL{\tilde{\sigma}}f)_i)(x)
\end{equation}  
for any $\tilde{\sigma}\in[0,1+s]$, with Remark \ref{re:maxfunc extended range for parameters in pw estiamtes} in mind.
Analogously, from (\ref{eq:maxfunc pw est of t^(1-s)nabla_x Fs}) we obtain
\begin{equation}\label{eq:maxfunc pw est of t^(2-s)nabla_xnabla_x Fs}
	\sup_{(y,t)\, ;\,|y-x|<t} |t^{2-\tilde{\sigma}}\nabla_x^2 F^s(y,t)|\lesssim \mathcal{M}((\mathcal{R}_i\mathcal{R}_j\fL{\tilde{\sigma}} f)_{i,j})(x)
\end{equation}
for any $\tilde{\sigma}\in[0,2]$.
\end{remark}

\subsection{Blackbox estimates from $\Rnpp$ to $\Rn$}\label{subsec:tt - blackbox estimates}
Having now collected a variety of estimates for the $s$-harmonic extension, we can combine these results with relatively little effort to more complex, custom estimates involving multiple $\CicRn$-functions.
We obtain three different estimates.
First, there is a Lorentz-estimate, where we obtain an estimate against the Lorentz-space norms of the occurring functions' fractional derivatives. 
Second, we have a $BMO$-estimate, where one $BMO$-norm appears next to the Lorentz space norms.
Analogous to this $BMO$-estimate, we also obtain a H\"older space estimate.

\begin{proposition}[Lorentz-estimate]\label{prop:tt - Lp est}
Let $s\in(0,2)$ and $2\leq k\in\N$.
Assume $p_i\in(1,\infty)$, $q_i\in[1,\infty]$ for $i\in\lbrace 1,\ldots k \rbrace$ with
\begin{equation*}
	\sum_{i=1}^k \frac{1}{p_i} = \sum_{i=1}^k \frac{1}{q_i} = 1.
\end{equation*}
Denote with $F^s_i(x,t)=P^s_t f_i(x)$ the extension of $f_i\in\CicRn$, $i\in\lbrace 1,\ldots k \rbrace$.
Then for $s_i\in(-n,0)$, $i= 1,2$, and $s_i\in(-n,0]$, $i\geq 3$, the following estimate holds.
\begin{mathex}
	\int_{\Rnpp} t^{-1-s_1-\ldots -s_k} |F^s_1(x,t)|\, |F^s_2(x,t)|\ldots |F^s_k(x,t)|\,dx\,dt\\
	\qquad\lesssim 
	\|I^{-s_1} f_1\|_{L^{(p_1,q_1)}}\, \|I^{-s_2} f_2\|_{L^{(p_2,q_2)}}\ldots \|I^{-s_k} f_k\|_{L^{(p_k,q_k)}}.
\end{mathex}
For any $i\in\lbrace 1,\ldots k \rbrace$, the following variations of the estimate hold.
Whenever choosing $s_i> 0$, replace $\|I^{-s_1} f_i\|_{L^{(p_i,q_i)}}$ with $\|\fL{s_1} f_i\|_{L^{(p_i,q_i)}}$.
\begin{enumerate}[label=(\alph*)]
\item Replacing $|F^s_i(x,t)|$ with $t|\partial_t F^s_i(x,t)|$ allows $s_i\in(-n,s)$ ($s_i\in(-n,s]$ for $i\geq 3$).
\item Replacing $|F^s_i(x,t)|$ with $t|\nabla_x F^s_i(x,t)|$ allows $s_i\in(-n,1)$ ($s_i\in(-n,1]$ for $i\geq 3$).
\item Replacing $|F^s_i(x,t)|$ with $t^2|\nabla_x\nabla_{\Rnp} F^s_i(x,t)|$ allows $s_i\in(-n,\min\lbrace 2,1+s\rbrace)$ for $i=1,2$ and $s_i\in(-n,\min\lbrace 2,1+s\rbrace]$ for $i\geq 3$.
\end{enumerate}
The estimate also holds if $(p_i,q_i)=(\infty,\infty)$ for any $i\geq 3$, except for variations (b) and (c).
\end{proposition}

\begin{proof}
There are two different options to treat the additional functions, i.e. the $F_i$-terms for $i\geq 3$.
First, ignoring the extended parameter range for $i\geq 3$, we may show the result via two Hölder estimates and the Bessel-potential space estimates from Subsection \ref{subsec:tt - building blocks}. 
\begin{mathex}[LCL]
	\mathcal{I} &\coloneqq& 
	\int_{\Rnpp} t^{-1-s_1-s_2-\ldots -s_k} |F^s_1(x,t)|\, |F^s_2(x,t)|\ldots |F^s_k(x,t)|\,dx\,dt\\
	&=& \int_{\Rn}\int_0^\infty t^{-\frac{1}{k}-s_1} |F^s_1(x,t)| \ldots t^{-\frac{1}{k}-s_k}|F^s_k(x,t)|\,dt\,dx\\
	&\leq&	\int_{\Rn}	\left(\int_0^\infty |t^{-\frac{1}{k}-s_1} F^s_1(x,t)^k\,dt \right)^\frac{1}{k}\ldots
						\left(\int_0^\infty |t^{-\frac{1}{k}-s_k} F^s_k(x,t)|^k\,dt \right)^\frac{1}{k} \,dx
\end{mathex}		
Applying Theorem \ref{th:tt lorentz holder}, the Hölder estimate for Lorentz spaces, we obtain
\begin{mathex}[LCL]
	\mathcal{I} &\lesssim&	
			\left\|x\mapsto \left(\int_0^\infty |t^{-\frac{1}{k}-s_1} F^s_1(x,t)|^k\,dt \right)^\frac{1}{k}\right\|_{L^{(p_1,q_1)}} 
			\ldots\left\|x\mapsto \left(\int_0^\infty |t^{-\frac{1}{k}-s_k} F^s_k(x,t)|^k\,dt \right)^\frac{1}{k}\right\|_{L^{(p_k,q_k)}}\\
	&\lesssim&	
			\|I^{-s_1} f_1\|_{L^{(p_1,q_1)}} \ldots \|I^{-s_k} f_k\|_{L^{(p_k,q_k)}}.
\end{mathex}
The last estimate is due to Remark \ref{re:tt bessel potential space estimates for larger parameters} and (\ref{eq:tt Bessel potential spaces via unmodified Fs}).
For variation (a), instead of (\ref{eq:tt Bessel potential spaces via unmodified Fs}) apply (\ref{eq:tt Bessel potential spaces via partial_tFs}), for variation (b) apply (\ref{eq:tt Bessel potential spaces via  nabla_xFs}) and for variation (c) apply the combination of (\ref{eq:tt Bessel potential spaces via  nabla_x^2Fs}) and (\ref{eq:tt Bessel potential spaces via partial_tnabla_xFs}). 

For the second option, we use the pointwise estimates from Proposition \ref{prop:maxfunc pointwise estimates} to treat the $F_i$-terms for $i\geq 3$, which makes the use of Remark \ref{re:tt bessel potential space estimates for larger parameters} obsolete.
For any $i\geq 3$  we apply (\ref{eq:maxfunc pw est of t^sigma Fs}) and  extract
\begin{equation*}
	\sup_{t>0} t^{-s_i}| F_i^s(x,t)| \leq \mathcal{M}(I^{-s_i}f_i)(x)
\end{equation*} 
from the inner integral before the first Hölder estimate.
Following the same line as above, we then have
\begin{equation*}
	\mathcal{I} \lesssim
	 \|I^{-s_1} f_1\|_{L^{(p_1,q_1)}} \|I^{-s_2} f_2\|_{L^{(p_2,q_2)}} 
	 \|\mathcal{M}(I^{-s_3}f_3)\|_{L^{(p_3,q_3)}}\ldots \|\mathcal{M}(I^{-s_k}f_k)\|_{L^{(p_k,q_k)}}.
\end{equation*}
We obtain the result with Theorem \ref{th:maxfunc - Lp estimate}. For variation (a), instead of (\ref{eq:maxfunc pw est of t^sigma Fs}) we use Remark \ref{re:maxfunc extended range for parameters in pw estiamtes} and (\ref{eq:maxfunc pw est of t^(1-s)partial_t Fs}). 
For variations (b) and (c), we use (\ref{eq:maxfunc pw est of t^(1-s)nabla_x Fs}) and a combination of (\ref{eq:maxfunc pw est of t^(2-s)partial_tnabla_x Fs}) and (\ref{eq:maxfunc pw est of t^(2-s)nabla_xnabla_x Fs}). 
The occurring Riesz transforms pose no problem due to their $L^p$-boundedness, see Proposition \ref{prop:Rt Lp bounded}.
The case $(p_i,q_i)=(\infty,\infty)$ is not compatible with the variations (b) and (c) because the $L^p$-boundedness of the Riesz transforms does not extend to $p=\infty$.
\end{proof}

In order to piece together the following BMO estimate, we need an additional estimate from \cite{Ste93}. This estimate will replace the $1$-dimensional Hölder estimates from the proofs of Proposition \ref{prop:tt - Lp est}, yielding terms related to Carleson-measures such as in the $BMO$-space characterizations and square functions such as in Proposition \ref{prop:tt square functions to Bessel potential spaces}. 
The result can be found in \cite[Chapter IV, §4.4, p.162, Proposition]{Ste93}.

\begin{lemma}\label{lem:tt carleson measure and square function estimate}
Let $\Phi,G$ be measurable functions given on $\Rnpp$. Then
\begin{mathex}
	\int_{\Rnpp} \Phi(x,t)G\,(x,t)\,dx\,dt\\
	\lesssim	
			\sup_{x\in\Rn,r>0} \left(\frac{1}{|B(x,r)|}\int_{T(B(x,r))}t |\Phi(y,t)|^2 dy\,dt\right)^\frac{1}{2}
			\int_\Rn\left(\int_{\lbrace (y,t);|y-x|<t\rbrace}|G(y,t)|^2 \frac{dtdy}{t^{n+1}}\right)^\frac{1}{2}dx
\end{mathex}
where $T(B(x,r))=\lbrace (y,t)\in\Rnpp\, :\,|y-x|<r-t\rbrace$ is the ``tent'' over $B(x,r)$.
\end{lemma}

Equipped with this additional tool, we can now tackle the $BMO$-estimates.

\begin{proposition}[$BMO$-estimate]\label{prop:tt - BMO est}
Let $s\in(0,2)$ and $2\leq k\in\N$.
Assume $p_i\in(1,\infty)$, $q_i\in[1,\infty]$ for $i\in\lbrace 1,\ldots,k\rbrace$ with
\begin{equation*}
	\sum_{i=1}^k \frac{1}{p_i} = \sum_{i=1}^k \frac{1}{q_i} = 1.
\end{equation*}
Denote with $\Phi^s(x,t)=P^s_t \phi(x)$, $F_i^s(x,t)=P^s_t f_i(x)$ for $\phi,f_i\in\CicRn$, $i\in\lbrace 1,\ldots,k\rbrace$.
Then for $s_1\in(-n,0)$ and $s_i\in(-n,0]$, $i\in\lbrace 2,\ldots,k\rbrace$, the following estimate holds.
\begin{mathex}
	\int_{\Rnpp} t^{-s_1-s_2-\ldots-s_k}|\nabla_{\Rnp}\Phi^s(x,t)|\, |F_1^s(x,t)|\, |F_2^s(x,t)|\ldots |F_k^s(x,t)| \,dx\,dt\\
	\qquad\lesssim [\phi]_{BMO} \|I^{-s_1}f_1\|_{L^{(p_1,q_1)}}\|I^{-s_2}f_2\|_{L^{(p_2,q_2)}}\ldots \|I^{-s_k}f_k\|_{L^{(p_k,q_k)}}.
\end{mathex}
The estimate also holds in the following variations.
\begin{enumerate}[label=(\alph*)]
	\item $|\nabla_{\Rnp}\Phi^s(x,t)|$ can be replaced with $|t^{\beta-1}P^s_t(\fL{\beta}\phi)(x)|$ for any $\beta>0$.
	\item $|\nabla_{\Rnp}\Phi^s(x,t)|$ can be replaced with $|t^{\beta}\nabla_{\Rnp}P^s_t(\fL{\beta}\phi)(x)|$ for any $\beta>0$.
\end{enumerate}
Additionally, for any $i\in\lbrace 1,\ldots k \rbrace$, the following variations of the estimate hold.
Whenever choosing $s_i\geq 0$, replace $\|I^{-s_1} f_i\|_{L^{(p_i,q_i)}}$ with $\|\fL{s_1} f_i\|_{L^{(p_i,q_i)}}$.
\begin{enumerate}[label=(\alph*)]
	\setcounter{enumi}{2}
	\item Replacing $|F_i^s(x,t)|$ with $t|\partial_t F_i^s(x,t)|$ allows $s_i\in(-n,s]$ ($s_1\in(-n,s)$ for $i=1$).
	\item Replacing $|F_i^s(x,t)|$ with $t|\nabla_x F_i^s(x,t)|$ allows $s_i\in(-n,1]$ ($s_1\in(-n,1)$ for $i=1$).
	\item Replacing $|F_i^s(x,t)|$ with $t^2|\nabla_x\nabla_{\Rnp}F_i^s(x,t)|$ allows $s_i\in(-n,\min\lbrace 1+s,2\rbrace]$ for $i\geq 2$, $s_1\in(-n,\min\lbrace 1+s,2\rbrace)$ for $i=1$.
\end{enumerate}
The estimate also holds if $(p_i,q_i)=(\infty,\infty)$ for any $i\geq 2$, except for variations (d) and (e).
\end{proposition}

\begin{proof}
To proof these estimates, we first apply Lemma \ref{lem:tt carleson measure and square function estimate}. 
The first term then corresponds exactly to the $BMO$-space characterizations.
The second term resembles the square function estimates from Subsection \ref{subsec:tt - building blocks}.
Before we apply these, similar to the proof of Proposition \ref{prop:tt - Lp est} we first extract the supremum of the $F_i$- terms, $i\geq 2$. \\
Denote $\Gamma(x)\coloneqq \lbrace (y,t);|y-x|<t\rbrace$ for $x\in\Rn$.
With Lemma \ref{lem:tt carleson measure and square function estimate} we have
\begin{mathex}[LCLL]
	\mathcal{I}&\coloneqq&
		\int_{\Rnpp} t^{-s_1-s_2-\ldots-s_k}|\nabla_{\Rnp}\Phi^s| |F_1^s(x,t)| |F_2^s(x,t)|\ldots |F_k^s(x,t)| \,dx\,dt\\
	&\lesssim& 
		\sup_{x\in\Rn,r>0} \left(\frac{1}{|B(x,r)|}\int_{T(B(x,r))}t |\nabla_{\Rnp}\Phi(y,t)|^2 dy\,dt\right)^\frac{1}{2}\\
	& &	\cdot
		\int_\Rn\left(\int_{\Gamma(x)}\left(t^{-s_1}|F_1^s(y,t)| t^{-s_2}|F_2^s(y,t)| \ldots t^{-s_k}|F_k^s(y,t)|\right)^2 \frac{dtdy}{t^{n+1}}\right)^\frac{1}{2}dx.
\end{mathex}
Due to (\ref{eq:maxfunc pw est of t^sigma Fs}) we have $\sup_{\Gamma(x)} t^{-s_i}|F_i^s(y,t)|\lesssim \mathcal{M}(I^{-s_i} f_i)(x)$ for all $i\in\lbrace 2,\ldots,k\rbrace$.  
Applying (\ref{eq:tt BMO characterization with nabla_Rnp}), we obtain
\begin{mathex}
	\mathcal{I}&\lesssim& 
		[\phi]_{BMO} \int_{\Rn} \mathcal{M}(I^{-s_2}f_2)(x)\ldots\mathcal{M}(I^{-s_k}f_k)(x)\left(\int_{\Gamma(x)}t^{-2s_1-n-1}|F_1^s(y,t)|^2 \,dt\,dy\right)^\frac{1}{2}dx\\
	&\lesssim&
		[\phi]_{BMO} \|\mathcal{M}(I^{s_2}f_2)\|_{L^{(p_2,q_2)}} \ldots \|\mathcal{M}(I^{s_k}f_k)\|_{L^{(p_k,q_k)}}
		\|I^{-s_1}f_1\|_{L^{(p_1,q_1)}}
\end{mathex}
where for the last estimate we applied Theorem \ref{th:tt lorentz holder}, the Hölder estimate for Lorentz spaces, and then immediately the square function estimate (\ref{eq:tt square function with unmodified kernel}).
We then obtain the result with Theorem \ref{th:maxfunc - Lp estimate}(ii).

The variations are obtained analogously by swapping out the respective estimates.
For variations (a) and (b), instead of (\ref{eq:tt BMO characterization with nabla_Rnp}) we use (\ref{eq:tt BMO characterization with fracLapl(beta)}) and (\ref{eq:tt BMO characterization with nabla_Rnpp fracLapl(beta)}) respectively.
For variations (c) or (d), instead of (\ref{eq:maxfunc pw est of t^sigma Fs}) and (\ref{eq:tt square function with unmodified kernel}) we use (\ref{eq:tt square function with partial_t}) and (\ref{eq:maxfunc pw est of t^(1-s)partial_t Fs}) or (\ref{eq:tt square function with nabla_x}) and (\ref{eq:maxfunc pw est of t^(1-s)nabla_x Fs}).
For variation (e), we use (\ref{eq:tt square function with double derivative}) and a combination of (\ref{eq:maxfunc pw est of t^(2-s)partial_tnabla_x Fs}) and (\ref{eq:maxfunc pw est of t^(2-s)nabla_xnabla_x Fs}).
Additionally, for variations (d) and (e), we need the $L^{(p,q)}$-boundedness of the Riesz transform, which we obtain via interpolation from Proposition \ref{prop:Rt Lp bounded}, the $L^p$-boundedness.

The annotation regarding $(p_i,q_i)=(\infty,\infty)$ is obvious. 
For variations (d) and (e) this does not work though, since the $L^p$-boundedness of the Riesz transforms only holds for $1<p<\infty$. 
\end{proof}

\begin{proposition}[Hölder-space-estimates]\label{prop:tt - Hölder est}
Let $s\in(0,2)$ and $2\leq k\in\N$.
Assume $p_i\in(1,\infty)$, $q_i\in[1,\infty]$ for $i\in\lbrace 1,\ldots k \rbrace$ with
\begin{equation*}
	\sum_{i=1}^k \frac{1}{p_i} = \sum_{i=1}^k \frac{1}{q_i} = 1.
\end{equation*}
Denote with $\Phi^s(x,t)=P^s_t \phi(x)$ and $F^s_i(x,t)=P^s_t f_i(x)$ the extensions of $ \phi,f_i\in\CicRn$, $i\in\lbrace 1,\ldots k \rbrace$.
Then for $\nu\in(0,\min\lbrace 1,s\rbrace)\setminus\N$, $s_i\in(-n,0)$ for $i=1,2$, and $s_i\in(-n,0]$ for $i\geq 3$, the following estimate holds.
\begin{mathex}
	\int_{\Rnpp} t^{-\nu-s_1-\ldots -s_k} |\nabla_{\Rnp}\Phi^s(x,t)|\, |F^s_1(x,t)|\ldots |F^s_k(x,t)|\,dx\,dt\\
	\qquad\lesssim 
	[\phi]_{\C^\nu} \,\|I^{-s_1} f_1\|_{L^{(p_1,q_1)}} \ldots \|I^{-s_k} f_k\|_{L^{(p_k,q_k)}}.
\end{mathex}
The estimate continues to hold under the following variations.
\begin{enumerate}[label=(\alph*)]
	\item For $\nu=0$ replace $[\phi]_{\C^\nu}$ with $\|\phi\|_{L^\infty}$.
	\item Replacing $|\nabla_{\Rnp} \Phi^s(x,t)|$ with $|\partial_t \Phi^s(x,t)|$ allows $\nu\in(0,s)\setminus\N$.
	\item Replacing $|\nabla_{\Rnp} \Phi^s(x,t)|$ with $|\nabla_x \Phi^s(x,t)|$ allows $\nu\in(0,1]$.
	\item Replacing $|\nabla_{\Rnp} \Phi^s(x,t)|$ with $t|\nabla_x\nabla_{\Rnp} \Phi^s(x,t)|$ allows $\nu\in(0,\min\lbrace 2,1+s\rbrace)\setminus\N$.
	\item Replacing $[\phi]_{\C^\nu}$ with $[\fL{\nu}\phi]_{BMO}$ allows to choose $s_2$ as the $s_i$ for $i\geq 3$.
\end{enumerate}
Additionally, the following replacements are possible for any $i\in\lbrace 1,\ldots k \rbrace$. 
Whenever choosing $s_i\geq 0$, replace $\|I^{-s_1} f_i\|_{L^{(p_i,q_i)}}$ with $\|\fL{s_1} f_i\|_{L^{(p_i,q_i)}}$.
\begin{enumerate}[label=(\alph*)]
\setcounter{enumi}{5}
\item Replacing $|F^s_i(x,t)|$ with $t|\partial_t F^s_i(x,t)|$ allows $s_i\in(-n,s)$ ($s_i\in(-n,s]$ for $i\geq 3$).
\item Replacing $|F^s_i(x,t)|$ with $t|\nabla_x F^s_i(x,t)|$ allows $s_i\in(-n,1)$ ($s_i\in(-n,1]$ for $i\geq 3$).
\item Replacing $|F^s_i(x,t)|$ with $t^2|\nabla_x\nabla_{\Rnp} F^s_i(x,t)|$ allows $s_i\in(-n,\min\lbrace 2,1+s\rbrace)$ for $i=1,2$ and $s_i\in(-n,\min\lbrace 2,1+s\rbrace]$ for $i\geq 3$.
\end{enumerate}
The estimate also holds if $(p_i,q_i)=(\infty,\infty)$ for any $i\geq 3$, except for variations (g) and (h).\\
\end{proposition}

\begin{proof}
Combining (\ref{eq:tt Hölder characterization with partial_t}) and (\ref{eq:tt Hölder characterization with nabla_x}), we obtain
\begin{mathex}
	\mathcal{I}&\coloneqq&
		\int_{\Rnpp} t^{1-\nu-1-s_1-\ldots -s_k} |\nabla_{\Rnp}\Phi^s(x,t)| |F^s_1(x,t)|\ldots |F^s_k(x,t)|\,dx\,dt\\
	&\lesssim& 
		[\phi]_{\C^\nu} \int_{\Rnpp} t^{-1-s_1-\ldots -s_k} |F^s_1(x,t)|\ldots |F^s_k(x,t)|\,dx\,dt.
\end{mathex}
The result then follows with Proposition \ref{prop:tt - Lp est}.
For variation (a), instead of (\ref{eq:tt Hölder characterization with partial_t}) and (\ref{eq:tt Hölder characterization with nabla_x}) apply (\ref{eq:maxfunc pw est of t^(1-s)partial_t Fs}) and (\ref{eq:maxfunc pw est of t^(1-s)nabla_x Fs}) followed up by Theorem \ref{th:maxfunc - Lp estimate}. 
For variation (b) apply only (\ref{eq:tt Hölder characterization with partial_t}), for variation (c) only (\ref{eq:tt Hölder characterization with nabla_x}) and for variation (d) apply (\ref{eq:tt Hölder characterization with nabla_Rnp nabla_x}).

Regarding variation (e), we follow the proof for Proposition \ref{prop:tt - BMO est}, but instead of the $BMO$-characterizations from Proposition \ref{prop:tt bmo char}, we choose a suitable estimate from Corollary \ref{cor:tt BMO(fLnu) characterization}.
\end{proof}

\newpage
\section{Outlook: Possible extensions of the method}\label{sec:outlook}

There seem to be two main options to further extend the method.

On the one hand, we may further modify the blackbox estimates, for example by implementing the Triebel-Lizorkin and Besov-Lipschitz characterizations in a more general form. 
However, there are some limitations. 
We would loose the extended parameter range for additional functions, which we originally obtained due to pointwise estimates, confer the proof of Proposition \ref{prop:tt - Lp est}.
Additionally, for the $BMO$-estimate we only used the Triebel-Litzorkin space characterization of $BMO$, which we might replace with a more general $F^s_{\infty,q}$ norm. 
If we wanted to replace other norms with more general equivalents though, we would need Triebel-Lizorkin space characterizations reminiscent of the nontangential square function estimates.
But since the regular square function estimates are special cases of the space characterizations which we obtained from the Bui-Candy result, it does not seem too unlikely that Triebel-Lizorkin spaces can be characterized by expressions similar to nontangential square functions.

On the other hand, we might modify the means of extending from $\Rn$ to $\Rnpp$.
For example, in Subsection \ref{sec:Rp comms}, we only obtained the Chanillo commutator estimate for Riesz potentials of order $<1$.
In general, the extension by Cafarelli and Silvestre in \cite{CaSi07} only works for fractional Laplacians of order $<2$. 
In order to prove this estimate for higher orders, we would also need a higer order extension, which has been introduced for example in \cite{RoSti16} or \cite{Yan13}.
It might also be possible to extend other operators with suitably chosen extensions. 
Lenzmann and Schikorra mention the Bessel potential operator $(1-\Delta)^\frac{s}{s}$ as possible example.
We have already seen this operator as the lifting operator for inhomogeneous Triebel-Lizorkin and Besov-Lipschitz spaces.

Besides all these possible modifications, we do not expect the applications, which we have shown here, to be exhaustive for the unmodified method.
There should be a variety of other estimates, which are possible to prove with this method. 
For example, Lenzmann and Schikorra suggest sharp limit space estimates as in \cite{BouLi14} and \cite{Li19} as particularly interesting possible applications.

\newpage
\appendix
\section{Appendix}\label{sec:jac}

\subsection{Stokes' Theorem for differential forms on $\Rnpp$}

\begin{theorem}[Stokes' Theorem for differential forms on $\Rnpp$]\label{th:Stoke Rnpp} 
Let $w$ be an $n$-form on $\Rnpp$ with
\begin{mathex}
	w = \sum_{1\leq i_1 \leq\ldots\leq i_n\leq n+1} w_{i_1,\ldots,i_n} ~ dx^{i_1}\wedge\ldots\wedge dx^{i_n}
\end{mathex}
where $w_{1,\ldots,n}$ is integrable on $\partial\Rnpp$ and all $w_{i_1,\ldots,i_n}$ as well as their derivatives are integrable on $\Rnpp$.
Then
\begin{mathex}
	\int_\Rnpp dw = \int_{\partial\Rnpp} w = \int_{\partial\Rnpp} w_{1,\ldots,n}.
\end{mathex} 
\end{theorem}

\begin{proof}
For compactly supported $w$ the statement is true. The according proof is well known since it is an essential part of the usual proof of Stoke's Theorem for compactly supported differential forms on arbitrary oriented smooth manifolds, see for example Theorem 16.11 in \cite{Lee12}.

For arbitrary $w$ we will approximate $w$ via compactly supported differential forms. 
We obtain these approximations by multiplying $w$ with a special compactly supported function.
Let $\tilde{\eta}\in\C^\infty([0,\infty))$ with 
\begin{mathex}
	\tilde{\eta}(r) = \begin{cases} 	1 & \text{for } 0\leq r \leq 1 \\
									0 & \text{for } 2\leq r \end{cases},
	\qquad	\tilde{\eta}(r) \in [0,1] \quad\text{for } 1\leq r \leq 2.
\end{mathex}
For $R>0$ we define $\eta_R\in\C_c^\infty(\Rnpp)$ by $\eta_R(z) \coloneqq \tilde{\eta}\left(\frac{|z|}{R}\right)$. 
Thus, we have 
\begin{mathex}
	\eta_R(z) = \begin{cases} 	1 & \text{for } 0\leq |z| \leq R \\
								0 & \text{for } 2R \leq |z|  \end{cases},
	\qquad	\eta_R(z) \in [0,1] \quad\text{for } R\leq |z| \leq 2R.
\end{mathex}
Regarding the derivatives of $\eta_R$, we obtain
\begin{equation*}
	\nabla_{\Rnp}\eta_R(z) = 0 \quad\text{for }0\leq |z| \leq R \text{ or } 2R \leq |z|,
	\qquad	\|\nabla_{\Rnp}\eta_R\|_{L^\infty} \leq \frac{\|\tilde{\eta}'\|_{L^\infty}}{R}.
\end{equation*}
Apparently, $(\eta_R w)_{R>0}$ are compactly supported differential forms that converge to $w$ pointwise.
Since $w_{1,\ldots,n}$ is integrable on $\partial\Rnpp$, by Lebesgue's dominated convergence theorem, 
\begin{equation*}
	\lim_{R\rightarrow\infty} \int_{\partial\Rnpp} \eta_R w = \lim_{R\rightarrow\infty} \int_{\partial\Rnpp} \eta_R w_{1,\ldots,n} 
	= \int_{\partial\Rnpp} w_{1,\ldots,n}.
\end{equation*}
Now turning to the term with the exterior derivative, we observe
\begin{equation*}
	\lim_{R\rightarrow\infty} \int_{\Rnpp} d(\eta_R w) 
	= \lim_{R\rightarrow\infty} \left(\int_{\Rnpp} d\eta_R \wedge w + \int_{\Rnpp} \eta_R\, dw \right).
\end{equation*}
Thanks to  the integrability of the $w_{i_1,\ldots,i_n}$ and the choice of $\eta_R$,
\begin{mathex}
	\left| \int_{\Rnpp} d\eta_R \wedge w \right| 
	&\leq \|\nabla_{\Rnp}\eta_R\|_{L^\infty} \sum_{1\leq i_1 \leq\ldots\leq i_n\leq n+1} 
													\int_{B_{2R}(0)\setminus B_R(0)} |w_{i_1,\ldots,i_n}| \\
	&\leq \frac{\|\tilde{\eta}'\|_{L^\infty}}{R} \sum_{1\leq i_1 \leq\ldots\leq i_n\leq n+1} \|w_{i_1,\ldots,i_n}\|_{L^1(\Rnpp)}
	\overset{R\rightarrow \infty}{\longrightarrow} 0.
\end{mathex}
Together with Lebesgue's dominated convergence theorem and the integrability of the derivatives of the $w_{i_1,\ldots,i_n}$ and therefore of $dw$, we obtain
\begin{equation*}
	\lim_{R\rightarrow\infty} \int_{\Rnpp} d(\eta_R w) = 0 + \lim_{R\rightarrow\infty} \int_{\Rnpp} \eta_R\, dw = \int_\Rnpp dw.
\end{equation*}
Since $\eta_R w$ are compactly supported differential forms for $R>0$,
\begin{equation*}
	\int_{\partial\Rnpp} w = \lim_{R\rightarrow\infty} \int_{\partial\Rnpp} \eta_R w 
	= \lim_{R\rightarrow\infty} \int_{\Rnpp} d(\eta_R w) = \int_\Rnpp dw.
\end{equation*}
\end{proof}

\newpage
\bibliographystyle{alpha}
\addcontentsline{toc}{section}{References}
\bibliography{LiteraturMA}

\end{document}